\documentclass[11pt,a4paper]{article}

\usepackage{authblk} 

\usepackage{mathrsfs}
\usepackage{latexsym,bm}
\usepackage{amsmath,amsfonts,amsmath,amssymb,amsthm}
\usepackage{CJKpunct,extarrows}

\usepackage{graphicx,subfigure,epstopdf,float}
\usepackage{enumerate,cases}
\usepackage[justification=centering]{caption}
\usepackage[margin=10pt,font=small,labelfont=bf,labelsep=period]{caption}

\usepackage{indentfirst}
\usepackage[top=25mm,bottom=20mm,left=25mm,right=20mm]{geometry}
\usepackage{cite}
\usepackage{listings}
\usepackage{booktabs}
\usepackage[bookmarksnumbered,colorlinks,linkcolor=blue,hyperindex]{hyperref}

\newtheorem{example}{Example}[section]

\usepackage{algorithm}
\usepackage{listings,xcolor} 
\usepackage{textcomp}
\usepackage[framed,numbered,autolinebreaks,useliterate]{mcode} 

\newcommand{\bb}{\boldsymbol}

\begin{document}

\title{\bf mVEM: A MATLAB Software Package for the Virtual Element Methods}
\author{Yue Yu\thanks{terenceyuyue@sjtu.edu.cn}}
\affil{School of Mathematical Sciences, Institute of Natural Sciences, MOE-LSC, Shanghai Jiao Tong University, Shanghai, 200240, P. R. China.}
\date{}
\maketitle

\begin{abstract}
  This paper summarizes the development of mVEM, a MATLAB software package containing efficient and easy-following codes for various virtual element methods (VEMs) published in the literature. We explain in detail the numerical implementation of the mixed VEMs for the Darcy problem and the three-dimensional linear VEMs for the Poisson equation. For other model problems, we present the construction of the discrete methods and only provide the implementation of the elliptic projection matrices. Some mesh related functions are also given in the package, including the mesh generation and refinement in two or three dimensions. mVEM is free and open source software.
\end{abstract}



\section{Introduction}

Developing the mimetic finite difference methods, Beir\~{a}o and Brezzi et al. proposed the virtual element method (VEM) in 2013, and established an abstract framework for error analysis \cite{Beirao-Brezzi-Cangiani-2013}. This method was further studied in \cite{Ahmad-Alsaedi-Brezzi-2013} for an extension to reaction-diffusion problems, where a crucial enhancement technique is introduced to construct a computable $L^2$ projection. The computer implementation has been further studied in \cite{Beirao-Brezzi-Marini-2014}.
The proposed finite dimensional space in \cite{Beirao-Brezzi-Cangiani-2013} has become the standard space for constructing conforming virtual element methods for second-order elliptic problems on polygonal meshes. The word \emph{virtual} comes from the fact that no explicit knowledge of the basis functions is necessary since the shape functions are piecewise continuous polynomials on the boundary of the element and are extended to the interior by assuming basis functions as solutions of local Laplace equations. The construction of the VEMs for elliptic problems is very natural and standard, which can be derived based on an integration by parts formula for the underlying differential operator. As a matter of fact, this idea is used to devise conforming and nonconforming VEMs for arbitrary order elliptic problems though the resulting formulation and theoretical analysis are rather involved \cite{Huang-2020,Chen-Huang-2020}.

VEMs have some advantages over standard finite element methods. For example, they are more convenient to handle partial differential equations on complex geometric domains or the ones associated with high-regularity admissible spaces. Until now, they have been successfully applied to solve various mathematical physical problems, such as the conforming and nonconforming VEMs for second-order elliptic equations \cite{Ahmad-Alsaedi-Brezzi-2013,Beirao-Brezzi-Cangiani-2013,Ayuso-Lipnikov-Manzini-2016,Cangiani-Manzini-Sutton-2017} and fourth-order elliptic equations \cite{Brezzi-Marini-2013,Chinosi-Marini-2016,Antonietti-Manzini-Verani-2018,Zhao-Zhang-Chen-2018}, the time-dependent problems
\cite{Vacca-Beirao-2015,Vacca-2017,Park-Chi-Paulino-2020,Adak-Nataraj-2020,Huang-Lin-2020}, the mixed formulation of the Darcy and Stokes problems \cite{Brezzi-Falk-Marini-2014,Beirao-Lovadina-Vacca-2017,Zhao-Zhang-Mao-2019} and the variational inequalities and hemivariational inequalities associated with the frictional contact problems \cite{Wang-Wei-2018,Wang-Zhao-2021,Wu-Wang-Han-2022,Wang-Wu-Han-2021}.

For second-order problems with variable coefficients and convection terms, direct use of the elliptic projection approximation of the gradient operator does not ensure the optimal convergence, so the external projection approximation has been introduced in the literature, see \cite{Beirao-Brezzi-Marini-Russo-2016d,Beirao-Brezzi-Marini-Russo-2016c} for example. This approximation technique is also commonly used for the construction of virtual element methods for complex problems such as elastic or inelastic mechanics
\cite{Beirao-Brezzi-Marini-2013,Artioli-Beirao-Lovadina-2017a,Artioli-Beirao-Lovadina-2017b}.
Since the virtual element space contains at least $k$-th order polynomials, the number of the degrees of freedom (d.o.f.s) in the virtual element space is generally higher than that in the classical finite element space when the polygonal element is degenerated into a triangle. These extra d.o.f.s are usually caused by internal moments and can be further reduced by exploiting the idea of building an incomplete finite element or the serendipity finite element \cite{BrennerScott2008}. In fact, Beir\~{a}o et al. has proposed the serendipity nodal VEM spaces in \cite{Beirao-Brezzi-Marini-Russo-2016b}.

In this paper, we are intended to develop a MATLAB software package for the VEMs in two or three dimensions, containing efficient and easy-following codes for various VEMs published in the literature. In particular, \cite{Beirao-Brezzi-Marini-2014} provided a detailed explanation of the formulation of the terms in the matrix equations for the high order virtual element method applied to such a problem in two dimensions, and \cite{Sutton-2017} presented a transparent MATLAB implementation of the conforming linear virtual element method for the Poisson equation in two dimensions. The construction of the VEMs for three-dimensional problems has been accomplished in many papers \cite{Gain-Talischi-Paulino-2014,Gain-Paulino-Leonardo-2015,Chi-Beirao-Paulino-2017,Beirao-Dassi-Russo-2017,
Beirao-Brezzi-Dassi-2018,Beirao-Dassi-Vacca-2020,Cihan-Hudobivnik-Aldakheel-2021}. However, to the best of knowledge, no related implementation is publicly available in the literature. As an extension of \cite{Sutton-2017} to three spatial dimensions, we have provided a clear and useable MATLAB implementation of the method for three-dimensional linear VEMs for the Poisson equation on general polyhedral meshes in the package with the detailed implementation given in Section \ref{sec:3D}. Although the current procedure is only for first-order virtual element spaces, the design idea can be directly generalized to higher-order cases.

The paper is organized as follows. In Section \ref{sec:Darcy}, we provide the complete and detailed implementation of the mixed VEMs for the Darcy problem as an example, which includes almost all of the programming techniques in mVEM, such as the construction of the data structure for polygonal meshes, the computation of elliptic and $L^2$ projection matrices, the treatment of boundary conditions, and examples that demonstrate the usage of running the codes, showing the solutions and meshes, computing the discrete errors and displaying the convergence rates. Section \ref{sec:mesh} summarizes the mesh related built-in functions, including the modified version of PolyMesher introduced in \cite{Talischi-Paulino-Pereira-2012} and generation of some special meshes for VEM tests. We also provide the generation of polygonal meshes by establishing the dual mesh of a Delaunay triangulation and some basic functions to show the polygonal meshes as well as a boundary setting function to identify the
Neumann and Dirichlet boundaries. Sections \ref{sec:Poisson}-\ref{sec:VI} discuss the construction of virtual element methods for various model problems and the computation of the corresponding elliptical projections, including the conforming and nonconforming VEMs for the Poisson equation or reaction-diffusion problems, the (locking-free) VEMs for the linear elasticity problems in the displacement/tensor type,  three VEMs for the fourth-order plate bending problems, the divergence-free mixed VEMs for the Stokes problems, the adaptive VEMs for the Poisson equation and the variational inequalities for the simplified friction problem. The paper ends with some concluding remarks in Section \ref{sec:conclude}.

\section{Mixed VEMs for the Darcy problem} \label{sec:Darcy}

Considering that the implementation of conforming or nonconforming VEMs has been publicly available in the literature \cite{Sutton-2017,Ortiz-Alvarez-Hitschfeld-2019}, we in this paper present the detailed implementation of the mixed VEMs for the Darcy problem to fill the gap in this regard.
The mixed VEM is first proposed in \cite{Brezzi-Falk-Marini-2014} for the classical model problem of Darcy flow in a porous medium. Let's briefly review the method in this section.

\subsection{Construction of the mixed VEMs}

\subsubsection{The model problem}

Given the (polygonal) computational domain $\Omega \subset \mathbb{R}^2$, let $f\in L^2(\Omega)$ and $g\in H^{1/2}(\partial \Omega)$. The Darcy problem is to find $p\in H^1(\Omega)$ such that
\begin{equation}\label{DarcyOriginal}
\begin{cases}
-{\rm div}(\mathbb{K}\nabla p) = f \quad & \mbox{in}~~\Omega, \\
(\mathbb{K}\nabla p)\cdot \boldsymbol{n} = g \quad & \mbox{on}~~\partial \Omega,
\end{cases}
\end{equation}
where $\mathbb{K}$ is a symmetric and positive definite tensor of size $2\times 2$. For simplicity, we assume that $\mathbb{K}$ is constant. The given data $f$ and $g$ satisfy the compatibility condition
\[\int_\Omega f {\rm d}x = \int_{\partial \Omega} g {\rm d}s.\]
To remove an additive constant, we additionally require that
\begin{equation}\label{pconstraint}
\int_\Omega p {\rm d}x = 0.
\end{equation}

Introducing the velocity variable $\boldsymbol{u} = \mathbb{K}\nabla p$, the above problem can be rewritten in the mixed form
\[
\begin{cases}
\boldsymbol{u} = \mathbb{K}\nabla p \quad & \mbox{in}~~\Omega, \\
{\rm div}\boldsymbol{u} = -f \quad & \mbox{in}~~\Omega, \\
\boldsymbol{u}\cdot\boldsymbol{n} = g \quad & \mbox{on}~~\partial \Omega.
\end{cases}
\]
Define
\[V_g = \{\boldsymbol{u}\in H({\rm div};\Omega): \boldsymbol{u}\cdot\boldsymbol{n} = g ~~\mbox{on}~~\partial \Omega\}, \qquad Q = L_0^2(\Omega), \qquad V = V_0.\]
The corresponding mixed variational problem is: Find $(\boldsymbol{u}, p)\in V_g\times Q$ such that
\begin{equation}\label{DarcyPro}
\begin{cases}
a(\bb{u},\bb{v}) + b(\bb{v}, p) & = 0, \quad  \bb{v}\in V,   \\
b(\bb{u},q) & = -(f,q), \quad   q\in Q,
\end{cases}
\end{equation}
where
\[a(\bb{u},\bb{v}) = (\mathbb{K}^{-1}\bb{u}, \bb{v}), \qquad b(\bb{v}, q) = ({\rm div}\bb{v}, q).\]
Note that the boundary condition is now related to the variable $\bb{u}$.

\subsubsection{The virtual element space and the elliptic projection}

The local virtual element space of $V$ is
\begin{align}
  V_k(K)
  & = \{\bb{v}\in H({\rm div};K)\cap H({\rm rot};K): \bb{v}\cdot\bb{n}|_e\in\mathbb{P}_k(e), \nonumber \\
  &  \hspace{3cm} {\rm div}\bb{v}|_K\in \mathbb{P}_{k-1}(K),~~ {\rm rot}\bb{v}|_K\in \mathbb{P}_{k-1}(K)\}, \label{DarcyVES}
\end{align}
where
\[{\rm div}\bb{v} = \partial_1v_1+\partial_2v_2, \quad {\rm rot}\bb{v} = \partial_1v_2-\partial_2v_1,\quad \bb{v}=(v_1,v_2)^T.\]

To present the degrees of freedom (d.o.f.s), we introduce a scaled monomial $\mathbb{M}_r(D)$ on a $d$-dimensional domain $D$
\[
\mathbb  M_{r} (D):= \Big \{ \Big ( \frac{\boldsymbol x -  \boldsymbol x_D}{h_D}\Big )^{\boldsymbol  s}, \quad |\boldsymbol  s|\le r\Big \},
\]
where $h_D$ is the diameter of $D$, $\boldsymbol  x_D$ the centroid of $D$, and $r$ a non-negative integer. For the multi-index ${\boldsymbol{s}} \in {\mathbb{N}^d}$, we follow the usual notation
\[\boldsymbol{x}^{\boldsymbol{s}} = x_1^{s_1} \cdots x_d^{s_d},\quad |\boldsymbol{s}| = s_1 +  \cdots  + s_d.\]
Conventionally, $\mathbb  M_r (D) =\{0\}$ for $r\le -1$.

The d.o.f.s can be given by
\begin{align}
&\int_e \bb{v}\cdot\bb{n} q{\rm d}s, \quad q\in \mathbb{M}_k(e),~~e\subset\partial K, \label{divdof1}\\
&\int_K \bb{v}\cdot\nabla q{\rm d}x, \quad q\in \mathbb{M}_{k-1}(K)\backslash\{1\}, \label{divdof2}\\
&\int_K {\rm rot}\bb{v}q{\rm d}x, \quad q\in \mathbb{M}_{k-1}(K) \label{divdof3}.
\end{align}
We remark that the moments on edges or elements are not divided by $|e|$ or $|K|$ since $a^K(\cdot,\cdot)$ is associated with the $L^2$ norm rather than the $H^1$ semi-norm for the case of Poisson equation.

For the Poisson equation, the elliptic projection maps from the virtual element space $V_k(K)$ into the polynomial space $\mathbb{P}_k(K)$ \cite{Beirao-Brezzi-Cangiani-2013,Ahmad-Alsaedi-Brezzi-2013,Beirao-Brezzi-Marini-2014}. For convenience, the image space is referred to as the elliptic projection space. For the Darcy problem, however, the elliptic projection space is now replaced by
\[\widehat{V}_k(K) = \{\widehat{\bb{v}}\in V_k(K): \widehat{\bb{v}} = \mathbb{K}\nabla \widehat{q}_{k+1}~~\mbox{for some}~~\widehat{q}_{k+1}\in \mathbb{P}_{k+1}(K)\}.\]
The elliptic projector $\widehat{\Pi}^K: V_k(K)\to \widehat{V}_k(K)$, $\bb{v}\mapsto\widehat{\Pi}^K\bb{v}$ is then defined by
\begin{equation}\label{ellipticProjDarcy}
a^K(\widehat{\Pi}^K\bb{v}, \widehat{\bb{w}}) = a^K(\bb{v}, \widehat{\bb{w}}), \quad \widehat{\bb{w}}\in \widehat{V}_k(K).
\end{equation}

We now consider the computability of the elliptic projection. In view of the symmetry of $\mathbb{K}$, the integration by parts gives
\begin{align}
a^K(\bb{v}, \widehat{\bb{w}})
& = \int_K\mathbb{K}^{-1}\bb{v}\cdot\widehat{\bb{w}}{\rm d}x = \int_K\bb{v}\cdot (\mathbb{K}^{-1}\widehat{\bb{w}}){\rm d}x \nonumber\\
& = \int_K\bb{v}\cdot \nabla \widehat{q}_{k+1} {\rm d}x = -\int_K \widehat{q}_{k+1}{\rm div}\bb{v} {\rm d}x + \int_{\partial K}\widehat{q}_{k+1}\bb{v}\cdot\bb{n}{\rm d}s. \label{integrationdarcy}
\end{align}
\begin{itemize}
  \item For the second term, since $\bb{v}\cdot\bb{n}|_e\in \mathbb{P}_k(e)$ we expand it in the scaled monomials on $e$ as
\[\bb{v}\cdot\bb{n}|_e(s) = c_1m_1^e(s) + \cdots + c_nm_n^e(s).\]
Clearly, the coefficients are uniquely determined by the d.o.f.s in \eqref{divdof1}.
  \item For the first term, noting that ${\rm div}\bb{v}|_K\in \mathbb{P}_{k-1}(K)$, we have
\[{\rm div}\bb{v}|_K(x) = c_1m_1^K(x) + \cdots + c_nm_n^K(x).\]
Let $q = m_i^K(x)$ and take inner product on both sides with respect to $q$. We obtain
\[\int_K {\rm div}\bb{v} q {\rm d}x = -\int_K \bb{v}\cdot\nabla q{\rm d}x + \int_{\partial K}\bb{v}\cdot\bb{n}q{\rm d}s, \quad q = m_i^K(x)\in \mathbb{P}_{k-1}(K).\]
Hence the coefficients are determined by the d.o.f.s in \eqref{divdof1} and \eqref{divdof2}.
\end{itemize}

In this paper, we only consider the implementation of the lowest order case $k=1$. At this time, only the d.o.f.s of the first and third types exist. The local d.o.f.s will be arranged as
\begin{align}
&\chi_i(\bb{v}) = \int_{e_i} (\bb{v}\cdot\bb{n}){\rm d}s, \quad i = 1,\cdots, N_v, \label{dof1} \\
&\chi_{N_v+i}(\bb{v}) = \int_{e_i} (\bb{v}\cdot\bb{n})\frac{s-s_{e_i}}{h_{e_i}}{\rm d}s, \quad i = 1,\cdots, N_v, \label{dof2}\\
&\chi_{2N_v+1}(\bb{v}) = \int_K {\rm rot}\bb{v}{\rm d}x, \label{dof3}
\end{align}
where $N_v$ is the number of the vertices of $K$, $h_e$ is the length of $e$, and $s$ is the natural parameter of $e$ with $s_e$ being the midpoint in the parametrization.

\subsubsection{The discrete problem}

In what follows, we denote the global virtual element space by $V_h$ associated with $V_k(K)$, and $Q_h$ by the discretization of $Q$, given as
\[Q_h:= \{q\in Q: q|_K\in \mathbb{P}_{k-1}(K),~~K\in \mathcal{T}_h.\}\]
In particular, $Q_h$ is piecewise constant for $k=1$.

The VEM approximation of $a^K(\bb{u},\bb{v})$ is
\[a_h^K(\bb{u},\bb{v}) = a^K(\widehat{\Pi}^K\bb{u},\widehat{\Pi}^K\bb{v}) + \|\mathbb{K}^{-1}\|S^K(\bb{u}-\widehat{\Pi}^K\bb{u},\bb{v}-\widehat{\Pi}^K\bb{v}),\]
where $\|\cdot\|$ is the Frobenius norm and the stabilization term is
\[S^K(\bb{v},\bb{w}) = \sum\limits_{i=1}^{2N_v+1} \chi_i(\bb{v}) \chi_i(\bb{w}). \]

The discrete mixed variational problem is: Find $(\boldsymbol{u}_h, p_h)\in V_h^g\times Q_h$ such that
\begin{equation}\label{discretemixedDarcy}
\begin{cases}
a_h(\bb{u}_h,\bb{v}_h) + b(\bb{v}_h, p_h) & = 0, \quad \bb{v}_h\in V_h, \\
b(\bb{u}_h,q_h) & = -(f,q_h), \quad q_h\in Q_h.
\end{cases}
\end{equation}
The constraint \eqref{pconstraint} is not naturally imposed in the above system. To this end, we introduce a Lagrange multiplier and consider the augmented variational formulation: Find $((\boldsymbol{u}_h, p_h),\lambda)\in V_h^g\times Q_h \times \mathbb{R}$ such that
\begin{equation}\label{augdiscretemixedDarcy}
\begin{cases}
a_h(\bb{u}_h,\bb{v}_h) + b(\bb{v}_h, p_h) & = 0, \quad \bb{v}_h\in V_h, \\
\displaystyle  b(\bb{u}_h,q_h) + \lambda \int_\Omega q_h{\rm d}x & = -(f,q_h), \quad q_h\in Q_h,\\
\displaystyle \mu\int_\Omega p_h{\rm d}x & = 0, \quad \mu\in \mathbb{R}.
\end{cases}
\end{equation}

It should be pointed out that the virtual element space $V_h$ cannot be understood as a vector or tensor-product space, which is different from the conforming or nonconforming VEMs for the linear elasticity problem. In the computation, $\bb{u}_h$ should be viewed as a scalar at this time.
Let $\bb{\varphi}_i$, $i=1,\cdots, N$ be the nodal basis functions of $V_h$, where $N$ is the dimension of $V_h$. We can write \[\bb{u} = \sum\limits_{i=1}^N \chi_i(\bb{u})\bb{\varphi}_i =: \bb{\varphi}^T\bb{\chi}(\bb{u}).\]
The basis functions of $Q_h$ are denoted by $\psi_l$, $l=1,\cdots,M$:
\[p_h = \sum\limits_{l=1}^Mp_l\psi_l.\]
Plug above equations in \eqref{discretemixedDarcy}, and take $\bb{v}_h = \bb{\varphi}_j$ and $q_h = \psi_l$. We have
\[
\begin{cases}
\sum\limits_{i=1}^Na_h(\bb{\varphi}_i,\bb{\varphi}_j)\chi_i + \sum\limits_{l=1}^M b(\bb{\varphi}_j, \psi_l)p_l & = 0, \quad j=1,\cdots,N, \\
\displaystyle \sum\limits_{i=1}^Nb(\bb{\varphi}_i,\psi_l)\chi_i + \lambda \int_\Omega \psi_l {\rm d}x & = -(f,\psi_l), \quad l=1,\cdots,M,\\
\displaystyle \sum\limits_{l=1}^M \int_\Omega \psi_l {\rm d}x p_l & = 0.
\end{cases}
\]
Let
\begin{equation}\label{dl}
d_l = \int_\Omega \psi_l {\rm d}x, \quad \bb{d} = [d_1,\cdots, d_M]^T.
\end{equation}
The linear system can be written in matrix form as
\begin{equation}\label{linearsystemDarcy}
\begin{bmatrix}
A   & B & \bb{0}\\
B^T & O & \bb{d}\\
\bb{0}^T & \bb{d}^T  & 0
\end{bmatrix}
\begin{bmatrix}
\bb{\chi} \\
\bb{p} \\
\lambda
\end{bmatrix}
=\begin{bmatrix}
\bb{0} \\
\bb{f} \\
0
\end{bmatrix},
\end{equation}
where
\[A=(a_h(\bb{\varphi}_j,\bb{\varphi}_i))_{N\times N}, \quad
B = (b(\bb{\varphi}_j, \psi_l))_{N\times M},\quad
\bb{f} = (-(f,\psi_l))_{M\times 1}.\]
For the global d.o.f.s in the unknown vector $\bb{\chi}$ in \eqref{linearsystemDarcy}, we shall arrange the first type in \eqref{dof1}, followed by the second and the third ones in \eqref{dof2} and \eqref{dof3}.

\subsection{Implementation}

\subsubsection{Overview of the code} \label{subsec:overview}

We first provide an overview of the test script.
\vspace{-0.8cm}
\begin{lstlisting}
%% Parameters
nameV = [32, 64, 128, 256, 512];
maxIt = length(nameV);
h = zeros(maxIt,1);   N = zeros(maxIt,1);
ErruL2 = zeros(maxIt,1);
ErrpL2 = zeros(maxIt,1);
ErrI = zeros(maxIt,1);

%% PDE data
pde = Darcydata;

%% Virtual element method
for k = 1:maxIt
    % load mesh
    fprintf('Mesh %d: \n', k);
    load( ['meshdata', num2str(nameV(k)), '.mat'] );
    % get boundary information
    bdStruct = setboundary(node,elem);
    % solve the problem
    [uh,ph,info] = Darcy_mixedVEM(node,elem,pde,bdStruct);
    % record and plot
    N(k) = length(uh);  h(k) = 1/sqrt(size(elem,1));
	[uhI,phI,nodeI,elemI] = ProjectionDarcy(node,elem,uh,ph,info,pde);
    figure(1);
    showresult(nodeI,elemI,pde.uexact,uhI);
    %showresult(nodeI,elemI,pde.pexact,phI);
    drawnow; %pause(0.1);
    % compute errors in discrete L2 norm
    [ErruL2(k),ErrpL2(k)] = getL2error_Darcy(node,elem,uh,ph,info,pde);
end

%% Plot convergence rates and display error table
figure,
showrateh(h,ErruL2,'r-*','||u-u_h||',  ErrpL2, 'b-s','||p-p_h||')

fprintf('\n');
disp('Table: Error')
colname = {'#Dof','h','||u-u_h||','||p-p_h||'};
disptable(colname,N,[],h,'%0.3e',ErruL2,'%0.5e',ErrpL2,'%0.5e');
\end{lstlisting}

In the \mcode{for} loop, we first load or generate the mesh data, which immediately returns the matrix \mcode{node} and the cell array \mcode{elem} defined later to the MATLAB workspace. Then we set up the boundary conditions to get the structural information of the boundary edges. The subroutine \mcode{Darcy\_mixedVEM.m} is the function file containing all source code to implement the VEM. When obtaining the numerical solutions, we can visualize the piecewise elliptic projection $\widehat{\Pi}^K \bb{u}_h$ by using the subroutines \mcode{ProjectionDarcy.m} and \mcode{showresult.m}. We then calculate
the discrete $L^2$ error defined as
\begin{equation} \label{ErrL2Darcy}
{\rm ErrL2} = \left( \sum\limits_{K \in \mathcal{T}_h} \| \bb{u} - \widehat{\Pi}^K \bb{u}_h \|_{0,E}^2 \right)^{1/2},
\end{equation}
through the subroutine \mcode{getL2error\_Darcy.m}. The procedure is completed by verifying the rate of convergence through \mcode{showrateh.m}.

The overall structure of a virtual element method implementation will be much the same as for a standard finite element method, as
outlined in Algorithm \ref{alg:vem}.
\begin{algorithm}[!htb]
  \caption{An overall structure of the implementation of a virtual element method \label{alg:vem}}
\textbf{Input}: Mesh data and PDE data
  \begin{enumerate}
    \item Get auxiliary data of the mesh, including some data structures and geometric quantities;
    \item Derive elliptic projections;
    \item Compute and assemble the linear system by looping over the elements;
    \item Apply the boundary conditions;
    \item Set solver and store information for computing errors.
  \end{enumerate}
\textbf{Output}: The numerical DoFs
\end{algorithm}

\subsubsection{Data structure} \label{subsec:datastructure}

\begin{figure}[!htb]
  \centering
  \subfigure[mesh]{\includegraphics[height=6cm,trim = 40 0 40 0,clip]{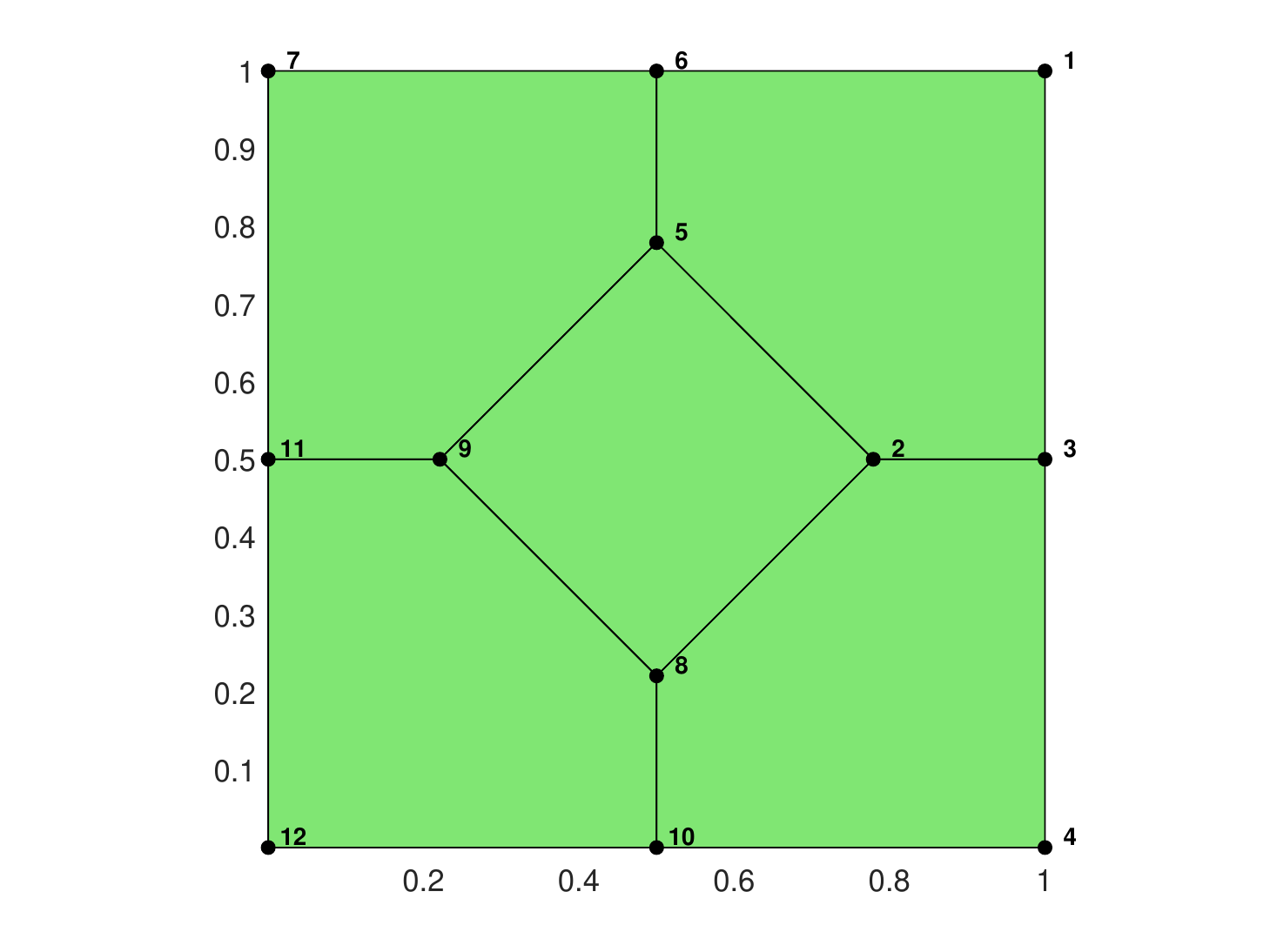}}
  \subfigure[node]{\includegraphics[height=5cm]{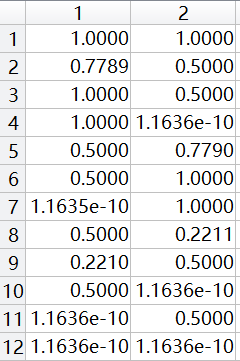}} \hspace{2em}
  \subfigure[elem]{\includegraphics[height=5cm,width=3cm]{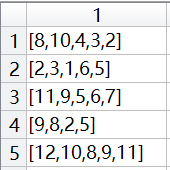}}\\
  \caption{Example of a polygonal mesh}\label{fig:polyh}
\end{figure}

We first discuss the data structure to represent polygonal meshes so as to facilitate the implementation. There are two basic data structures \mcode{node} and \mcode{elem}, where \mcode{node} is a matrix with the first and second columns contain $x$- and $y$-coordinates of the nodes in the mesh, and \mcode{elem} is a cell array recording the vertex indices of each element in a counterclockwise order as shown in Fig.~\ref{fig:polyh}. The mesh can be displayed by using \mcode{showmesh.m}.

Using the basic data structures, we can extract the topological or combinatorial structure of a polygonal mesh. These data are referred to as the auxiliary data structures as given in Tab.~\ref{auxsts}. The combinatorial structure will benefit the implementation of virtual element methods. The idea stems from the treatment of triangulation in $i$FEM for finite element methods \cite{ChenL-iFEM-2009}, which is generalized to polygonal meshes with certain modifications.

\begin{minipage}{0.47\textwidth}
\begin{table}[H]
  \centering
  \caption{Auxiliary data structures}\label{auxsts}
  \begin{tabular}{ll}
    \hline
    \mcode{edge}         \\
    \mcode{elem2edge}    \\
    \mcode{bdEdge}     \\
    \mcode{edge2elem}    \\
    \mcode{neighbor}    \\
    \mcode{node2elem}   \\
    \hline
  \end{tabular}
\end{table}
\end{minipage}
\begin{minipage}{0.47\textwidth}
\centering
\begin{figure}[H]
  \centering
  \includegraphics[scale=0.5]{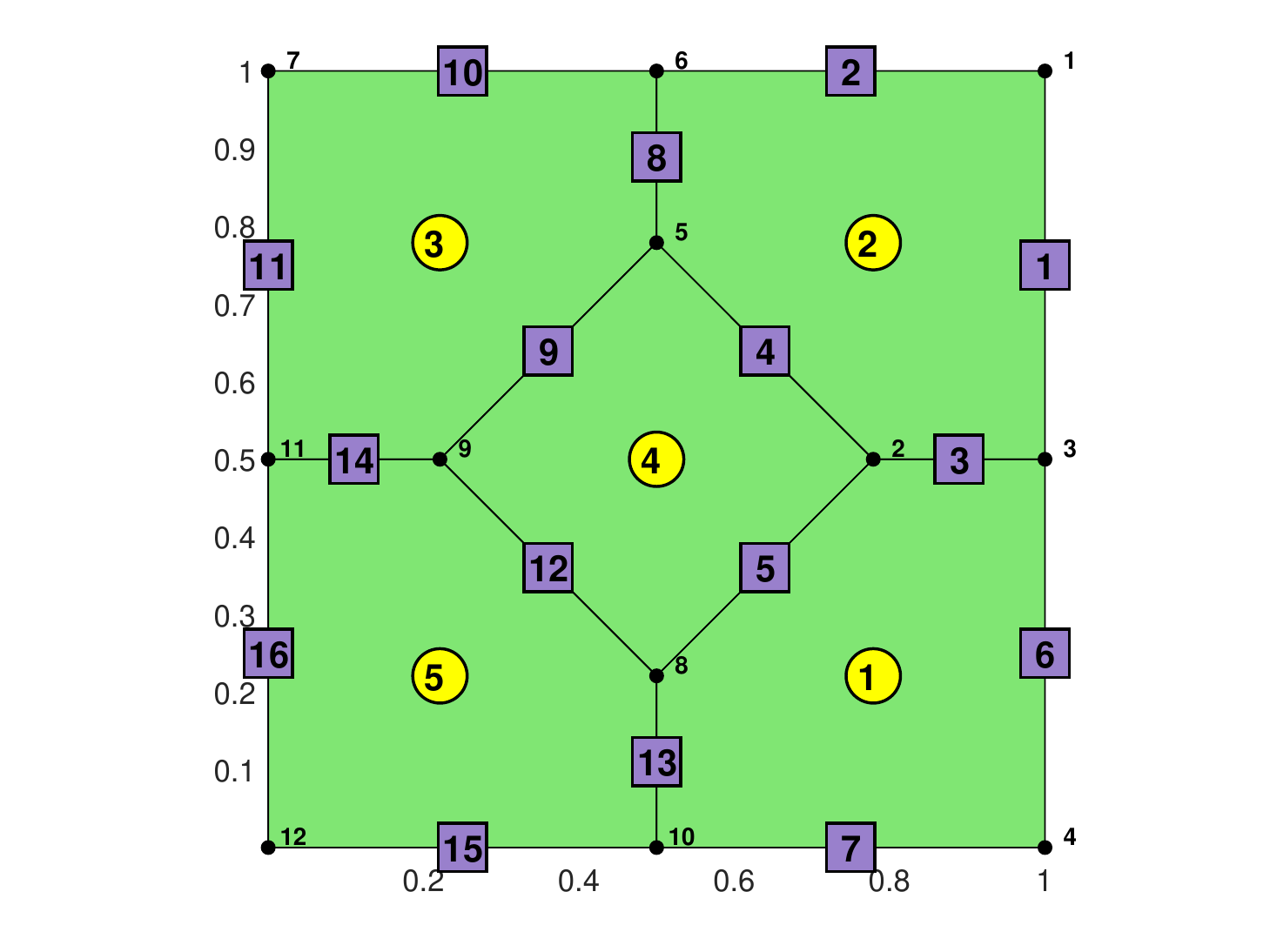}\\
  \caption{Illustration of the auxiliary data structures.}\label{Fig:auxstructure}
\end{figure}
\end{minipage}

\textbf{edge}. The d.o.f.s in \eqref{dof1} and \eqref{dof2} are edge-oriented. In the matrix \mcode{edge(1:NE,1:2)}, the first and second rows contain indices of the starting and ending points. The column is sorted in the way that for the $k$-th edge, \mcode{edge(k,1)<edge(k,2)}. The indices of these edges are are marked with purple boxes as shown in Fig.~\ref{Fig:auxstructure}.

Following \cite{ChenL-iFEM-2009}, we shall use the name convention a2b to represent the link form a to b. This link is
usually the map from the local index set to the global index set. Throughout this paper, we use the symbols \mcode{N}, \mcode{NT}, and \mcode{NE} to represent the number of nodes, elements, and edges.

\textbf{elem2edge}. The cell array \mcode{elem2edge} establishes the map of local index of edges in each polygon to its global index in matrix \mcode{edge}. For instance, \mcode{elem2edge\{1\} = [13,7,6,3,5]} for the mesh in Fig.~\ref{Fig:auxstructure}.

\textbf{bdEdge}. This matrix exacts the boundary edges from \mcode{edge}.

\textbf{edge2elem}. The matrix \mcode{edge(1:NE,1:2)} records the neighboring polygons for each edge. In Fig.~\ref{Fig:auxstructure}, \mcode{edge(3,1:2) = [1,2]}. For a boundary edge, the outside is specified to the element sharing it as an edge.

\textbf{neighbor}. We use the cell array \mcode{neighbor} to record the neighboring polygons for each element. For example, \mcode{neighbor\{4\} = [5,1,2,3]}, where the $i$-th entry corresponds to the $i$-th edge of the current element. Note that for a boundary edge, the neighboring polygon is specified to the current element, for instance, \mcode{neighbor\{1\} = [5,1,1,2,4]}.

\textbf{node2elem}. This cell array finds the elements sharing a common nodes. For example, \mcode{node2elem\{2\} = [1,2,4]}.

In addition, we provide a subroutine \mcode{auxgeometry.m} to compute some useful geometric quantities, such as the barycenter \mcode{centroid}, the diameter \mcode{diameter} and the area \mcode{area} of each element.

\subsubsection{Computation of the elliptic projection}

\subsubsection*{Transition matrix}

The shape functions of $V_k(K)$ are written in the following compact notation
\[\bb{\phi}^T = (\bb{\phi}_1, \bb{\phi}_2, \cdots, \bb{\phi}_{N_k}),\]
where $N_k = 2N_v+1$ is the cardinality of the local basis set. The basis of $\widehat{V}_k(K)$ is given by
\[\bb{\widehat{\bb{m}}}^T = (\bb{\widehat{\bb{m}}}_1, \bb{\widehat{\bb{m}}}_2, \cdots, \bb{\widehat{\bb{m}}}_{\widehat{N}_p}).\]
Noting that $\widehat{V}_k(K)\subset V_k(K)$, we set
$\bb{\widehat{\bb{m}}}^T = \bb{\phi}^T\bb{D}$,
where $\bb{D}$ is referred to as the transition matrix from the elliptic projection space $\widehat{V}_k(K)$ to the virtual element space $V_k(K)$. By the definition of the d.o.f.s,
\[\widehat{\bb{m}}_\alpha = \sum\limits_{i=1}^{N_k}\bb{\phi}_i\bb{D}_{i\alpha}, \quad \bb{D}_{i\alpha} = \chi_i(\widehat{\bb{m}}_\alpha).\]

For every $\bb{v}\in \widehat{V}_k(K)$, by definition, $\bb{v}\in \mathbb{K}\nabla \mathbb{P}_{k+1}(K)$. For $k=1$, the scaled monomials of $\mathbb{P}_{k+1}(K) = \mathbb{P}_2(K)$ are
\[m^T = (m_1,m_2,\cdots,m_{N_m}), \quad N_m = 6,\]
with
\[m_1(x,y) = 1,~~m_2(x,y) = \frac{x - x_K}{h_K},~~m_3(x,y) = \frac{y - y_K}{h_K},\]
\[m_4(x,y) = \frac{(x - x_K)^2}{h_K^2},~~m_5(x,y) = \frac{(x - x_K)(y - y_K)}{{h_K^2}},~~
m_6(x,y) = \frac{(y-y_K)^2}{h_K^2},\]
where $(x_K, y_K)$ and $h_K$ are the barycenter and diameter of $K$, respectively.
We then choose
\begin{equation}\label{widehatm}
\widehat{\bb{m}}_\alpha = \mathbb{K}\nabla(h_K m_{\alpha+1}), \quad \alpha = 1,\cdots, \widehat{N}_m=N_m-1=5,
\end{equation}
with the explicit expressions given by
\begin{align*}
\widehat{\bb{m}}_1 = \mathbb{K}\begin{bmatrix} 1 \\ 0 \end{bmatrix}, \quad
\widehat{\bb{m}}_2 = \mathbb{K}\begin{bmatrix} 0 \\ 1 \end{bmatrix}, \quad
\widehat{\bb{m}}_3 = \mathbb{K}\begin{bmatrix} 2m_2 \\ 0 \end{bmatrix}, \quad
\widehat{\bb{m}}_4 = \mathbb{K}\begin{bmatrix} m_3 \\ m_2 \end{bmatrix}, \quad
\widehat{\bb{m}}_5 = \mathbb{K}\begin{bmatrix} 0 \\ 2m_3 \end{bmatrix}.
\end{align*}

We now compute the transition matrix.
For $\chi_i$ with $i=1,\cdots,N_v$ in \eqref{dof1}, noting that $\bb{v}\cdot\bb{n}|_e\in \mathbb{P}_1(e)$, the trapezoidal rule gives
\begin{equation}\label{chidof1}
\chi_i(\bb{v}) = \int_{e_i} (\bb{v}\cdot\bb{n}){\rm d}s = \frac{1}{2}(\bb{v}(z_i)+\bb{v}(z_{i+1}))\cdot(h_{e_i}\bb{n}_{e_i}), \quad \bb{v}=\widehat{\bb{m}}_\alpha.
\end{equation}
For $\chi_{N_v+i}$ with $i=1,\cdots,N_v$ in \eqref{dof2}, using the Simpson formula yields
\begin{align}
\chi_{N_v+i}(\bb{v})
& = \int_{e_i} (\bb{v}\cdot\bb{n})\frac{s-s_{e_i}}{h_{e_i}}{\rm d}s =
\frac{h_{e_i}}{6}(f(z_i)+4f(s_{e_i})+f(z_{i+1})) \nonumber\\
& = \frac{h_{e_i}}{6}(f(z_i)+f(z_{i+1})), \quad f = (\bb{v}\cdot\bb{n})\frac{s-s_{e_i}}{h_{e_i}}\nonumber\\
& = \frac{1}{6}\Big(\bb{v}(z_i)\cdot\bb{n}_{e_i}(s_i-s_{e_i})+\bb{v}(z_{i+1})\cdot\bb{n}_{e_i}(s_{i+1}-s_{e_i})\Big)\nonumber\\
& = \frac{1}{6}\Big(\bb{v}(z_i)\cdot\bb{n}_{e_i}\frac{-1}{2}h_{e_i}+\bb{v}(z_{i+1})\cdot\bb{n}_{e_i}\frac{1}{2}h_{e_i}\Big)\nonumber\\
& = \frac{1}{12}\Big(\bb{v}(z_{i+1})-\bb{v}(z_i)\Big)\cdot(\bb{n}_{e_i}h_{e_i}), \quad \bb{v}=\widehat{\bb{m}}_\alpha. \label{chidof2}
\end{align}
For $\chi_{2N_v+1}$ in \eqref{dof3}, the integration by parts gives
\begin{align*}
\chi_{2N_v+1}(\bb{v})
& = \int_K {\rm rot}\bb{v}{\rm d}x = \int_K (\partial_1v_2-\partial_2v_1){\rm d}x \\
& = \int_{\partial K}(v_2n_1-v_1n_2){\rm d}s = \int_{\partial K}(\bb{v}\cdot\bb{t}){\rm d}s\\
& = \sum\limits_{i=1}^{N_v} \frac{1}{2}(\bb{v}(z_i)+\bb{v}(z_{i+1}))\cdot(h_{e_i}\bb{t}_{e_i}), \quad \bb{v}=\widehat{\bb{m}}_\alpha.
\end{align*}

According to the above discussion, the transition matrix is calculated as follows.
\vspace{-0.8cm}
\begin{lstlisting}
K = pde.K; % coefficient matrix
    % ------- element information ----------
    index = elem{iel};     Nv = length(index);
    xK = centroid(iel,1); yK = centroid(iel,2); hK = diameter(iel);
    x = node(index,1); y = node(index,2);
    v1 = 1:Nv; v2 = [2:Nv,1]; % loop index for vertices or edges
    xe = (x(v1)+x(v2))/2;  ye = (y(v1)+y(v2))/2; % mid-edge points
    Ne = [y(v2)-y(v1), x(v1)-x(v2)]; % he*ne
    Te = [-Ne(:,2), Ne(:,1)]; % he*te

    % ------- scaled monomials --------
    m2 = @(x,y) (x-xK)./hK;
    m3 = @(x,y) (y-yK)./hK;
    m4 = @(x,y) (x-xK).^2./hK^2;
    m5 = @(x,y) (x-xK).*(y-yK)./hK^2;
    m6 = @(x,y) (y-yK).^2./hK^2;
    % \hat{m}_a = K*grad(hK*m_{a+1})
    mh1 = @(x,y) [K(1,1)+0*x; K(2,1)+0*x];
    mh2 = @(x,y) [K(1,2)+0*x; K(2,2)+0*x];
    mh3 = @(x,y) [2*K(1,1)*m2(x,y); 2*K(2,1)*m2(x,y)];
    mh4 = @(x,y) [K(1,1)*m3(x,y)+K(1,2)*m2(x,y); K(2,1)*m3(x,y)+K(2,2)*m2(x,y)];
    mh5 = @(x,y) [2*K(1,2)*m3(x,y); 2*K(2,2)*m3(x,y)];
    mh = @(x,y) [mh1(x,y), mh2(x,y), mh3(x,y), mh4(x,y), mh5(x,y)];

    % -------- transition matrix ----------
    NdofA = 2*Nv+1; Nmh = 5;
    D = zeros(NdofA,Nmh);
    for i = 1:Nv % loop of edges
        % v at z_i, z_{i+1} for v = [mK1,...,mK5]
        va = mh(x(v1(i)),y(v1(i))); vb = mh(x(v2(i)),y(v2(i)));
        % chi_i, i = 1,...,Nv
        D(i,:) =  1/2*Ne(i,:)*(va+vb);
        % chi_{Nv+i}, i = 1,...,Nv
        D(Nv+i,:) =  1/12*Ne(i,:)*(vb-va);
        % chi_{2Nv+1}
        D(end,:) = D(end,:) + 1/2*Te(i,:)*(va+vb);
    end
\end{lstlisting}
Here, \mcode{K} is for $\mathbb{K}$ and \mcode{mhat} is for $\bb{\widehat{\bb{m}}}^T$.

\subsubsection*{Elliptic projection matrices}

We denote the matrix representation of the $\widehat{\Pi}^K$-projection by $\bb{\widehat{\Pi}}^K$ in the sense that
\[\widehat{\Pi}^K(\bb{\phi}_1, \bb{\phi}_2, \cdots, \bb{\phi}_{N_k})
= (\bb{\phi}_1, \bb{\phi}_2, \cdots, \bb{\phi}_{N_k})\bb{\widehat{\Pi}}^K \quad \mbox{or} \quad
\widehat{\Pi}^K\bb{\phi}^T = \bb{\phi}^T\bb{\widehat{\Pi}}^K.\]
By the definition of d.o.f.s, the $j$-th column of $\bb{\widehat{\Pi}}^K$ is the d.o.f vector of $\widehat{\Pi}^K\bb{\phi}_j$, i.e.,
$\bb{\widehat{\Pi}}^K = \Big( \chi_i(\widehat{\Pi}^K\bb{\phi}_j) \Big)$.
The elliptic projection vector $\widehat{\Pi}^K\bb{\phi}^T$ can be expanded in the basis $\bb{\widehat{\bb{m}}}^T$ of the elliptic projection space $\widehat{V}_k(K)$ as
$\widehat{\Pi}^K\bb{\phi}^T = \widehat{\bb{m}}^T\bb{\widehat{\Pi}}_*^K$.
It is easy to check that
$\bb{\widehat{\Pi}}^K = \bb{D}\bb{\widehat{\Pi}}_*^K$.

The definition \ref{ellipticProjDarcy} is equivalent to
\[
a^K(\widehat{\bb{m}}, \widehat{\Pi}^K\bb{\phi}^T) = a^K(\widehat{\bb{m}}, \bb{\phi}^T)
 \quad  \mbox{or} \quad
\widehat{\bb{G}}\bb{\widehat{\Pi}}_*^K = \widehat{\bb{B}},\]
where
\[\widehat{\bb{G}} = a^K(\widehat{\bb{m}}, \widehat{\bb{m}}^T), \quad \widehat{\bb{B}} = a^K(\widehat{\bb{m}}, \bb{\phi}^T).\]
We also have the consistency relation $\widehat{\bb{G}} = \widehat{\bb{B}}\bb{D}$.

We now compute $\widehat{\bb{B}}$. From \eqref{integrationdarcy} and \eqref{widehatm}, one has
\begin{align*}
\widehat{\bb{B}}_{\alpha i}
& = a^K(\widehat{\bb{m}}_\alpha, \bb{\phi}_i) = - h_K\int_K m_{\alpha+1}{\rm div}\bb{\phi}_i {\rm d}x + h_K\int_{\partial K}m_{\alpha+1}\bb{\phi}_i\cdot\bb{n}{\rm d}s \\
& =: -h_K I_1(\alpha,i) + h_K I_2(\alpha,i),
\end{align*}
where
\[I_1(\alpha,i) = \int_K m_{\alpha+1}{\rm div}\bb{\phi}_i {\rm d}x, \quad I_2(\alpha,i) = \int_{\partial K}m_{\alpha+1}\bb{\phi}_i\cdot\bb{n}{\rm d}s.\]
For $I_1$, noting that ${\rm div}\bb{\phi}_i$ is constant when $k=1$, set
${\rm div}\bb{\phi}_i = c_i m_1 = c_i$, which gives
\[
c_i = |K|^{-1}\int_K {\rm div}\bb{\phi}_i{\rm d}x = |K|^{-1}\int_{\partial K} \bb{\phi}_i\cdot\bb{n}{\rm d}s=
\begin{cases}
|K|^{-1}, \quad & i = 1,\cdots, N_v, \\
0, \quad & i>N_v,
\end{cases}
\]
and hence
\[I_1(\alpha,i) = \int_K m_{\alpha+1}{\rm div}\bb{\phi}_i {\rm d}x = c_i\int_K m_{\alpha+1}{\rm d}x. \]
The first term is now computed in MATLAB as follows.
\vspace{-0.8cm}
\begin{lstlisting}
    % first term
    nodeT = [node(index,:);centroid(iel,:)];
    elemT = [(Nv+1)*ones(Nv,1),(1:Nv)',[2:Nv,1]'];
    m = @(x,y) [m2(x,y), m3(x,y), m4(x,y), m5(x,y), m6(x,y)]; % m_{a+1},...
    ci = zeros(1,NdofA); ci(1:Nv) = 1/area(iel);
    Intm = integralTri(m,3,nodeT,elemT);
    I1 = Intm'*ci;
\end{lstlisting}
The subroutine \mcode{integralTri.m} calculates the integral on a polygonal element which is triangulated with the basic data structures \mcode{nodeT} and \mcode{elemT}.

For $I_2$, since $\bb{\phi}_i\cdot\bb{n}|_{e_j}\in\mathbb{P}_1({e_j})$~($i=1,\cdots,2N_v+1$), set
\[\bb{\phi}_i\cdot\bb{n}_{e_j}|_{e_j} = c_0^i+c_1^i\frac{s-s_{e_j}}{h_{e_j}},\]
where $s$ it the natural parameter of the edge $e_j$ and $s_{e_j}$ is the mid-point in the local parametrization.
Noting that
\[\int_e\Big(\frac{s-s_e}{h_e}\Big)^\alpha \Big(\frac{s-s_e}{h_e}\Big)^\beta {\rm d}s = \frac{h_e}{(\alpha+\beta+1)}
\Big( \frac{1}{2^{\alpha+\beta+1}} - \frac{1}{(-2)^{\alpha+\beta+1}}\Big),\]
we then obtain
\[
\begin{cases}
h_{e_j}c_0^i = \int_{e_j}\bb{\phi}_i\cdot\bb{n}{\rm d}s=:r_j^i,  \\
\frac{1}{12}h_{e_j}c_1^i = \int_{e_j}\bb{\phi}_i\cdot\bb{n}\frac{s-s_{e_j}}{h_{e_j}}{\rm d}s=:s_j^i.
\end{cases}
\]
By definition, the nonzero elements of $r_j^i$ and $s_j^i$ are
\[r_j^j = s_j^{j+N_v} = 1, \quad j=1,\cdots, N_v.\]
We have
\[c_0^i = \frac{1}{h_{e_j}}r_j^i, \quad c_1^i = \frac{12}{h_{e_j}}s_j^i,\]
and hence
\[\bb{r}_j^i(s):=\bb{\phi}_i\cdot(h_{e_j}\bb{n}_{e_j})|_{e_j} = r_j^i+12s_j^i\frac{s-s_{e_j}}{h_{e_j}}.\]
It is obvious that
\[\bb{r}_j^i(s_a)=r_j^i-6s_j^i, \quad \bb{r}_j^i(s_e)=r_j^i, \quad \bb{r}_j^i(s_b)= r_j^i+6s_j^i,\]
where, $s_a$, $s_e$ and $s_b$ correspond to the starting point, middle point and ending point in the local parameterization, respectively.

Observing that the Simpson formula is accurate for cubic polynomials, we have
\begin{align*}
I_2(\alpha,i)
& = \sum\limits_{j=1}^{N_v}\int_{e_j}m_{\alpha+1}\bb{\phi}_i\cdot\bb{n}{\rm d}s \\
& = \sum\limits_{j=1}^{N_v}\frac{h_{e_j}}{6}(m_{\alpha+1}\bb{\phi}_i(z_j) + 4m_{\alpha+1}\bb{\phi}_i(z_{e_j}) + m_{\alpha+1}\bb{\phi}_i(z_{j+1})\cdot\bb{n}_{e_j}\\
& = \frac{1}{6}\sum\limits_{j=1}^{N_v}m_{\alpha+1}(z_j)\bb{r}_j^i(s_a) + 4m_{\alpha+1}(z_{e_j})\bb{r}_j^i(s_e)+ m_{\alpha+1}(z_{j+1})\bb{r}_j^i(s_b)
\end{align*}
or
\[I_2(\alpha,:) = \frac{1}{6}\sum\limits_{j=1}^{N_v}m_{\alpha+1}(z_j)\bb{r}_j(s_a) + 4m_{\alpha+1}(z_{e_j})\bb{r}_j(s_e)+ m_{\alpha+1}(z_{j+1})\bb{r}_j(s_b).\]
Then the second term and the matrices $\widehat{\bb{B}}$ and $\widehat{\bb{G}}$ can be computed as follows.
\vspace{-0.8cm}
\begin{lstlisting}
    % second term
    rij = zeros(Ndof,Nv); rij(1:Nv,:) = eye(Nv);
    sij = zeros(Ndof,Nv); sij(Nv+1:2*Nv,:) = eye(Nv);
    rija = rij- 6*sij;
    rije = rij;
    rijb = rij + 6*sij;
    I2 = zeros(np,Ndof);
    for j = 1:Nv
        ma = m(x(v1(j)),y(v1(j)));
        me = m(xe(j),ye(j));
        mb = m(x(v2(j)),y(v2(j)));
        rja = rija(:,j)';  rje = rije(:,j)';  rjb = rijb(:,j)';
        I2 = I2 + 1/6*(ma'*rja + 4*me'*rje + mb'*rjb);
    end
    % \hat{B} and \hat{G}
    Bs = hK*(-I1+I2); Gs = Bs*D
\end{lstlisting}

\subsubsection{Computation and assembly of the stiffness matrix and load vector}

\subsubsection*{The stiffness matrix}

The local stiffness matrix of $A$ in \eqref{linearsystemDarcy} is $A^K = A^K_1 + A^K_2$, where
\begin{align*}
A^K_1(i,j) & = a^K(\widehat{\Pi}^K\bb{\phi}_j,\widehat{\Pi}^K\bb{\phi}_i), \\
A^K_2(i,j)
 & = \|\mathbb{K}^{-1}\|S^K(\bb{\phi}_j-\widehat{\Pi}^K\bb{\phi}_j,\bb{\phi}_i-\widehat{\Pi}^K\bb{\phi}_i) \\
 & = \|\mathbb{K}^{-1}\|\sum\limits_{i=1}^{N_k}\chi_r(\bb{\phi}_j-\widehat{\Pi}^K\bb{\phi}_j)
\chi_r(\bb{\phi}_i-\widehat{\Pi}^K\bb{\phi}_i).
\end{align*}
The first term is written in matrix form as
\[A_K^1 = a^K(\widehat{\Pi}^K\bb{\phi},\widehat{\Pi}^K\bb{\phi}^T)
=(\bb{\widehat{\Pi}}_*^K)^Ta^K(\widehat{\bb{m}},\widehat{\bb{m}}^T)\bb{\widehat{\Pi}}_*^K = (\bb{\widehat{\Pi}}_*^K)^T\widehat{\bb{G}}\bb{\widehat{\Pi}}_*^K.\]
For the second term, from $\chi_r(\widehat{\Pi}^K\bb{\phi}_i) = (\bb{\widehat{\Pi}}^K)_{ri}$ one has
\[A_K^2 = \|\mathbb{K}^{-1}\|({\bf{I}}-\bb{\widehat{\Pi}}^K)^T({\bf{I}}-\bb{\widehat{\Pi}}^K).\]

\begin{figure}[!htb]
  \centering
  \includegraphics[scale=0.8]{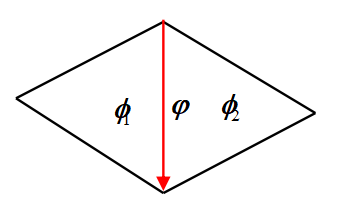}\\
  \caption{Illustration of the global and local basis functions}\label{fig:signBase}
\end{figure}

Note that the signs of the d.o.f.s in \eqref{dof1} vary with the edge orientation, which means the global basis function $\bb{\varphi}_i$ restriction to the element may have the opposite sign with the local basis function $\bb{\phi}_i$. As shown in Fig.~\ref{fig:signBase}, $e$ is an interior edge with the global orientation given by the arrow. Let $\bb{\varphi}$ be the global basis function, and let $\bb{\phi}_1$ and $\bb{\phi}_2$ be the local basis functions corresponding to the left and right elements, respectively. Then one has
\[\bb{\phi}_1 = - \bb{\varphi}|_{K_1}, \quad \bb{\phi}_2 = \bb{\varphi}|_{K_2}.\]
Obviously, the sign can be determined by computing the difference between the two indices of the end points of $e$. For this reason, one only needs to compute the elementwise signs of the edges of $\partial K$. The code is given as
\vspace{-0.8cm}
\begin{lstlisting}
%% Get elementwise signs of basis functions
bdEdgeIdx = bdStruct.bdEdgeIdx;  E = false(NE,1); E(bdEdgeIdx) = 1;
sgnBase = cell(NT,1);
for iel = 1:NT
    index = elem{iel};   Nv = length(index);   NdofA = 2*Nv+1;
    sgnedge = sign(diff(index([1:Nv,1])));
    id = elem2edge{iel}; sgnbd = E(id); sgnedge(sgnbd) = 1;
    sgnelem = ones(NdofA,1); sgnelem(1:Nv) = sgnedge;
    sgnBase{iel} = sgnelem;
end
\end{lstlisting}
Note that the positive signs of the boundary edges are recovered in the above code.

Since
$a^K(\pm\bb{\phi}_i,\pm\bb{\phi}_j) = \pm \cdot \pm a^K(\bb{\phi}_i,\bb{\phi}_j)$,
we can introduce a signed stiffness matrix \mcode{sgnK} to add the correct signs to the local stiffness matrix. The matrix $A^K$ can be computed in the following way.
\vspace{-0.8cm}
\begin{lstlisting}
    % -------- sign matrix and sign vector -------
    sgnelem = sgnBase{iel};
    sgnK = sgnelem*sgnelem';

    % ------------- stiffness matrix -------------
    % Projection
    Pis = Gs\Bs;   Pi  = D*Pis;   I = eye(size(Pi));
    % Stiffness matrix A
    AK  = Pis'*Gs*Pis + norm(inv(K),'fro')*(I-Pi)'*(I-Pi);  % G = Gs
    AK = AK.*sgnK;
    AK = reshape(AK',1,[]); % straighten as row vector for easy assembly
\end{lstlisting}

By the definition of $Q_h$, on each element $K$, the basis function $\psi_l\in \mathbb{P}_{k-1}(K)$. Only one basis $\psi_l = m_1 = 1$ for $k=1$. Hence, the element matrix $B^K$ of $B$ in \eqref{linearsystemDarcy} is a column vector of size $N_k\times 1$ (Note that $B$ is of size $\mcode{(2NE+NT)} \times \mcode{NT}$), and
\begin{align*}
B^K_{j1}
& = b^K(\bb{\phi}_j, m_1)=\int_Km_1{\rm div}\bb{\phi}_j{\rm d}x = \int_K{\rm div}\bb{\phi}_j{\rm d}x\\
& = \int_{\partial K}\bb{\phi}_j\cdot\bb{n}{\rm d}s =
\begin{cases}
1, \quad & 1\le j\le N_v, \\
0, \quad & j>N_v.
\end{cases}
\end{align*}
The sign of each entry is adjusted by using
$b^K(\pm\bb{\phi}_i,\psi_j) = \pm b^K(\bb{\phi}_i,\psi_j)$.
\vspace{-0.8cm}
\begin{lstlisting}
    % Stiffness matrix B
    BK = zeros(NdofA,1); BK(1:Nv) = 1;
    BK = BK.*sgnelem;
    BK = reshape(BK',1,[]); % straighten as row vector for easy assembly
\end{lstlisting}

For $k=1$, $Q_h$ is piecewise constant, and hence has \mcode{NT} basis functions given by
\[\psi_l(x) =
\begin{cases}
1, \quad & x \in K_l, \\
0, \quad & \mbox{otherwise}.
\end{cases}\]
Obviously, the vector in \eqref{dl} is
\[d_l = \int_{K_l} {\rm d}x = |K_l|, \quad l = 1,\cdots,M = \mcode{NT}.\]

We compute the elliptic projections and provide the assembly index by looping over the elements. The assembly index for the matrices $A$ and $B$ is given by
\vspace{-0.8cm}
\begin{lstlisting}
    % ------ assembly index for bilinear forms --------
    NdofA = 2*Nv+1; NdofB = 1;
    indexDofA = [elem2edge{iel},elem2edge{iel}+NE,iel+2*NE];
    indexDofB = iel;
    iiA(idA+1:idA+NdofA^2) = reshape(repmat(indexDofA, NdofA,1), [], 1);
    jjA(idA+1:idA+NdofA^2) = repmat(indexDofA(:), NdofA, 1);
    ssA(idA+1:idA+NdofA^2) = AK(:);
    idA = idA + NdofA^2;
    iiB(idB+1:idB+NdofA*NdofB) = reshape(repmat(indexDofA, NdofB,1), [], 1);
    jjB(idB+1:idB+NdofA*NdofB) = repmat(indexDofB(:), NdofA, 1);
    ssB(idB+1:idB+NdofA*NdofB) = BK(:);
    idB = idB + NdofA*NdofB;
\end{lstlisting}
Afterwards, we can assemble the matrices $A$ and $B$ using the MATLAB functions \mcode{sparse}.
\vspace{-0.8cm}
\begin{lstlisting}
A = sparse(iiA,jjA,ssA,NNdofA,NNdofA);
B = sparse(iiB,jjB,ssB,NNdofA,NNdofB);
d = area;  % for Lagrange multiplier
\end{lstlisting}

\subsubsection*{The load vector}

For the right-hand side, in view of \eqref{linearsystemDarcy}, we only need to compute $\bb{f} = (-(f,\psi_l))$. For $k=1$, the local vector is
\[\bb{f}_K = -(f,\psi_l)_K = -\int_K f {\rm d}x,\]
with the realization reading
\vspace{-0.8cm}
\begin{lstlisting}
    % -------- load vector f ----------
    fxy = @(x,y) pde.f([x,y]); % f(p) = f([x,y])
    rhs = integralTri(fxy,3,nodeT,elemT); rhs = rhs';
    fK = -rhs;
\end{lstlisting}

The assembly index for the vector $\bb{f}$ is given by
\vspace{-0.8cm}
\begin{lstlisting}
    % ------- assembly index for rhs -------
    elemb(ib+1:ib+NdofB) = indexDofB(:);
    Fb(ib+1:ib+NdofB) = fK(:);
    ib = ib + NdofB;
\end{lstlisting}
Then $\bb{f}$ can be assembled using the MATLAB functions \mcode{accumarray}.
\vspace{-0.8cm}
\begin{lstlisting}
    FB = accumarray(elemb,ffB,[NNdofB 1]);
\end{lstlisting}
To sum up, the linear system without boundary conditions imposed is given by
\vspace{-0.8cm}
\begin{lstlisting}
%% Get block linear system
kk = sparse(NNdof+1,NNdof+1);  ff = zeros(NNdof+1,1);
kk(1:NNdofA,1:NNdofA) = A;
kk(1:NNdofA, (1:NNdofB)+NNdofA) = B;
kk((1:NNdofB)+NNdofA, 1:NNdofA) = B';
kk((1:NNdofB)+NNdofA, end) = d;
kk(end, (1:NNdofB)+NNdofA) = d';
ff((1:NNdofB)+NNdofA) = FB;
\end{lstlisting}
Note that the extra variable $\lambda$ leads to \mcode{NNdof+1} rows or columns.

\subsubsection{Applying the boundary conditions}

The boundary condition $\bb{u}\cdot\bb{n} = g$ is now viewed as a Dirichlet condition for $\bb{u}$, which provides the values of the first two types of d.o.f.s, i.e.,
\begin{align*}
&\chi_i(\bb{u}) = \int_{e_i} (\bb{u}\cdot\bb{n}){\rm d}s, \quad i = 1,\cdots, N_v, \\
&\chi_{N_v+i}(\bb{u}) = \int_{e_i} (\bb{u}\cdot\bb{n})\frac{s-s_{e_i}}{h_{e_i}}{\rm d}s, \quad i = 1,\cdots, N_v
\end{align*}
The Simpson rule is used to approximate the exact ones as follows.
\vspace{-0.8cm}
\begin{lstlisting}
%% Apply Dirichlet boundary conditions
% bdDof, freeDof
bdEdge = bdStruct.bdEdge;  bdEdgeIdx = bdStruct.bdEdgeIdx;
id = [bdEdgeIdx; bdEdgeIdx+NE];
isBdDof = false(NNdof+1,1); isBdDof(id) = true;
bdDof = (isBdDof); freeDof = (~isBdDof);
% bdval
u = pde.uexact;
z1 = node(bdEdge(:,1),:); z2 = node(bdEdge(:,2),:); ze = (z1+z2)./2;
e = z1-z2;  % e = z2-z1
Ne = [-e(:,2),e(:,1)];
chi1 = 1/6*sum((u(z1)+4*u(ze)+u(z2)).*Ne,2); % u*n = g
chi2 = 1/12*sum((u(z2)-u(z1)).*Ne,2);
bdval = [chi1;chi2];
% sol
sol = zeros(NNdof+1,1);
sol(bdDof) = bdval;
ff = ff - kk*sol;

%% Set solver
sol(freeDof) = kk(freeDof,freeDof)\ff(freeDof);
uh = sol(1:NNdofA); % u = [u1,u2]
ph = sol(NNdofA+1:end-1);

%% Store information for computing errors
info.Ph = Ph; info.elem2dof = elem2dof;
info.kk = kk; %info.freeDof = freeDof;
\end{lstlisting}

To compute the discrete errors, one can store the matrix representation $\bb{\widehat{\Pi}}_*^K$ and the assembly index \mcode{elem2dof} in the function file.

\subsubsection{Numerical example}

In this paper, all examples are implemented in MATLAB R2019b. Our code is available from GitHub (\url{https://github.com/Terenceyuyue/mVEM}). The subroutine \mcode{Darcy\_mixedVEM.m} is used to compute the numerical solutions and the test script \mcode{main\_Darcy\_mixedVEM.m} verifies the convergence rates. The PDE data is generated by \mcode{Darcydata.m}.

\begin{example}\label{ex:Darcy}
Let $\mathbb{K}$ be the identity matrix. The right-hand side $f$ and the boundary conditions are chosen in such a way that the exact solution of \eqref{DarcyOriginal} is
$p(x,y) = \sin(\pi x) \cos(\pi y)$.
\end{example}

\begin{figure}[!htb]
  \centering
  \includegraphics[scale=0.8,trim=80 0 80 0,clip]{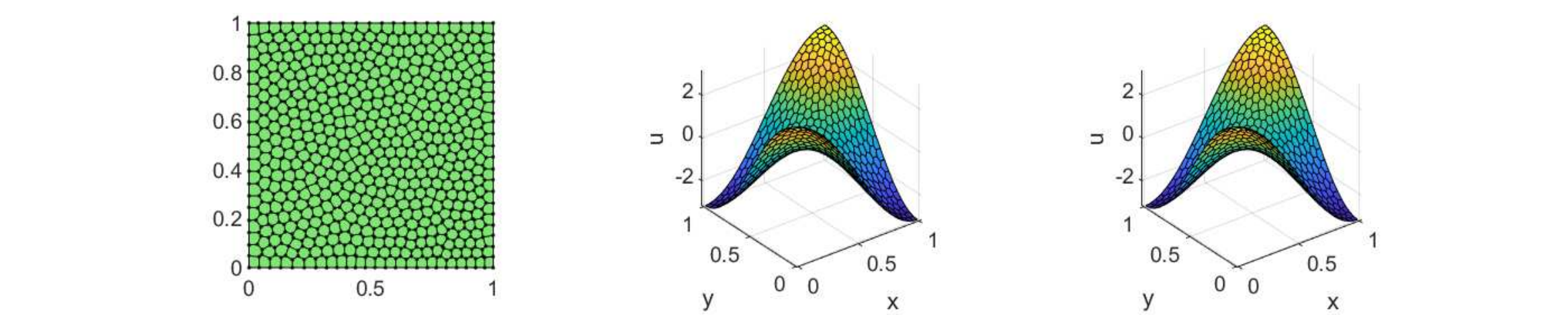}\\
  \caption{Numerical and exact results for the Darcy problem}\label{fig:Darcy}
\end{figure}

\begin{figure}[!htb]
  \centering
  \includegraphics[scale=0.45]{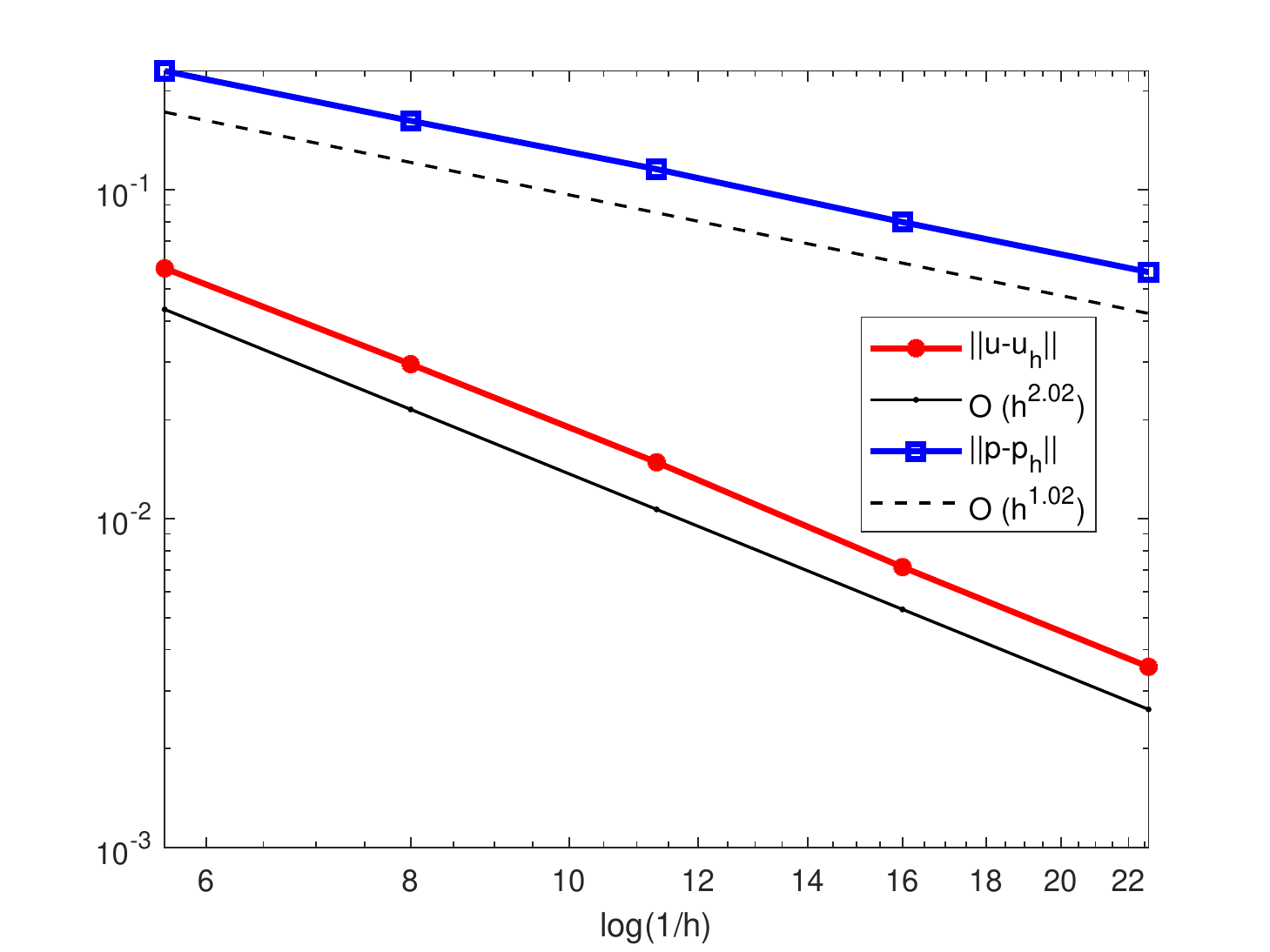}\\
  \caption{The convergence rates for the Darcy problem}\label{fig:Darcyrate}
\end{figure}

To test the accuracy of the proposed method we consider a sequence of polygonal meshes, which is a Centroidal Voronoi Tessellation of the unit square in 32, 64, 128, 256 and 512 polygons. These meshes are generated by the MATLAB toolbox - PolyMesher introduced in \cite{Talischi-Paulino-Pereira-2012}. we report the nodal values of the exact solution and the piecewise elliptic projection $\widehat{\Pi}^K \bb{u}_h$ in Fig.~\ref{fig:Darcy}.
The convergence orders of the errors against the mesh size $h$ are shown in Fig.~\ref{fig:Darcyrate}. Generally speaking, $h$ is proportional to $N^{-1/2}$, where $N$ is the total number of elements in the mesh. The convergence rate with respect to $h$ is estimated by assuming $\text{ErrL2}(h) = ch^{\alpha}$, and by computing a least squares fit to this log-linear relation.
As observed from Fig.~\ref{fig:Darcyrate}, the convergence rate of $p$ is linear with respect to the $L^2$ norm, and the VEM ensures the quadratic convergence for $\bb{u}$ in the $L^2$ norm, which is consistent with the theoretical prediction in \cite{Beirao-Brezzi-Cangiani-2014}.

\subsubsection{A lifting mixed virtual element method}

We can obtain the second-order convergence of the variable $p$ using the lifting technique. To do so, we modify the original virtual element space to the following lifting space
\begin{align*}
  V_k(K)
  & = \{\bb{v}\in H({\rm div};K)\cap H({\rm rot};K): \bb{v}\cdot\bb{n}|_e\in\mathbb{P}_k(e), \\
  &  \hspace{3cm} {\rm div}\bb{v}|_K\in \mathbb{P}_k(K),~~ {\rm rot}\bb{v}|_K\in \mathbb{P}_{k-1}(K)\},
\end{align*}
with the d.o.f.s given by
\begin{align*}
&\int_e \bb{v}\cdot\bb{n} q{\rm d}s, \quad q\in \mathbb{M}_k(e),~~e\subset\partial K, \\
&\int_K \bb{v}\cdot\nabla q{\rm d}x, \quad q\in  \mathbb{M}_k(K)\backslash\{1\}, \\
&\int_K {\rm rot}\bb{v}q{\rm d}x, \quad q\in \mathbb{M}_{k-1}(K) .
\end{align*}
In the lowest order case $k=1$, the local d.o.f.s are arranged as
\begin{align*}
&\chi_i(\bb{v}) = \int_{e_i} (\bb{v}\cdot\bb{n}){\rm d}s, \quad i = 1,\cdots, N_v, \\
&\chi_{N_v+i}(\bb{v}) = \int_{e_i} (\bb{v}\cdot\bb{n})\frac{s-s_{e_i}}{h_{e_i}}{\rm d}s, \quad i = 1,\cdots, N_v,\\
& \chi_{2N_v+1}(\bb{v}) = \int_K \bb{v}\cdot\nabla m_2{\rm d}x= \frac{1}{h_K}\int_K \bb{v}_1{\rm d}x,\\
& \chi_{2N_v+2}(\bb{v}) = \int_K \bb{v}\cdot\nabla m_3{\rm d}x = \frac{1}{h_K}\int_K \bb{v}_2{\rm d}x,\\
&\chi_{2N_v+3}(\bb{v}) = \int_K {\rm rot}\bb{v}{\rm d}x.
\end{align*}
In this case, $Q_h$ will be replaced by the piecewise linear space.

\begin{figure}[!htb]
  \centering
  \includegraphics[scale=0.8,trim=80 0 80 0,clip]{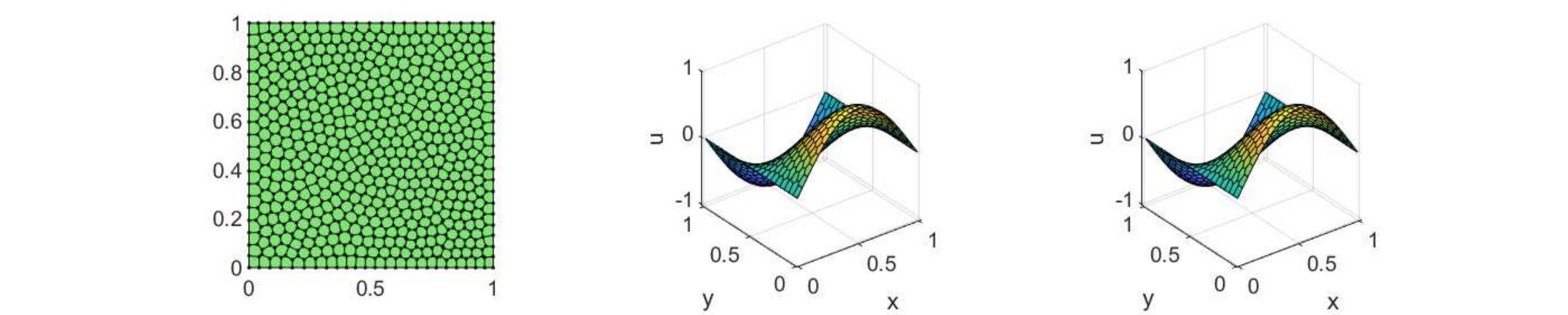}\\
  \caption{$p_h$ and $p$ for the lifting VEM of the Darcy problem}\label{fig:Darcyp}
\end{figure}

\begin{figure}[!htb]
  \centering
  \includegraphics[scale=0.45]{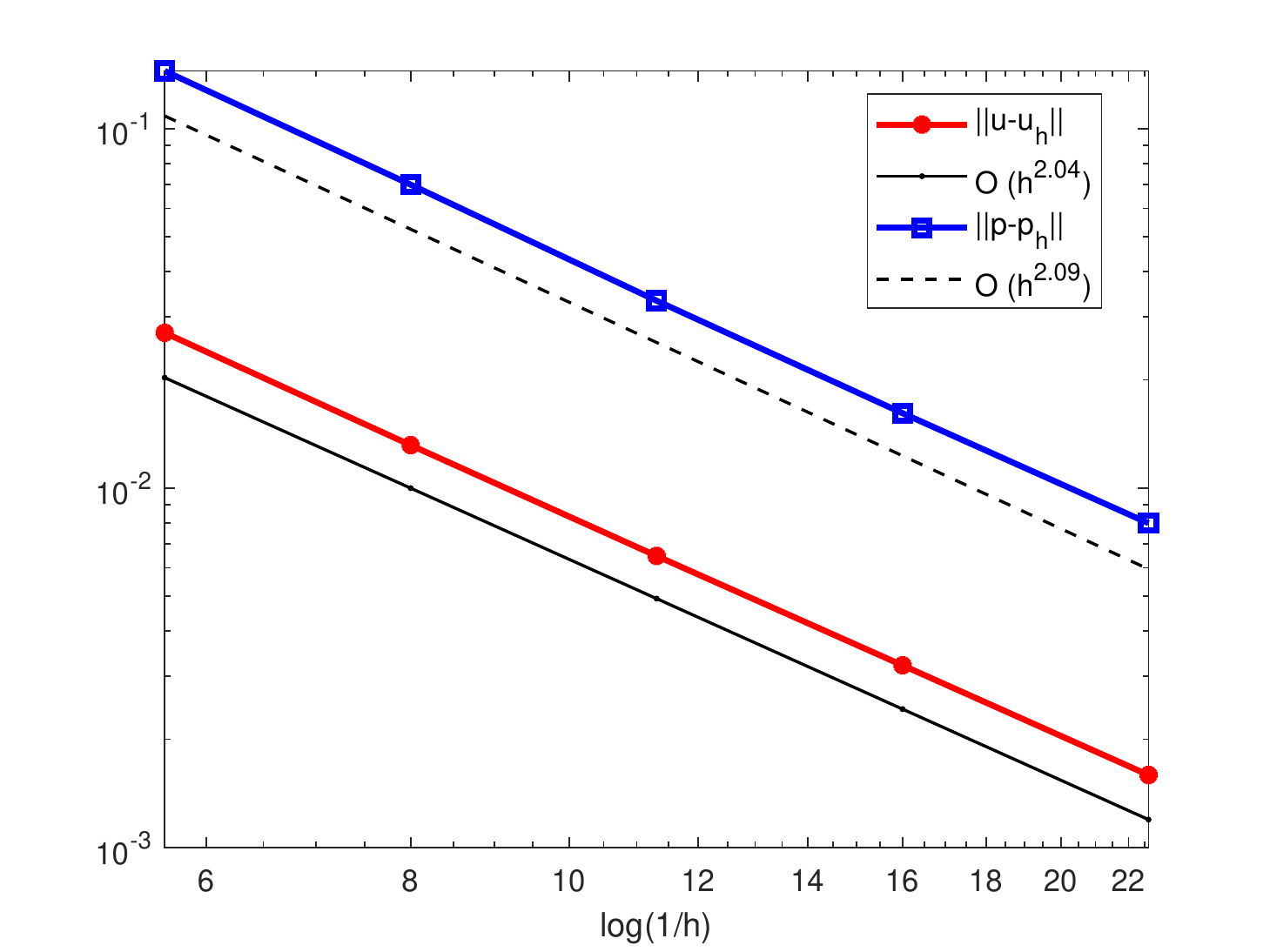}\\
  \caption{The convergence rates for the lifting VEM of the Darcy problem}\label{fig:DarcyrateLifting}
\end{figure}

We repeat the test in Example \ref{ex:Darcy}. The subroutine \mcode{Darcy\_LiftingmixedVEM.m} is used to compute the numerical solutions and the test script \mcode{main\_Darcy\_LiftingmixedVEM.m} verifies the convergence rates. The exact and numerical solutions for the variable $p$ are shown in Fig.~\ref{fig:Darcyp}. The corresponding convergence rates are displayed in Fig.~\ref{fig:DarcyrateLifting}, from which we observe the second-order convergence for both variables.

\section{Polygonal mesh generation and refinement} \label{sec:mesh}

\begin{figure}[!htb]
  \centering
  \includegraphics[scale=0.5,trim=0 80 0 80,clip]{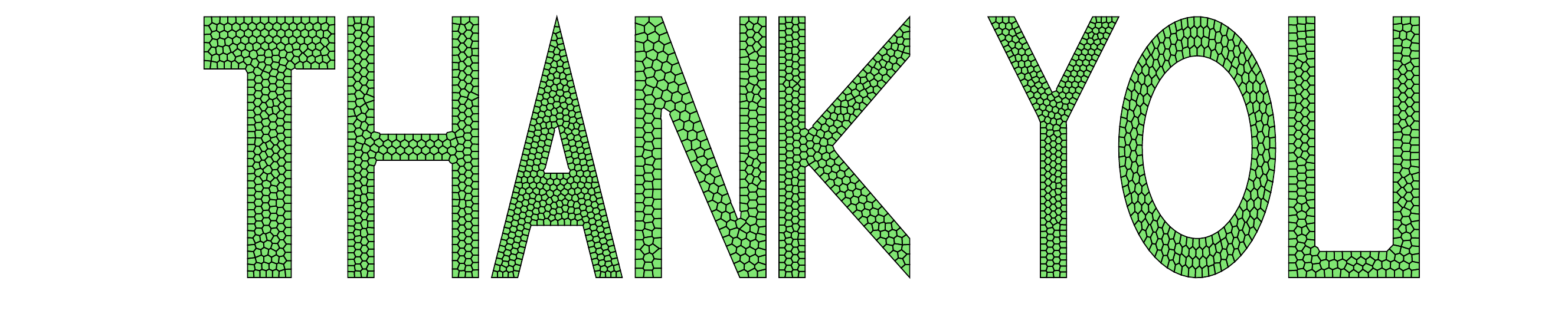}\\
  \caption{A polygonal mesh generated by PolyMesher for a complex geometry.}\label{fig:thankU}
\end{figure}

The polygonal meshes can be generated by using the MATLAB toolbox - PolyMesher introduced in \cite{Talischi-Paulino-Pereira-2012}. We provide a modified version in our package and give a very detailed description in the document. The test script is \mcode{meshfun.m} and we present a mesh for a complex geometry in Fig.~\ref{fig:thankU}.

We present some basic functions to show the polygonal meshes, including marking of the nodes, elements and (boundary) edges, see \mcode{showmesh.m}, \mcode{findnode.m}, \mcode{findelem.m} and \mcode{findedge.m}. For the convenience of the computation, some auxiliary mesh data are introduced (see Subsection \ref{subsec:datastructure}). The idea stems from the treatment of triangulation in $i$FEM, which is generalized to polygonal meshes with certain modifications. We also provide a boundary setting function to identify the Neumann and Dirichlet boundaries (See \mcode{setboundary.m}).

\begin{figure}[!htb]
  \centering
  \subfigure[]{\includegraphics[scale=0.45]{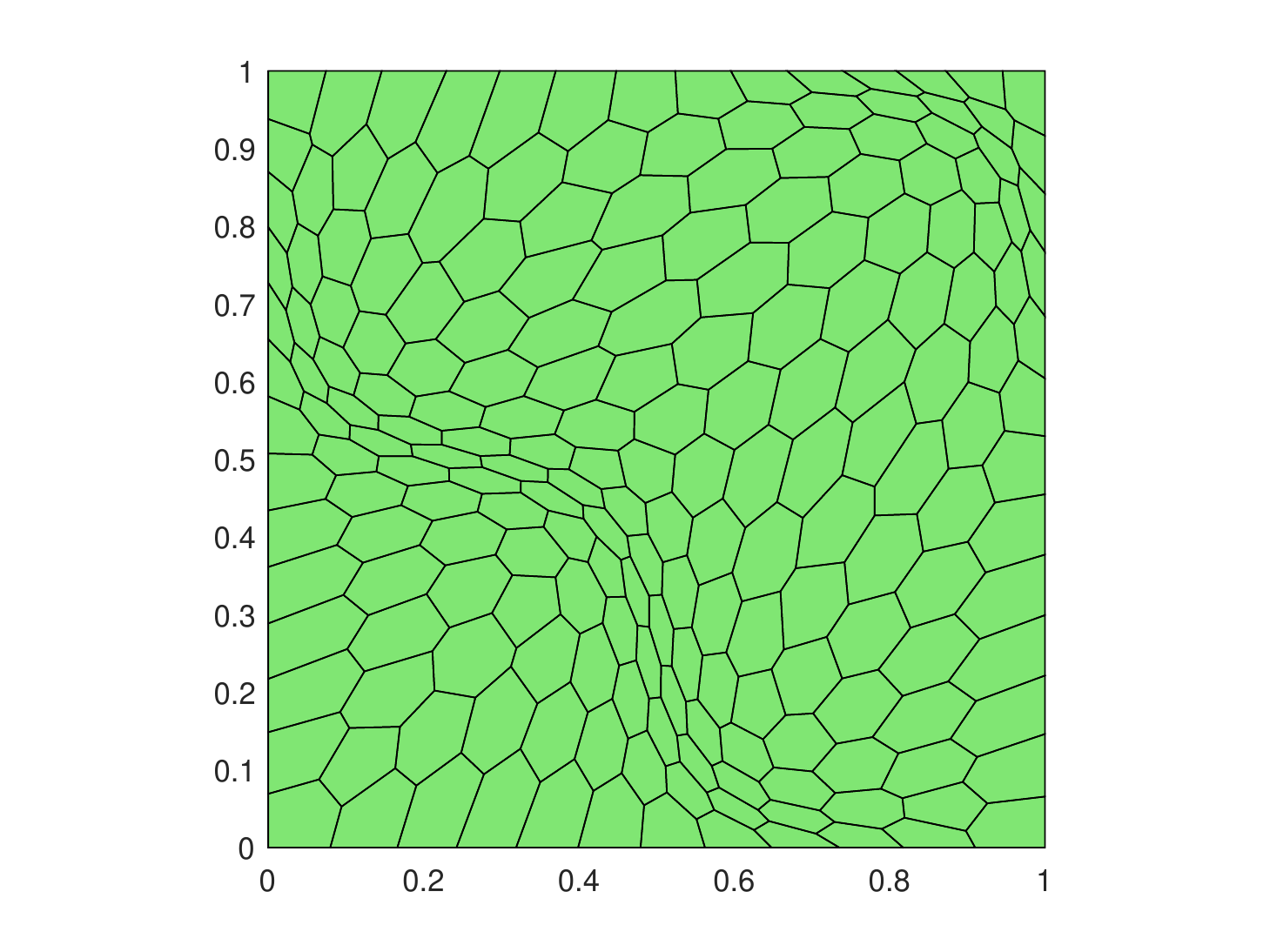}}
  \subfigure[]{\includegraphics[scale=0.45]{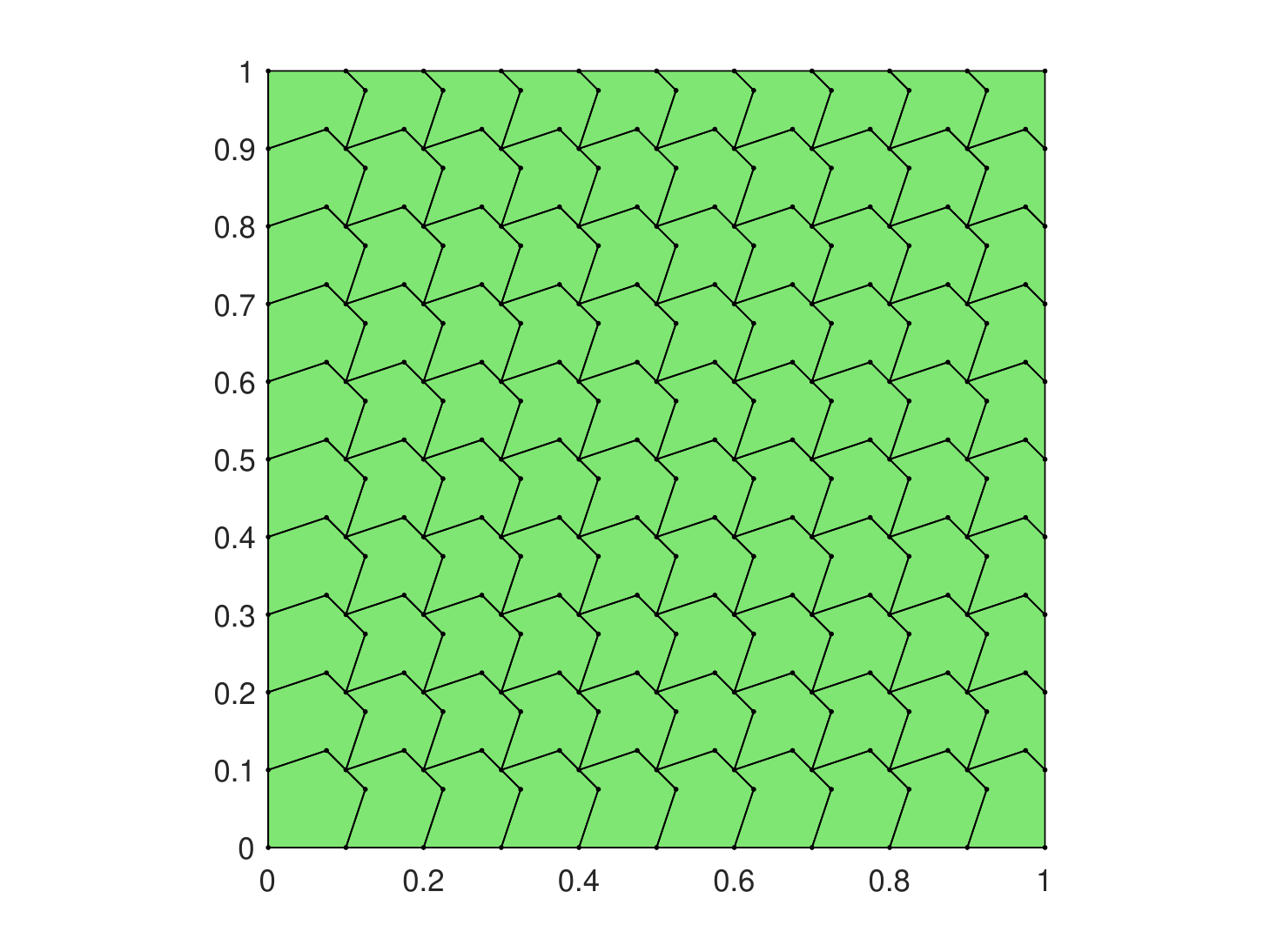}}\\
  \caption{A distorted mesh and a non-convex octagonal mesh.}\label{fig:distortion}
\end{figure}

The routine \mcode{distortionmesh.m} is used to generate a distorted mesh. Let $(\xi, \eta)$ be the coordinates on the original mesh. The nodes of the distorted mesh are obtained by the following transformation
\[
x = \xi + t_c \sin(2\pi\xi)\sin(2\pi\eta), \quad y = \eta + t_c \sin(2\pi\xi)\sin(2\pi\eta),
\]
where $(x,y)$ is the coordinate of new nodal points; $t_c$, taken as $0.1$ in the computation, is the distortion parameter. Such an example is displayed in Fig.~\ref{fig:distortion}(a). In the literature of VEMs (see \cite{Cangiani-Manzini-Sutton-2017} for example), one often finds the test for the non-convex mesh in Fig.~\ref{fig:distortion}(b), which is generated by \mcode{nonConvexMesh.m} in our package.

\begin{figure}[!htb]
  \centering
  \subfigure[]{\includegraphics[scale=0.5,trim=0 60 0 60,clip]{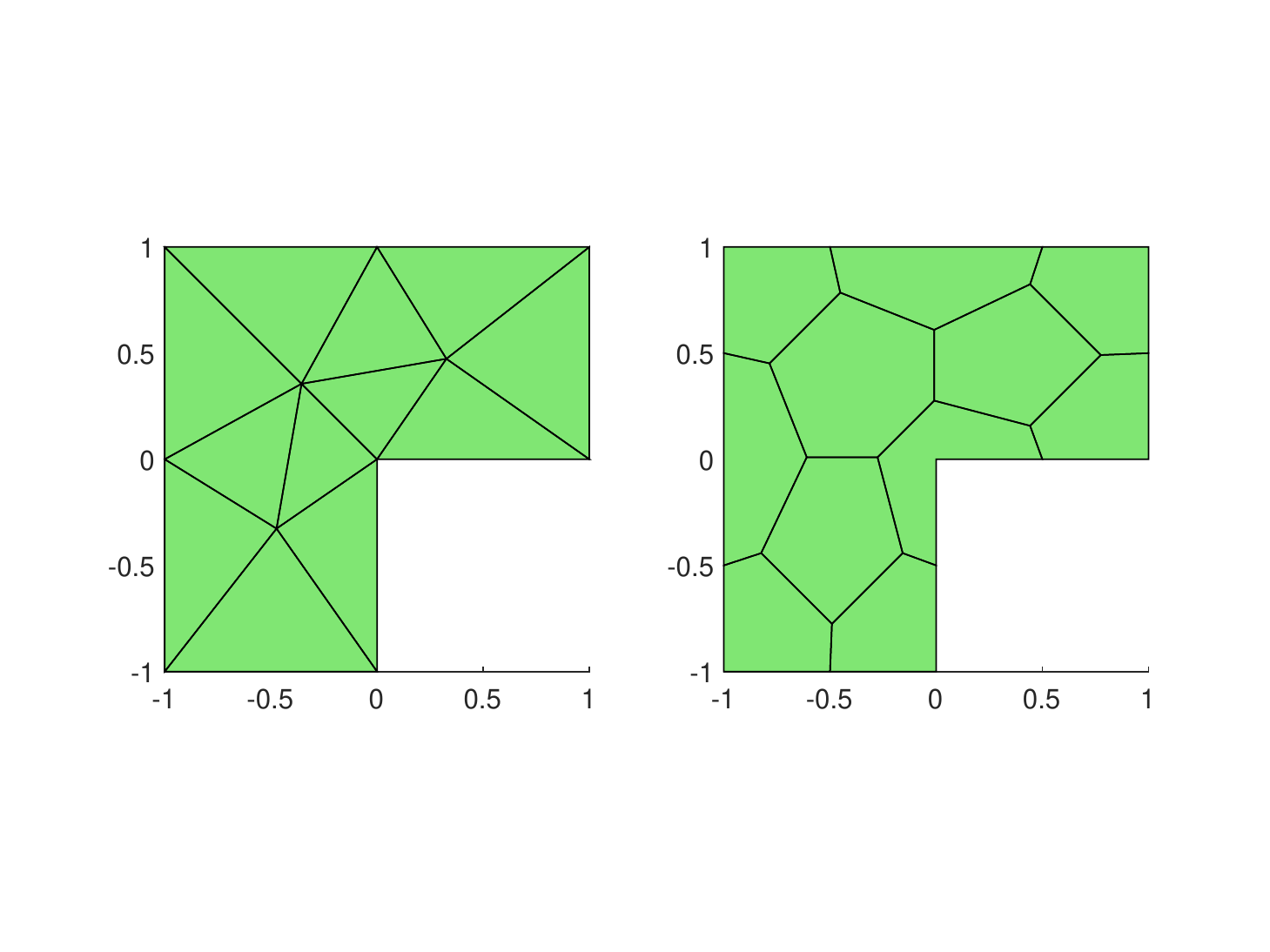}}
  \subfigure[]{\includegraphics[scale=0.5,trim=0 60 0 60,clip]{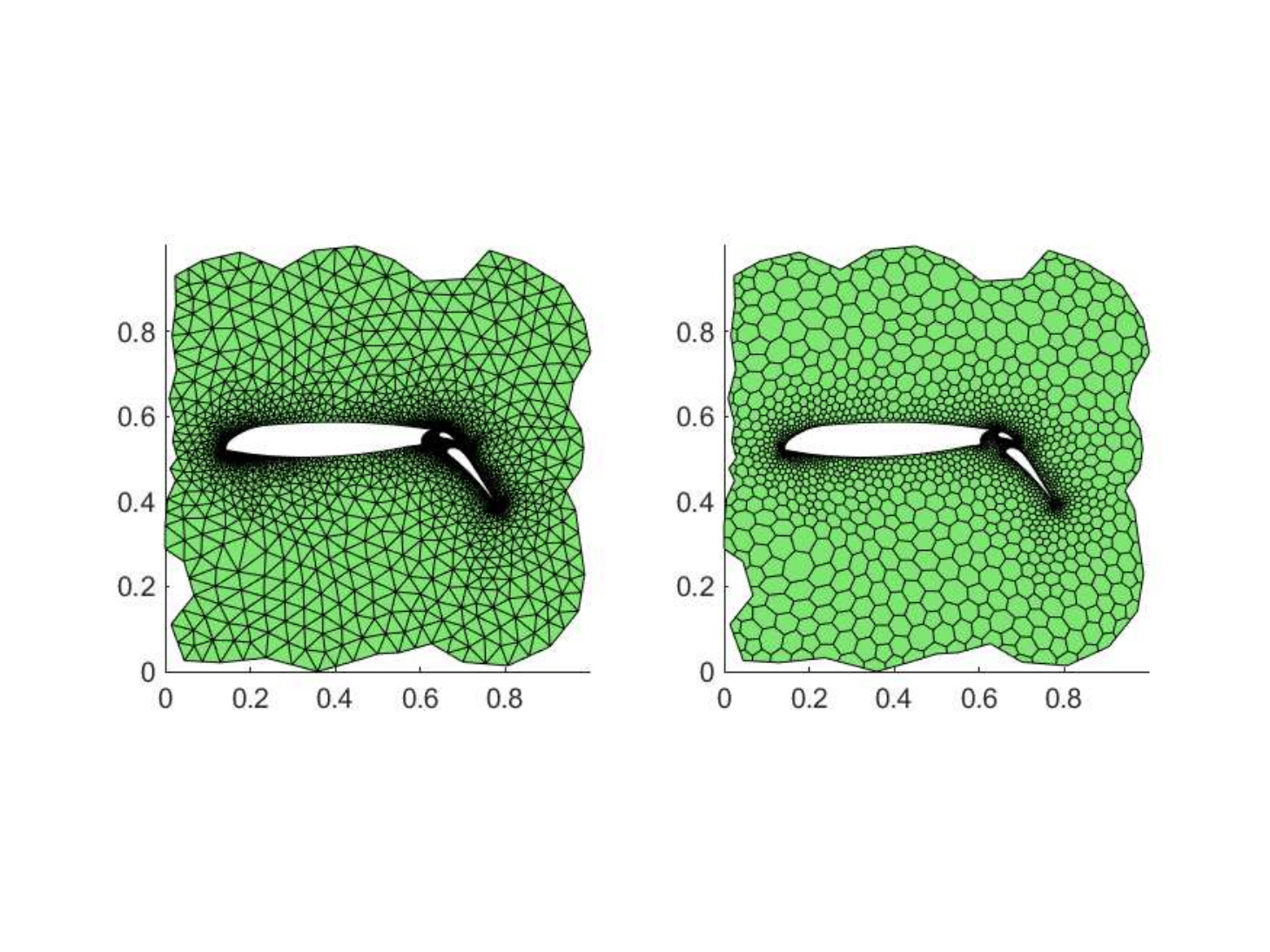}}\\
  \subfigure[]{\includegraphics[scale=0.5,trim=0 60 0 60,clip]{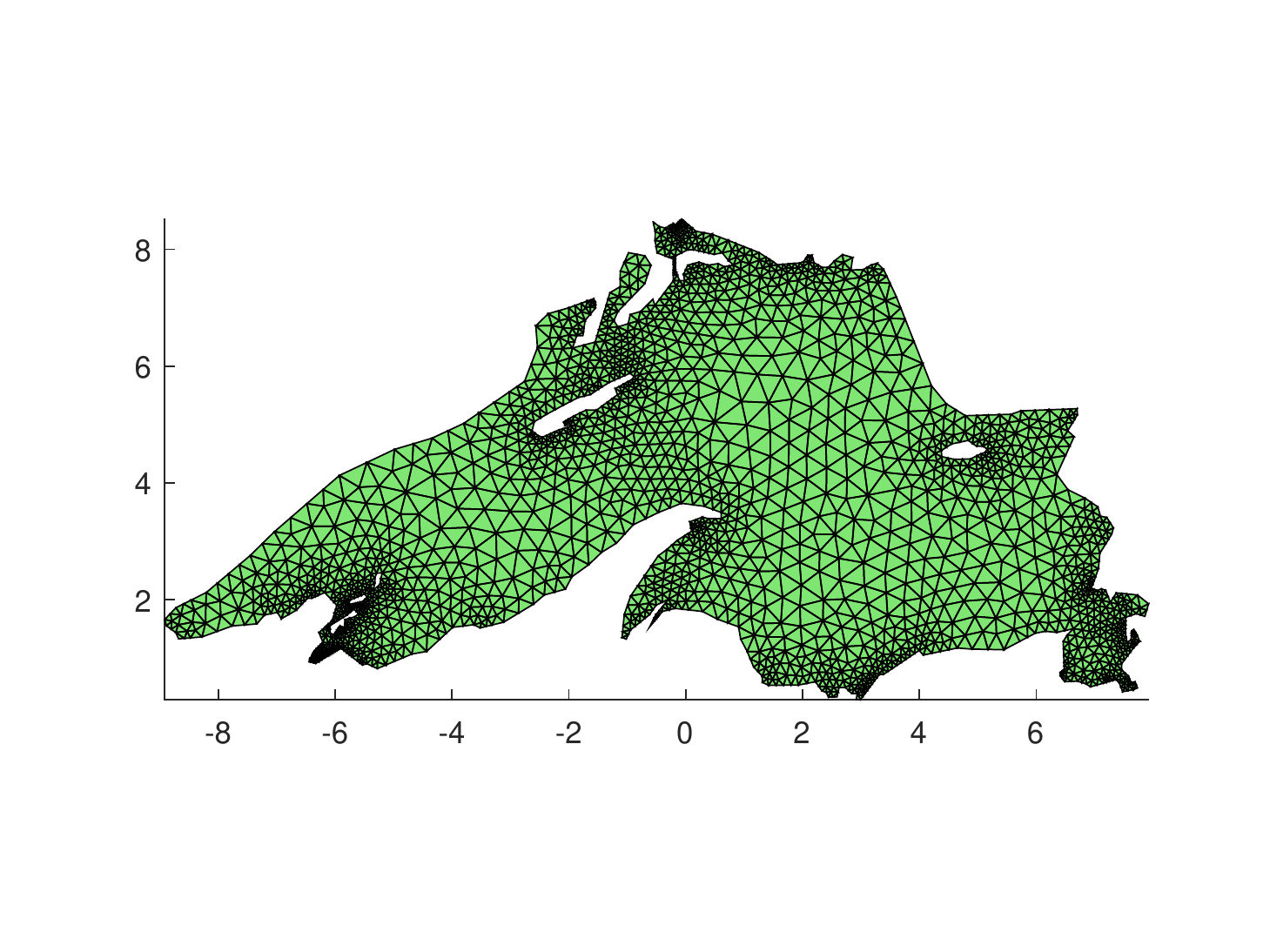}}
  \subfigure[]{\includegraphics[scale=0.5,trim=0 60 0 60,clip]{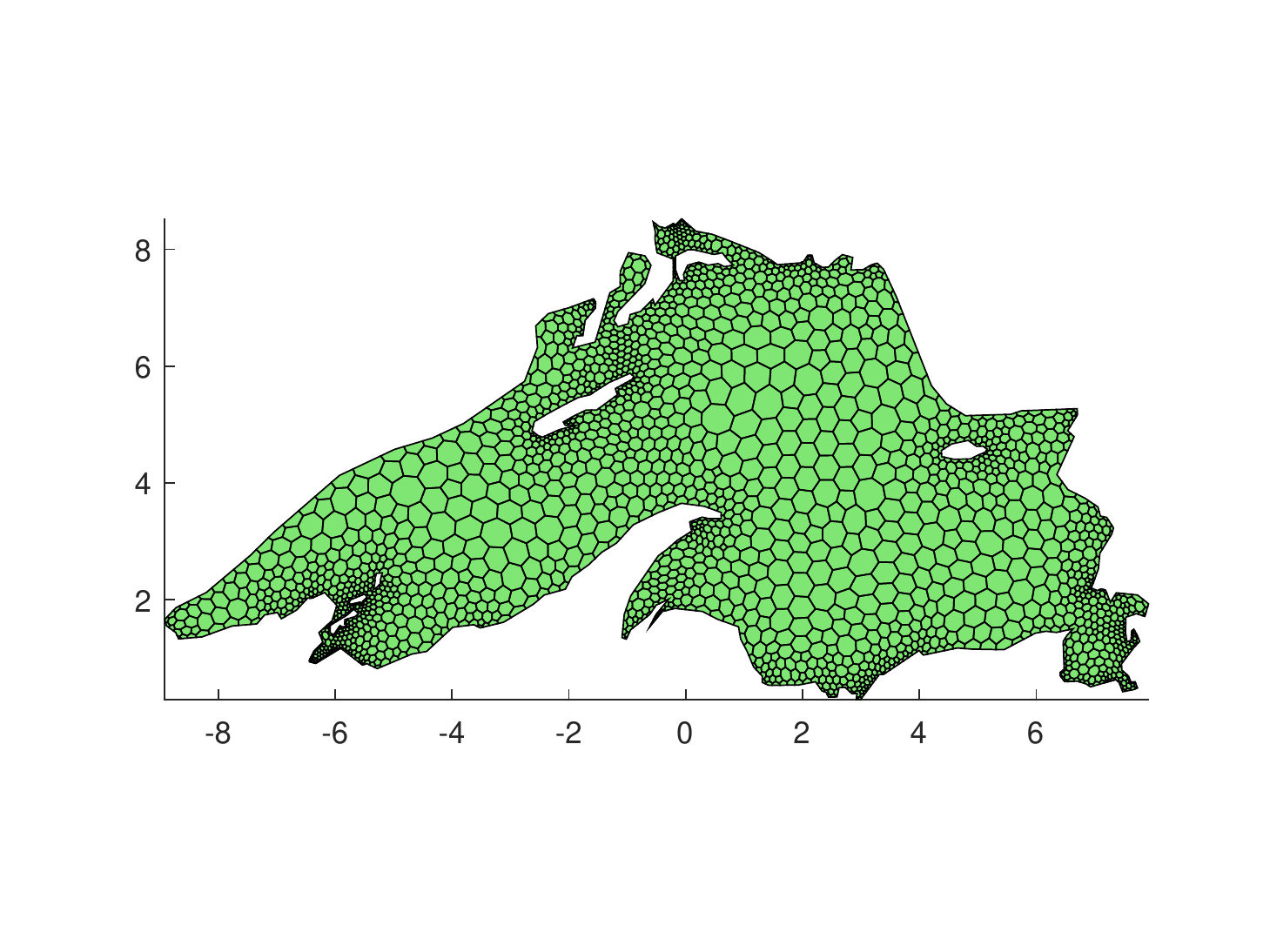}}\\
  \caption{Dual meshes generated by \mcode{dualMesh.m}.}\label{fig:dualMesh}
\end{figure}

The polygonal meshes can also be obtained from a Delaunay triangulation by establishing its dual mesh. The implementation is given in \mcode{dualMesh.m}. Several examples are given in \mcode{main\_dualMesh.m} as shown in Fig.~\ref{fig:dualMesh}.

Due to the large flexibility of the meshes, researchers have focused on the a posterior error analysis of the VEMs and made some progress in recent years \cite{Beirao-Manzini-2015,Cangiani-Georgoulis-Pryer-Sutton-2017,Berrone-Borio-2017a,Chi-Beirao-Paulino-2019,Beirao-Manzini-Mascotto-2019}.
We present an efficient implementation of the mesh refinement for polygonal meshes \cite{PolyMeshRefine} (see \mcode{PolyMeshRefine.m}). To the best of our knowledge, this is the first publicly available implementation of the polygonal mesh refinement algorithms. We divide elements by connecting the midpoint of each edge to its barycenter, which may be the most natural partition frequently used in VEM papers. To remove small edges, some additional neighboring polygons of the marked elements are included in the refinement set by requiring the one-hanging-node rule: limit the mesh to have just one hanging node per edge. The current implementation requires that the barycenter of each element is an interior point.

\section{Poisson equation} \label{sec:Poisson}

In this section we focus on the virtual element methods for the reaction-diffusion problem
\begin{equation}\label{ReactionDiffusion}
\begin{cases}
   - \Delta u + \alpha u = f \quad & \mbox{in}~~\Omega ,  \\
  u = 0 \quad & \mbox{on}~~\partial \Omega ,
\end{cases}
\end{equation}
where $\alpha$ is a nonnegative constant and $\Omega$ is a polygonal domain.

\subsection{Conforming VEMs}

Consider the first virtual element space \cite{Beirao-Brezzi-Cangiani-2013}
\begin{equation}
\label{first-VEMspace}
V_k(K): = \left\{ v \in H^1(K): \Delta v \in \mathbb{P}_{k - 2}(K)~~{\rm in}~~K,\quad  v|_{\partial K} \in \mathbb{B}_k(\partial K) \right\},
\end{equation}
where
\begin{equation*}
\mathbb{B}_k(\partial K): = \left\{ v \in C(\partial K):  v|_e \in \mathbb{P}_k(e),~~e \subset \partial K \right\}.
\end{equation*}
Equipped with the d.o.f.s in \cite{Beirao-Brezzi-Cangiani-2013}, a computable approximate bilinear form can be constructed by using the elliptic projector $\Pi_k^\nabla$. To ensure the accuracy and the well-posedness of the discrete method, a natural candidate to approximate the reaction term $(u,v)_K$ in the local bilinear form is $(\Pi _k^0u,\Pi _k^0v)_K$, where $\Pi _k^0$ is the $L^2$ projection onto $\mathbb{P}_k(K)$ and $(\cdot,\cdot)_K$ the usual $L^2(K)$ inner product. Nevertheless, $\Pi _k^0$ can not be computed in terms of d.o.f.s attached to $V_k(K)$. Therefore, Ahmad \emph{et al.}~\cite{Ahmad-Alsaedi-Brezzi-2013} modified the VEM space \eqref{first-VEMspace} to a local enhancement space
\begin{equation}
\label{second-VEMspace}
W_k(K) = \left\{ w \in \widetilde{V}_k(K): (w - \Pi _k^\nabla w,q)_K = 0,~~ q \in \mathbb{M}_k(K) \backslash \mathbb{M}_{k - 2}(K) \right\},
\end{equation}
where the lifting space is
\begin{equation*}
\widetilde{V}_k(K):= \left\{ v \in H^1(K):  v|_{\partial K} \in \mathbb{B}_k(\partial K),~~ \Delta v \in \mathbb{P}_k(K)~~\mbox{in}~K \right\}.
\end{equation*}
In this space, the operator $\Pi _k^0$ can be easily computed using $\Pi _k^\nabla $ and the local d.o.f.s related to $W_k(K)$:
\begin{itemize}
\item $\chi_a$: the values at the vertices of $K$,
\[\chi_i(v) = v(a_i),\quad a_i \mbox{ is a vertex of }K.\]
\item $\chi_e^{k-2}$: the moments on edges up to degree $k-2$,
\[\chi_e (v) = |e|^{-1}(m_e, v)_{e},\quad m_e\in \mathbb  M_{k-2}(e), ~~e\subset \partial K.\]
\item $\chi_K^{k-2}$: the moments on element $K$ up to degree $k-2$,
\[\chi_{K}(v) = |K|^{-1}(m_K, v)_K ,\quad m_K\in \mathbb  M_{k-2}(K).\]
\end{itemize}

In what follows, denote by $W_{k,h}$ a finite dimensional subspace of $V$, produced by combining all $W_k(K)$ for $K\in \mathcal{T}_h$ in a standard way. We define a local $H^1$-elliptic projection operator $\Pi_k^\nabla:{H^1}(K)\to\mathbb{P}_k(K)$ by the relations
\begin{equation}\label{projector}
(\nabla \Pi_k^\nabla v,\nabla p)_K=(\nabla v,\nabla p)_K, \quad  p\in \mathbb{P}_k(K).
\end{equation}
Since $(\nabla \cdot, \nabla \cdot)$ is only semi-positive definite, the constraints
\begin{align*}
\int_{\partial K} v {\rm d} s=\int_{\partial K} \Pi_{1}^\nabla  v {\rm d} s, \quad  k = 1; \qquad
\int_{K} v {\rm d} x=\int_{K} \Pi_{k}^\nabla  v {\rm d} x,  \quad  k\geq 2
\end{align*}
should be imposed.
The VEM is to find $u_h \in W_{k,h}$ such that
 \begin{equation}\label{apprxobil}
 {a_h^c}(u_h,v_h) = \langle f_h,v_h\rangle ,\quad v_h \in W_{k,h},
 \end{equation}
 where
\begin{equation*}
{a_h^c}(u_h,v_h):= \sum\limits_{K \in \mathcal{T}_h} {a_h^{c,K}(u_h,v_h)}
\end{equation*}
and $\langle \cdot,\cdot\rangle$ is the duality pairing between $W_{k,h}'$ and $W_{k,h}$. The local approximate bilinear form is
\begin{equation*}
a_h^{c,K}(v,w) = a^K(\Pi _k^\nabla v, \Pi _k^\nabla w) + \alpha (\Pi _k^0v,\Pi _k^0w)_K + S^K(v - \Pi _k^\nabla v,w - \Pi _k^\nabla w)
\end{equation*}
where $a^K(v,w) = (\nabla v, \nabla w)_K$ and the stabilization term is realized as
\begin{equation*}
S^K(v- \Pi_k^\nabla v, w- \Pi_k^\nabla w) := (1+\alpha h_K^2)\boldsymbol  \chi (v- \Pi_k^\nabla v)\cdot \boldsymbol  \chi (w-\Pi_k^\nabla w),
\end{equation*}
where $\boldsymbol \chi$ is referred to as the d.o.f vector of a VEM function on $K$.

Let $\phi^T  = (\phi_1,\phi_2, \cdots ,\phi_{N_k})$ be the basis of $W_k(K)$ and let $m^T = (m_1,m_2, \cdots ,m_{N_m})$ be the scaled monomials. Then the transition matrix $D$ can be defined through $m^T = \phi^TD$. The matrix expressions of $\Pi_k^\nabla \phi^T$ in the basis $\phi^T$ and $m^T$ are defined by
\[\Pi_k^\nabla \phi^T = \phi^T\boldsymbol\Pi_k^\nabla, \quad  \Pi_k^\nabla \phi^T = m^T\boldsymbol\Pi_{k*}^\nabla .\]
One easily finds that $\boldsymbol\Pi_k^\nabla  = D\boldsymbol\Pi_{k*}^\nabla$. In vector form, the definition of \eqref{projector} can be rewritten as
\[
\begin{cases}
a^K(m , \Pi_k^\nabla \phi^T) = a^K(m, \phi^T ), \\
P_0(\Pi_k^\nabla \phi^T) = P_0( \phi^T),
\end{cases}
\]
where
\[P_0(v) = \int_{\partial K} v {\rm d}s, \quad k=1;  \qquad  P_0(v) = \int_K v {\rm d}x, \quad k\ge 2.\]
Define
\begin{equation}\label{GB}
G = a^K(m , m^T), \quad B = a^K(m,\phi^T).
\end{equation}
Then
\[\begin{cases}
  G\boldsymbol\Pi_{k*}^\nabla  = B, \\
P_0( m^T)\boldsymbol\Pi_{k*}^\nabla  = P_0( \phi^T)
\end{cases} \quad \mbox{or} \quad
\tilde G\boldsymbol\Pi_{k*}^\nabla  = \tilde B.\]
Note that one always has the consistency relations $G = BD$ and $\tilde G = \tilde BD$.

We remark that in the implementation of the VEM, the most involved step is to compute the boundary part of $B$ resulted from the integration by parts:
\begin{equation}\label{Bibp}
B = \int_K \nabla m  \cdot \nabla \phi^T\mathrm{d}x
 =  - \int_K \Delta m   \cdot \phi^T\mathrm{d}x + \sum\limits_{e \subset \partial K} \int_e (\nabla  m  \cdot {{\boldsymbol n}_e})\phi^T \mathrm{d}s.
\end{equation}
For the conforming VEMs, however, we can compute the second term of \eqref{Bibp} by using the assembling technique for finite element methods since the VEM function restriction to the boundary $\partial K$ is piecewise polynomial. For example, the computation for $k=2$ reads
\vspace{-0.8cm}
\begin{lstlisting}
% second term
I2 = zeros(Nm, Ndof);
v1 = 1:Nv; v2 = [2:Nv,1];
elem1 = [v1(:), v2(:), v1(:)+Nv];
for im = 1:Nm
    gradmc = Gradmc{im};
    F1 = 1/6*sum(gradmc(x(v1), y(v1)).*Ne, 2);
    F2 = 1/6*sum(gradmc(x(v2), y(v2)).*Ne, 2);
    F3 = 4/6*sum(gradmc(xe, ye).*Ne, 2);
    F = [F1, F2, F3];
    I2(im, :) = accumarray(elem1(:), F(:), [Ndof 1]);
end
\end{lstlisting}

\begin{example}\label{ex:ReactionDiffusion}
Let $\alpha=1$ and $\Omega  = (0,1)^2$. The Neumann boundary condition is imposed on $x=0$ and $x=1$ with the exact solution given by $u(x,y)=\sin(2x+0.5)\cos(y+0.3)+\log(1+xy)$.
\end{example}

\begin{figure}[!htb]
  \centering
  \subfigure[$k=1$]{\includegraphics[scale=0.35]{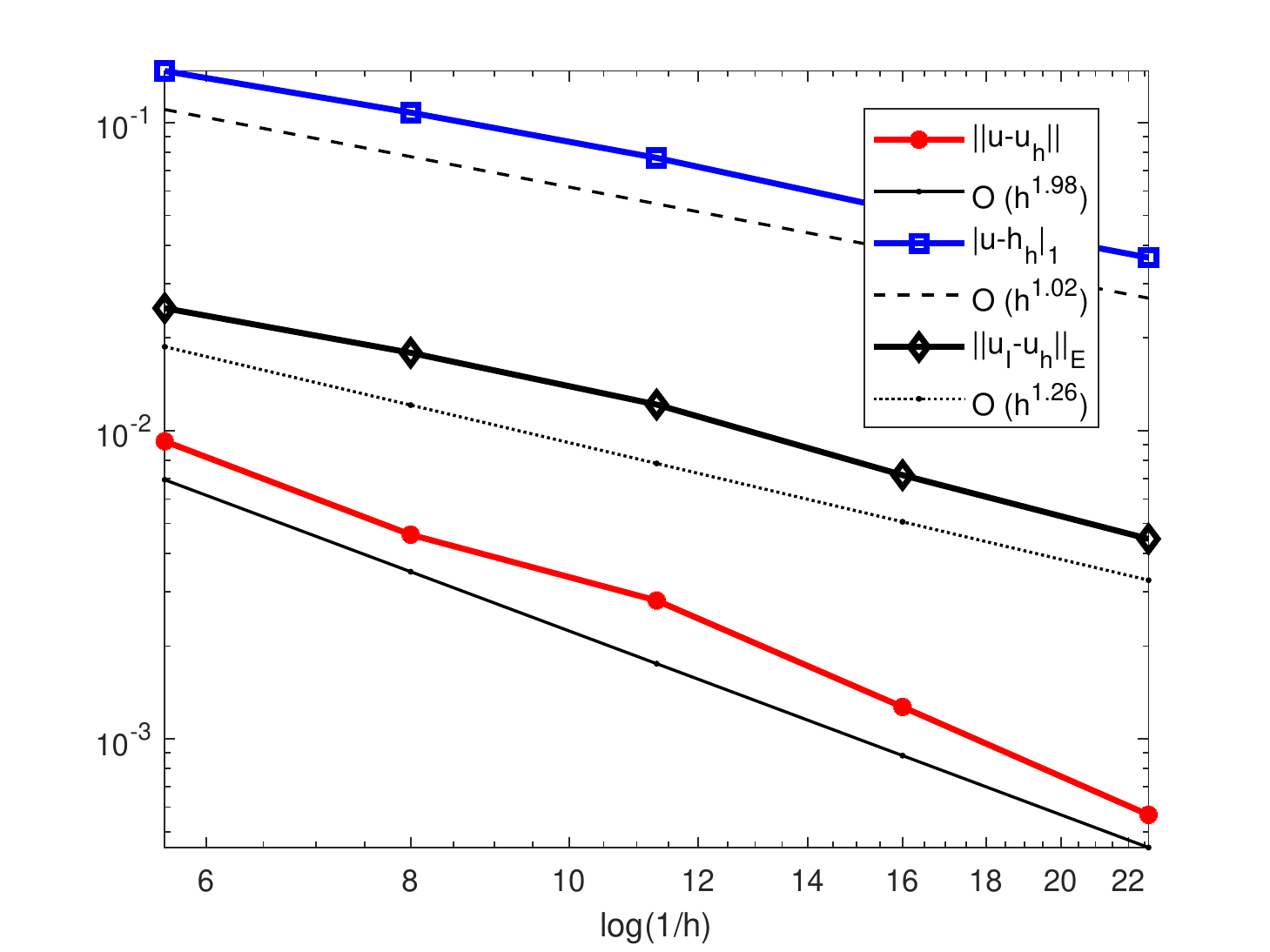}}
  \subfigure[$k=2$]{\includegraphics[scale=0.35]{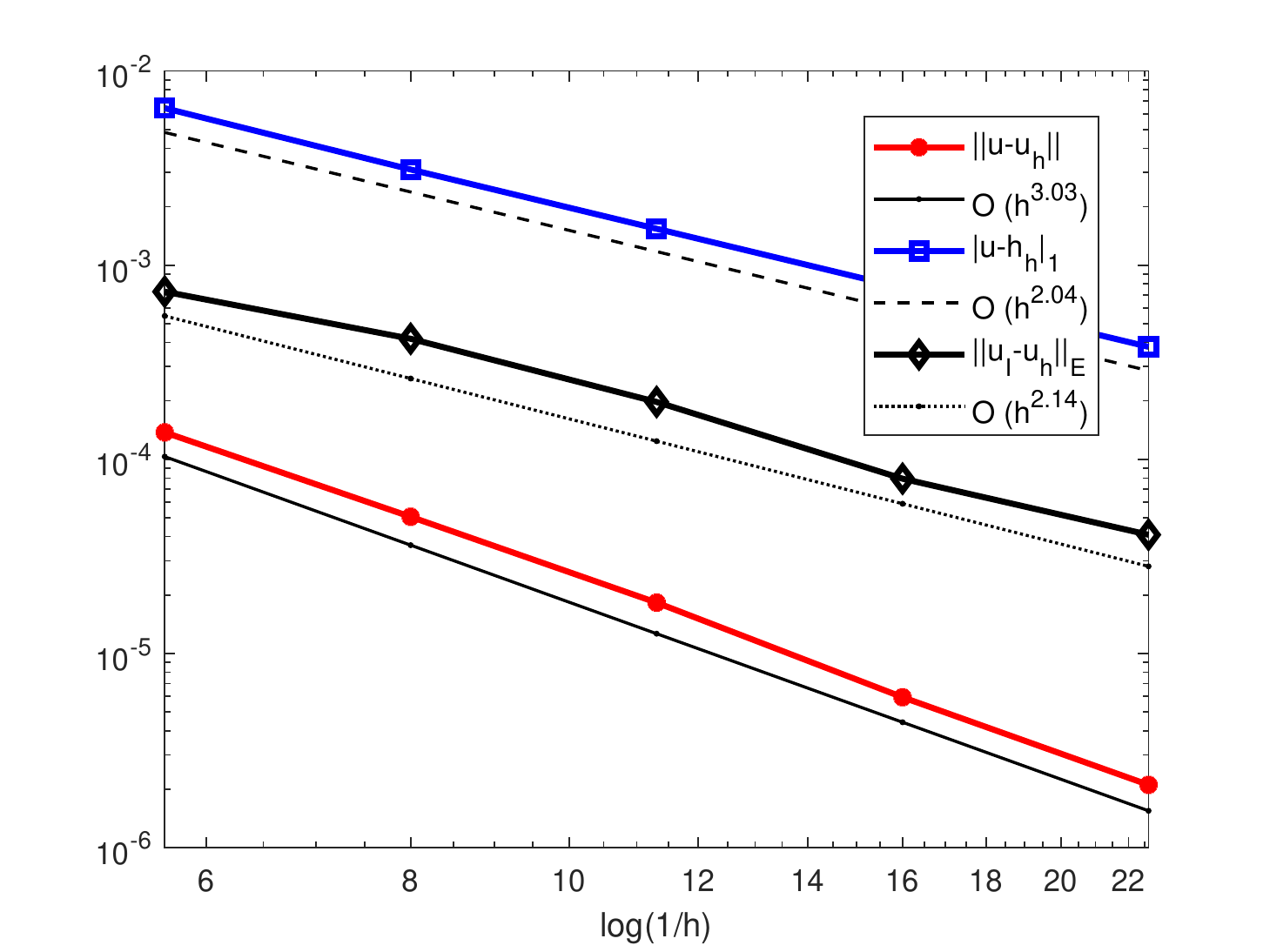}}
  \subfigure[$k=3$]{\includegraphics[scale=0.35]{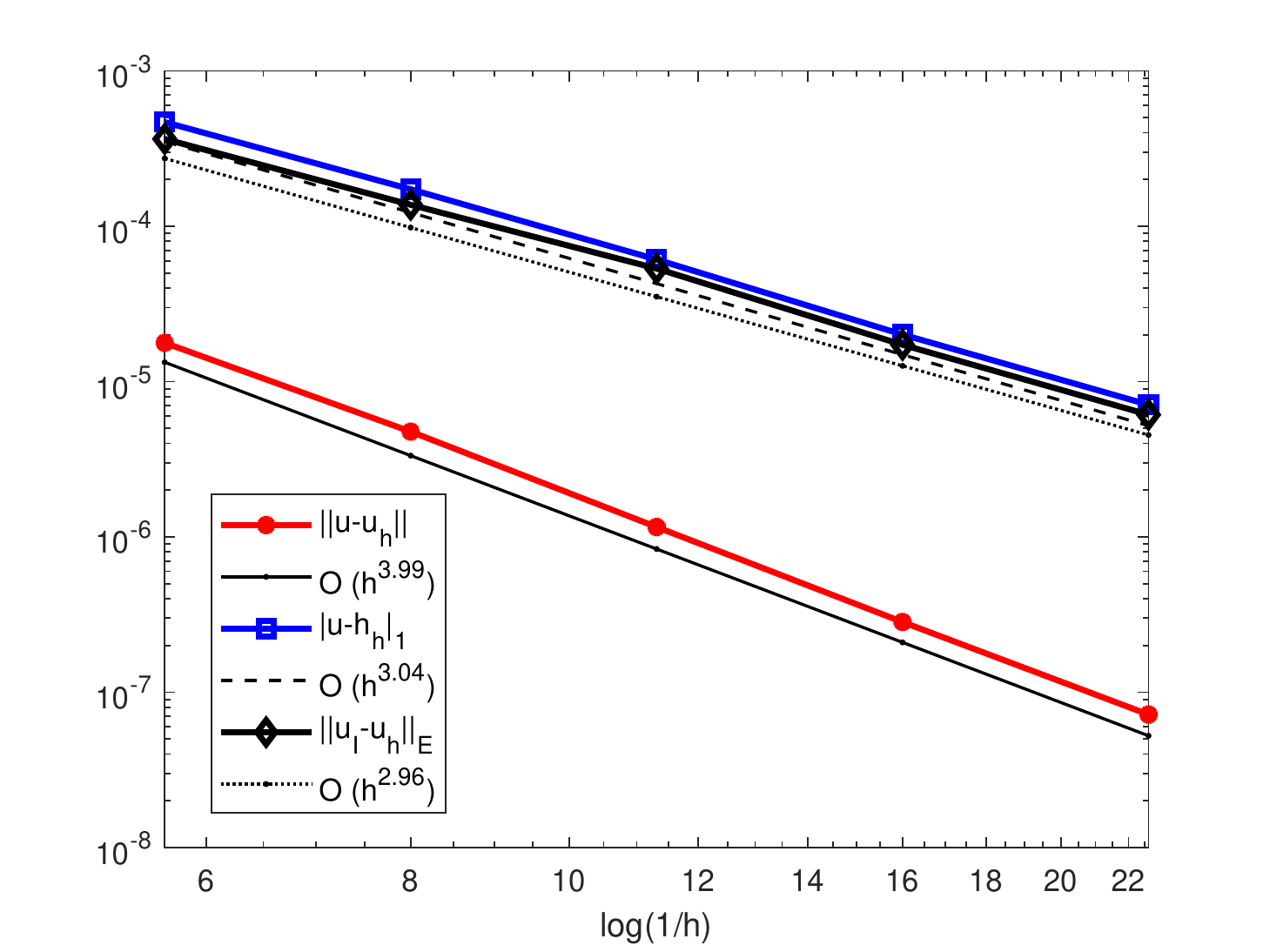}}\\
  \caption{Convergence rates for Example \ref{ex:ReactionDiffusion}.}\label{fig:ReactionDiffusion}
\end{figure}

The test script is \mcode{main\_PoissonVEMk3.m} for $k=3$. The optimal rate of convergence of the $H^1$-norm (3rd order), $L^2$-norm (4th order) and energy norm (3rd order) is observed for $k = 3$ from Fig.~\ref{fig:ReactionDiffusion}.

\begin{example}\label{ex:circle}
 In this example, the exact solution is $u(x,y) = y^2\sin (\pi x)$. The domain $\Omega$ is taken as a unit disk.
\end{example}

\begin{figure}[!htb]
  \centering
  \includegraphics[scale=0.3]{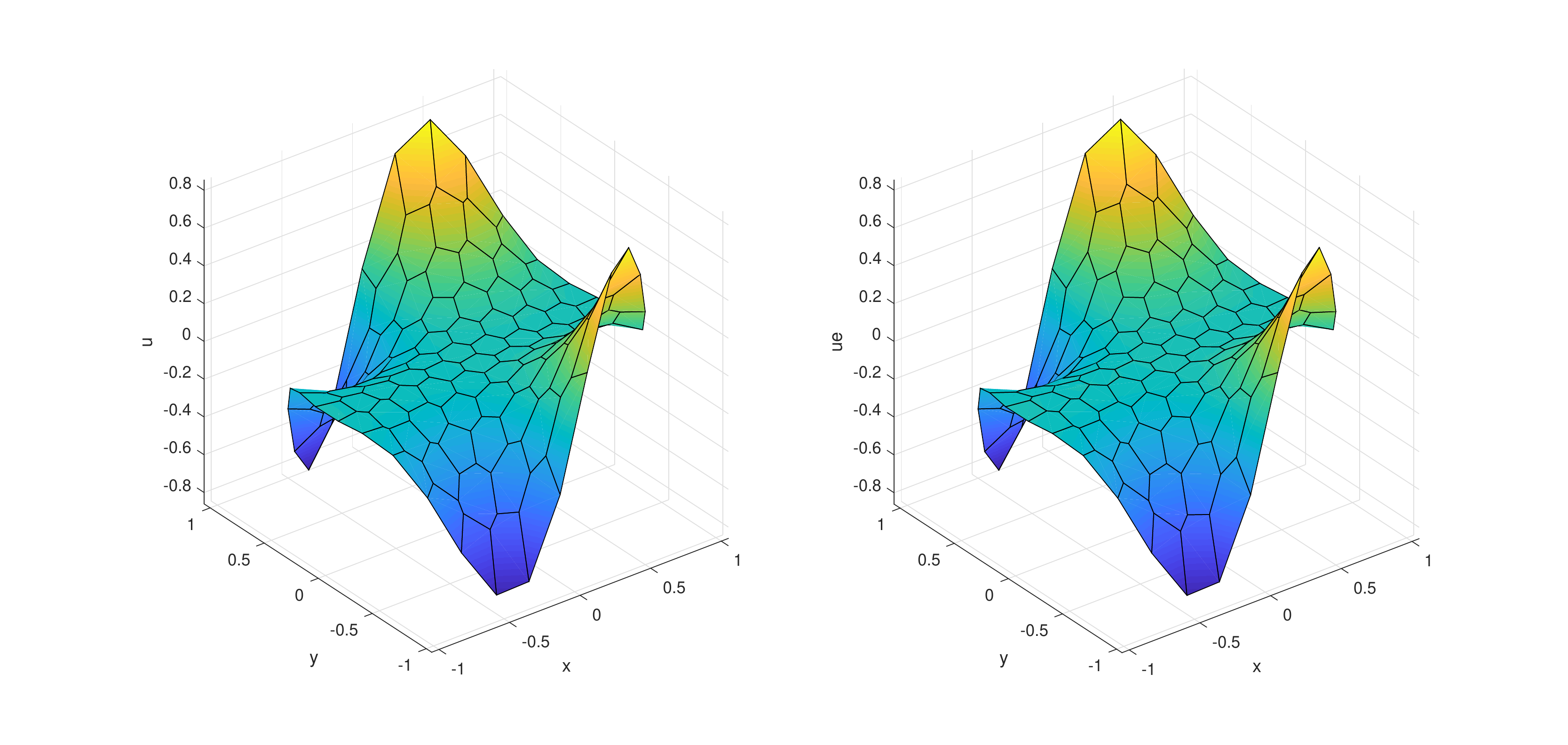}\\
  \caption{Numerical and exact solutions for Example \ref{ex:circle} ($k=1$)}\label{fig:solutionCircle}
\end{figure}

The Neumann boundary condition is imposed on the boundary of the upper semicircle. The nodal values are displayed in Fig.~\ref{fig:solutionCircle}.

\subsection{Nonconforming VEMs}

\subsubsection{The standard treatment of the domain boundary}

Now we consider the nonconforming VEMs proposed in \cite{Ayuso-Lipnikov-Manzini-2016} for the problem \eqref{ReactionDiffusion} with $\alpha=0$. The local nonconforming virtual element space is defined by
\[V_k^{nc}(K) = \left\{ v \in H^1(K): \Delta v \in \mathbb{P}_{k - 2}(K),~~\partial_{\boldsymbol n}v|_e \in \mathbb{P}_{k - 1}(e),~e \subset \partial K \right\}, \quad k\ge 1,\]
and the d.o.f.s can be chosen as
\begin{itemize}
  \item $\boldsymbol\chi_{\partial K}(v)$: the moments on edges up to degree $k-1$,
  \[\chi_e(v) =|e|^{-1} (m_e,v)_e,\quad m_e \in \mathbb{M}_{k - 1}(e),~~e\subset \partial K.\]
  \item $\boldsymbol\chi_K(v)$: the moments on $K$ up to degree $k-2$,
  \[\chi_K(v) =|K|^{-1}(m_K,v)_K,\quad m_K \in \mathbb{M}_{k - 2}(K).\]
  \end{itemize}

The elliptic projection and the discrete problem can be constructed in the same way as before. We briefly discuss the implementation of $k=1$. In this case, the first term for the matrix $B$ in \eqref{Bibp} vanishes and the second term is obtained from the Kronecher's property $\chi_i(\phi_j) = \sigma_{ij}$, which is realized as follows.
\vspace{-0.8cm}
\begin{lstlisting}
% B
B = zeros(Nm,Nv);
for i = 1:Nv % loop of edges
    gi = sum(Gradm.*repmat(Ne(i,:),3,1), 2); % gradm*Ne
    B(:,i) = gi;
end
\end{lstlisting}

For completeness, we also present the application of the Neumann boundary conditions. Let $e$ be a boundary edge of the Neumann boundary $\Gamma_N$. Then the local boundary term is approximated as
\[F_e = \int_e \Pi_{0,e}^0 g_N \phi_e \mathrm{d}s,\]
where $g_N = \partial_{\boldsymbol n}u = \nabla u \cdot \bb{n}_e$,  $\Pi_{0,e}^0$ is the $L^2$ projection on $e$, and $\phi_e$ is the local basis function associated with the d.o.f $\chi_e(v) = |e|^{-1}(m_e,v)_e$. For problems with known explicit solutions, we provide the gradient $g_N = \nabla u$ in the PDE data instead and compute the true $g_N$ in the M-file. Obviously, the $L^2$ projection $\Pi_{0,e}^0 g_N$ can be computed by the mid-point formula or the trapezoidal rule. The computation reads
\vspace{-0.8cm}
\begin{lstlisting}
%% Assemble Neumann boundary conditions
bdEdgeN = bdStruct.bdEdgeN;  bdEdgeIdxN = bdStruct.bdEdgeIdxN;
if ~isempty(bdEdgeN)
    g_N = pde.Du;
    z1 = node(bdEdgeN(:,1),:); z2 = node(bdEdgeN(:,2),:);
    e = z1-z2;  % e = z2-z1
    Ne = [-e(:,2),e(:,1)]; % scaled ne
    F1 = sum(Ne.*g_N(z1),2);
    F2 = sum(Ne.*g_N(z2),2);
    FN = (F1+F2)/2;
    ff = ff + accumarray(bdEdgeIdxN(:), FN(:),[NNdof 1]);
end
\end{lstlisting}
Note that the structure data \mcode{bdStruct} stores all necessary information of boundary edges, which is obtained by using the subroutine \mcode{setboundary.m}. In the above code, \mcode{bdEdgeN} gives the data structure \mcode{bdEdge} in Subsection \ref{subsec:datastructure} for Neumann edges and \mcode{bdEdgeIdxN} provides their indices in the data structure \mcode{edge}.

\begin{figure}[!htb]
  \centering
  \includegraphics[scale=0.45]{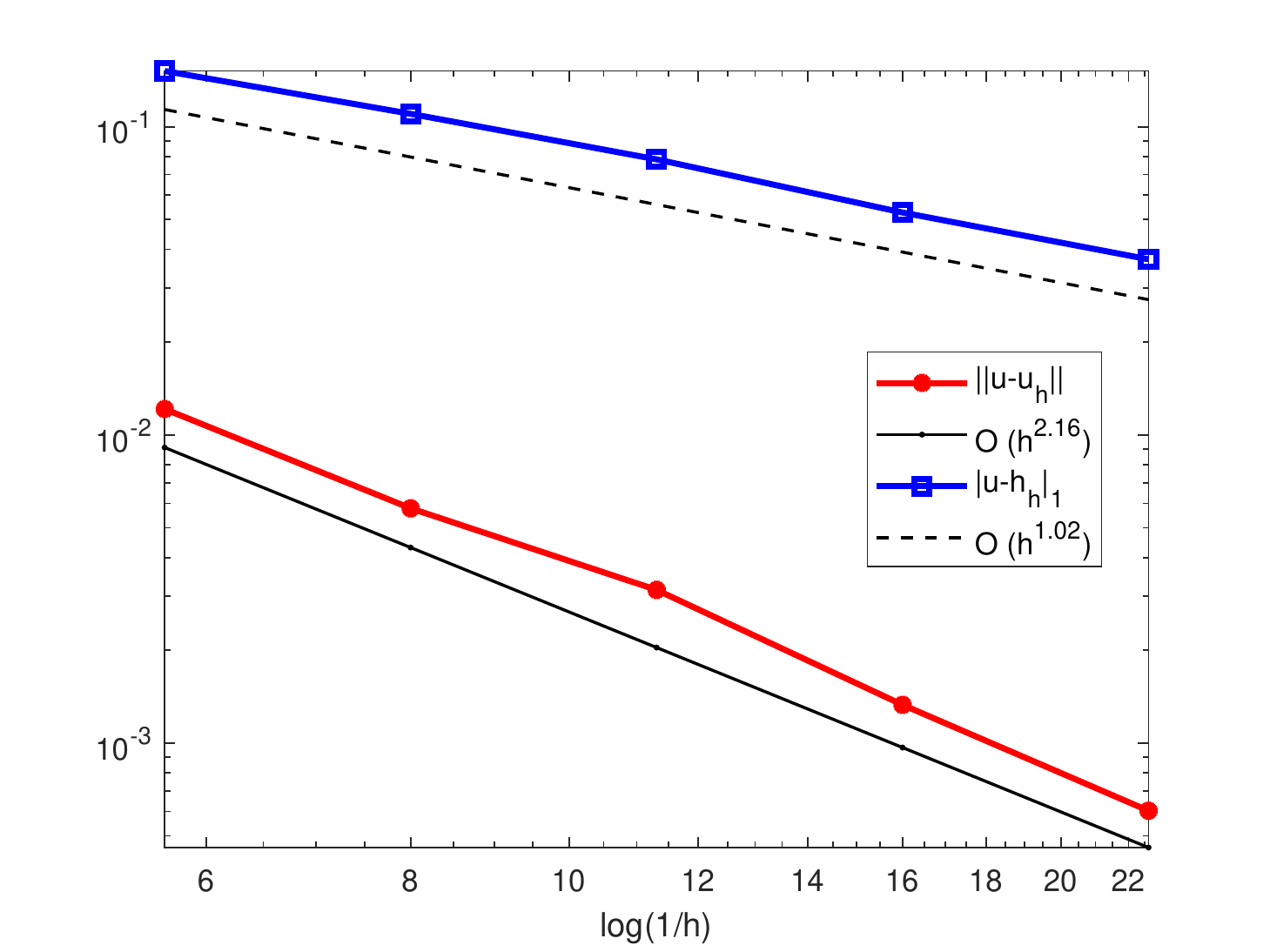}\\
  \caption{Convergence rates of the Poisson equation for the nonconforming VEM ($k=1$)}\label{fig:PoissonNonconformingrate}
\end{figure}

We still consider Example \ref{ex:ReactionDiffusion} and display the convergence rates in the discrete $L^2$ and $H^1$ norms in Fig.~\ref{fig:PoissonNonconformingrate}. The optimal rate of convergence is observed for both norms and the test script is \mcode{main\_PoissonVEM\_NC.m}.

\subsubsection{A continuous treatment of the domain boundary}

In some cases it may be not convenient to deal with the moments on the domain boundary $\partial \Omega$. To do so, we can modify the corresponding local virtual element spaces to the following one
\begin{align*}
V_k^{ncb}(K)
& = \Big\{ v \in H^1(K): \Delta v \in \mathbb{P}_{k - 2}(K),~~\partial_{\boldsymbol n}v|_e \in \mathbb{P}_{k - 1}(e),~~e ~ \mbox{is an interior edge}  \\
& \hspace{3cm} v|_e \in \mathbb{P}_k(e), ~~ e ~\mbox{ is a boudary edge}, \quad e\subset \partial K \Big \}.
\end{align*}

We repeat the test in Fig.~\ref{fig:PoissonNonconformingrate} and the discrete $L^2$ and $H^1$ errors are listed in Tab.~\ref{tab:modifiedNcVEM}. Please refer to \mcode{main\_PoissonVEM\_NCb.m} for the test script.

\begin{table}[H]
 \centering
 \caption{The discrete errors for the modified nonconforming VEM} \label{tab:modifiedNcVEM}
\begin{tabular}{rcccccccccc}
  \hline
  $\sharp {\rm Dof}$   &     $h$    &   ErrL2   &  ErrH1 \\
  \hline
  97 &  1.768e-01 &  1.17632e-02  & 1.51836e-01 \\
 193 &  1.250e-01 &  5.37635e-03  &  1.10077e-01 \\
 383 &  8.839e-02 &  3.14973e-03  & 7.82566e-02 \\
 760 &  6.250e-02 &  1.43956e-03  & 5.26535e-02 \\
1522 &  4.419e-02 &  6.48454e-04  &  3.71405e-02 \\
  \hline
\end{tabular}
\end{table}

\section{Linear elasticity problems}

The linear elasticity problem is
\begin{equation}\label{linearElasticity}
\begin{cases}
   - {\rm div} \bb{\sigma } = \bb{f} \quad & \text{in}~~~\Omega ,  \\
  \bb{u} = {\bb{0}} \quad & \text{on}~~~\Gamma_0,  \\
  \bb{\sigma n} = \bb{g} \quad & \text{on}~~~\Gamma_1,
\end{cases}
\end{equation}
where $\bb{n} = (n_1,n_2)^T$ denotes the outer unit vector normal to $\partial\Omega$. The constitutive relation for linear elasticity is
\[\bb{\sigma}(\bb{u}) = 2\mu\bb{\varepsilon}(\bb{u}) + \lambda({\rm div} \bb{u})\bb{I},\]
where $\bb{\sigma} = (\sigma_{ij})$ and $\bb{\varepsilon} = (\varepsilon_{ij})$ are the second order stress and strain tensors, respectively, satisfying $\varepsilon_{ij} = \frac{1}{2}(\partial_i u_j + \partial_j u_i)$, $\lambda$ and $\mu$ are the Lam\'{e} constants, $\bb{I}$ is the identity matrix, and ${\rm div}\bb{u} = \partial_1 u_1 + \partial_2 u_2$.

\subsection{Conforming VEMs}

\subsubsection{The displacement type}\label{subsec:displacement}

The equilibrium equation in \eqref{linearElasticity} can also be written in the form
\begin{equation}\label{displacementType}
 - \mu \Delta \bb{u} - (\lambda  + \mu ){\text{grad}}({\rm div} \bb{u}) = \bb{f} \quad {\text{in}}~~\Omega,
\end{equation}
which is referred to as the displacement type or the Navier type in what follows. In this case, we only consider $\Gamma_0 = \Gamma := \partial \Omega$. The first term $\Delta \bb{u}$ can be treated as the vector case of the Poisson equation, which is clearly observed in the implementation.

The continuous variational problem is to find $\bb{u}\in \bb{V}:=\bb{H}_0^1(\Omega)$ such that
\[a(\bb{u},\bb{v}) = (\bb{f}, \bb{v}), \quad \bb{v}\in \bb{V}, \]
where the bilinear form restriction to $K$ is
\begin{align*}
a^K(\bb{u},\bb{v})
& = \mu \int_K {\nabla \bb{u} \cdot \nabla {\bb{v}}} \mathrm{d}x + (\lambda  + \mu )\int_K ({\rm div} \bb{u})({\rm div} {\bb{v}}) \mathrm{d}x  \\
& = :\mu a_\nabla^K(\bb{u},{\bb{v}}) + (\lambda  + \mu )a_{{\rm div}}^K(\bb{u},{\bb{v}})
\end{align*}
with
\[a_\nabla^K(\bb{u},{\bb{v}}) = \int_K {\nabla \bb{u} \cdot \nabla {\bb{v}}} \mathrm{d}x,\quad a_{{\rm div}}^K(\bb{u},{\bb{v}}) = \int_K ({\rm div} \bb{u})({\rm div} {\bb{v}}) \mathrm{d}x.\]

The local virtual element space $\bb{V}_k(K)$ can be simply taken as the tensor-product of $V_k(K)$ defined in \eqref{first-VEMspace}, i.e., $\bb{V}_k(K) = (V_k(K))^2$. Let $\Pi^\nabla: \boldsymbol{V}_k(K) \to (\mathbb{P}_k(K))^2$ be the elliptic projector induced by $a_\nabla^K$. Then the approximate bilinear form can be split as
\[a_h^K(\bb{u},{\bb{v}}) = \mu a_{h,\nabla }^K(\bb{u},{\bb{v}}) + (\lambda+\mu) a_{h,{\rm div} }^K(\bb{u},{\bb{v}}),\]
where
\begin{align*}
& a_{h,\nabla }^K(\bb{u},{\bb{v}}) = a_\nabla^K({\Pi^\nabla}\bb{u},{\Pi^\nabla}{\bb{v}}) + S^K(\bb{u} - \Pi^\nabla\bb{u},\bb{v} - \Pi^\nabla\bb{v}),\\
& a_{h,{\rm div} }^K(\bb{u},{\bb{v}}) = (\Pi_{k - 1}^0{\rm div} \bb{u}, \Pi_{k - 1}^0{\rm div} \bb{v})_K,
\end{align*}
and the stabilization term for $k=1$ is given by
\[S^K(\bb{u},{\bb{v}}) = \sum\limits_{i = 1}^{N_v} \chi_i(\bb{u}) \cdot \chi_i(\bb{v}), \quad
\chi_i(\bb{u}) = [\chi_i(u_1), \chi_i(u_2)]^T.\]
The discrete problem is: Find $\bb{u}_h\in\bb{V}_h$ such that
\begin{equation}\label{VEMLinearElasticity}
a_h(\bb{u}_h,\,\bb{v}_h) = \langle \bb{f}_h,\,\bb{v}_h \rangle, \quad  \bb{v}_h\in\bb{V}_h,
\end{equation}
where
\[
a_h(\bb{u}_h,\,\bb{v}_h) = \sum\limits_{K\in \mathcal{T}_h} a_h^K(\bb{u}_h,{\bb{v}}_h)
\]
and the approximation of the right hand side is given by \cite{Zhang-Zhao-Yang-Chen-2019}
\begin{equation}\label{rhsLinearElasticity}
\langle \bb{f}_h, \bb{v}_h \rangle = (\bb{f}, P_h \bb{v}_h),
\end{equation}
where
\[P_h \bb{v}_h|_K = \frac{1}{|\partial K|} \int_{\partial K} \bb{v}_h {\rm d}s.\]

Let $\phi_1, \cdots ,\phi_N$ be the basis functions of the scalar space $V_k(K)$. Then the basis functions of the vector space $\bb{V}_k(K)$ can be defined by
\[{\overline \phi_1}, \cdots ,{\overline \phi_N},{\underline \phi_1}, \cdots ,{\underline \phi_N},\]
where
\[{\overline \phi_i} = \begin{bmatrix}
  \phi_i \\
  0
\end{bmatrix} ,\qquad
{\underline \phi_i} = \begin{bmatrix}
  0 \\
  \phi_i
\end{bmatrix},\quad i = 1, \cdots ,N.\]
Similarly, we introduce the notation
\begin{equation}\label{overunderm}
{\overline m_\alpha } = \begin{bmatrix}
  m_\alpha  \\
  0
\end{bmatrix},\qquad
{\underline m_\alpha } = \begin{bmatrix}
  0 \\
  m_\alpha
\end{bmatrix},\quad \alpha  = 1,\cdots,N_m.
\end{equation}
These vector functions will be written in a compact form as
\[\bb{m}^T = [{\overline m_1},\cdots,{\overline m_{N_m}},{\underline m_1},\cdots,{\underline m_{N_m}}] = :[\overline m^T,\underline m^T],\]
\[\bb \phi^T = [{\overline \phi_1}, \cdots ,{\overline \phi_{N_k}},{\underline \phi_1}, \cdots ,{\underline \phi_{N_k}}] = :[\overline \phi^T ,\underline \phi^T ].\]

We introduce the transition matrix $\bb{D}$ such that $\boldsymbol{m}^T = {\boldsymbol \phi^T}{\boldsymbol{D}}$. One easily finds that $\bb{D} = {\rm diag}(D, D)$, where $D$ is the transition matrix for the scalar case, i.e., $m^T = \phi^T D$. This block structure is also valid for the elliptic projection matrices. That is,
\[\bb{G} = {\rm diag}(G, G),  \quad \bb{B} = {\rm diag}(B, B), \quad \bb{G} = \bb{B}\bb{D},\]
where $G$ and $B$ are the same ones given in \eqref{GB}.

In the following, we consider the matrix expression of the $L^2$ projector $\Pi_{k-1}^0 {\rm div}$ satisfying
\[\int_K {\Pi_{k - 1}^0({\rm div} {\boldsymbol{v}})p} \mathrm{d}x = \int_K {({\rm div} {\boldsymbol{v}})p} \mathrm{d}x,~~p \in {\mathbb{P}_{k - 1}}(K).\]
The vector form can be written as
\[\int_K {{m^0}\Pi_{k - 1}^0({\rm div} {\boldsymbol \phi^T})} \mathrm{d}x = \int_K {{m^0}({\rm div} {\boldsymbol \phi^T})} \mathrm{d}x,\]
where $m^0$ is the vector consisting of the scaled monomials of order $\le k-1$. Let $\bb{\Pi}_{0*}$ be the matrix expression of $\Pi_{k - 1}^0({\rm div} {\boldsymbol \phi^T})$ in the basis $(m^0)^T$, i.e.,
\[\Pi_{k - 1}^0({\rm div} {\boldsymbol \phi^T}) = (m^0)^T\bb{\Pi}_{0*}.\]
Then $H_0 \bb{\Pi}_{0*} = C_0$, where
\[H_0 = \int_K {{m^0}{{({m^0})}^T}} \mathrm{d}x,\quad ~C_0 = \int_K {{m^0}({\rm div} {\boldsymbol \phi^T})} \mathrm{d}x.\]
For $k=1$, one has $m^0 = m_1 = 1$ and
\[
C_0=\int_K {{m^0}({\rm div} {\boldsymbol \phi^T})} \mathrm{d}x = \int_K {{\rm div} {\boldsymbol \phi^T}} \mathrm{d}x = \int_{\partial K} {{\boldsymbol \phi^T} \cdot {\boldsymbol n}} \mathrm{d}s,
\]
where
\begin{equation}\label{intphi}
\int_{\partial K} {{\boldsymbol \phi^T} \cdot {\boldsymbol n}} \mathrm{d}s = \int_{\partial K} {[\overline \phi^T   \cdot {\boldsymbol n},\underline \phi^T  \cdot {\boldsymbol n}]} \mathrm{d}s = \int_{\partial K} {[\phi^T  \cdot {n_x},\phi^T \cdot {n_y}]} \mathrm{d}s .
\end{equation}
We remark that \eqref{intphi} can be assembled along the boundary $\partial K$ using the technique in finite element methods.
\vspace{-0.8cm}
\begin{lstlisting}
v1 = 1:Nv;  v2 = [2:Nv,1];
H0 = area(iel);
C0 = zeros(1,2*Nv);
F = 1/2*[(1*Ne); (1*Ne)]; % [he*n1, he*n2]
C0(:) = accumarray([elem1(:);elem1(:)+Nv], F(:), [2*Nv 1]);
\end{lstlisting}

\begin{example}\label{ex:linearElasticity}
We consider a typical example to check the locking-free property of the proposed method.
The right-hand side $\bb{f}$ and the boundary conditions are chosen in such a way that the exact solution is
  \[\bb{u}(x,y) =  \begin{bmatrix}
  ( - 1 + \cos 2\pi x)\sin2\pi y \\
   - ( - 1 + \cos 2\pi y)\sin2\pi x
\end{bmatrix}  + \frac{1}{1 + \lambda}\sin \pi x\sin \pi y \begin{bmatrix}
  1 \\
  1
\end{bmatrix}.\]
\end{example}

\begin{figure}[!htb]
  \centering
  \includegraphics[scale=0.7]{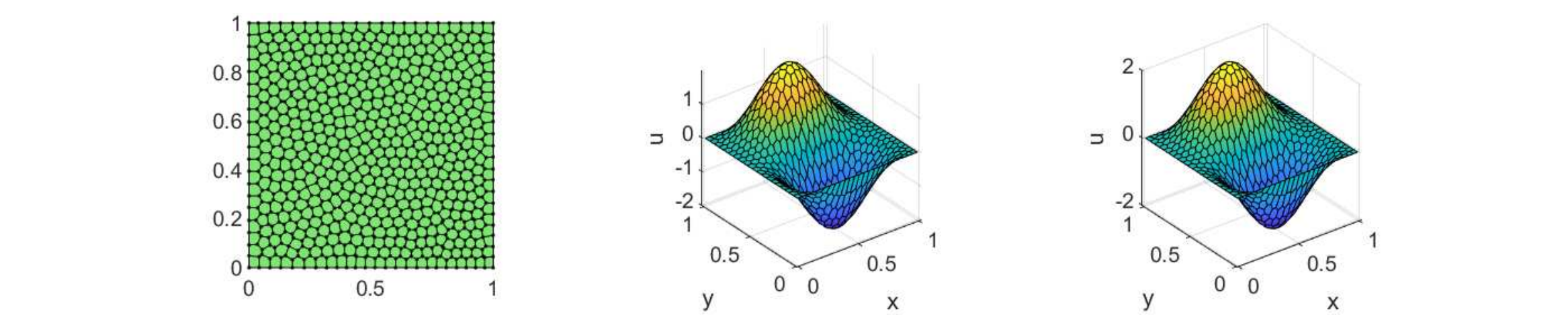}\\
  \caption{The exact and numerical nodal values for the linear elasticity problem of displacement type with $\lambda = 10^8$ and $\mu = 1$.}\label{fig:LockingfreeValue}
\end{figure}

\begin{figure}[!htb]
  \centering
  \includegraphics[scale=0.45]{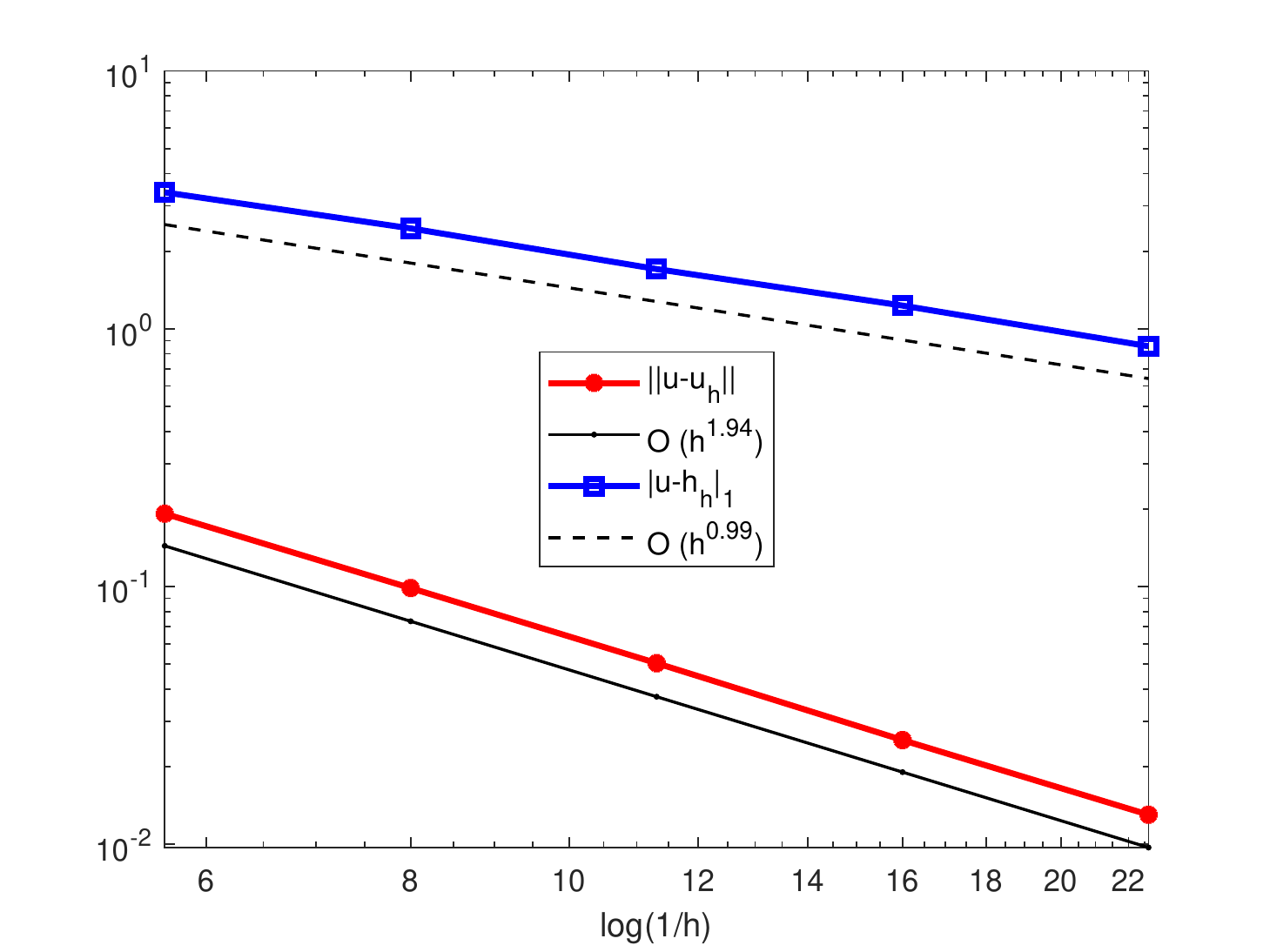}\\
  \caption{The performance of the VEM for the linear elasticity problem of displacement type with $\lambda = 10^8$ and $\mu = 1$.}\label{fig:Lockingfree}
\end{figure}

One easily finds that the proposed VEM is exactly the $\mathbb{P}_1$-Lagrange element when the polygonal mesh degenerates into a triangulation, in which case the method cannot get a uniform convergence with respect to the Lam\'{e} constant $\lambda$. However, for general polygonal meshes, the method seems to be robust with respect to $\lambda$ (see Fig.~\ref{fig:LockingfreeValue}), and
the optimal rates of convergence are achieved in the nearly incompressible case as shown in Fig.~\ref{fig:Lockingfree}, although we have not justified it in a theoretical way or found a counter example. Please refer to \mcode{main\_elasticityVEM\_Navier.m} for the test script.

\subsubsection{The tensor type}

Let $\bb{f}\in\bb{L}^2(\Omega)$. For simplicity, we still consider $\Gamma_0 = \partial \Omega$. The tensor-type variational formulation of \eqref{linearElasticity} is to find $\bb{u}\in\bb{V}$ such that
\begin{equation}
\label{Elasticity_TensorEq_Var}
a(\bb{u}, \bb{v}) = (\bb{f}, \bb{v}), \quad \bb{v} \in \bb{V},
\end{equation}
where
\begin{equation}\label{bilinear}
a(\bb{u}, \bb{v})
= 2\mu (\bb{\varepsilon}(\bb{u}), \bb{\varepsilon}(\bb{v}))
+ \lambda ({\rm div} \bb{u}, {\rm div} \bb{v}),
\end{equation}
and
\[
(\bb{f}, \bb{v})
= \int_{\Omega} \bb{f} \bb{v} \mathrm{d}x .
\]
For the ease of the presentation, we introduce the following notation
\[a_\mu(\bb{u}, \bb{v}) = (\bb{\varepsilon}(\bb{u}), \bb{\varepsilon}(\bb{v})), \qquad
a_\lambda(\bb{u}, \bb{v}) = ({\rm div} \bb{u}, {\rm div} \bb{v}).\]
The conforming VEM for \eqref{Elasticity_TensorEq_Var} is proposed in \cite{Beirao-Brezzi-Marini-2013},
where a locking-free analysis is carried out for the virtual element spaces of order $k\ge 2$. The order requirement is to
ensure the so-called discrete inf-sup condition and the optimal convergence.

Introduce the elliptic projection
\[{\Pi^a}:~{{\boldsymbol{V}}_k}(K) \to {(\mathbb{P}_k(K))^2},\quad {\boldsymbol{v}} \mapsto {\Pi^a}{\boldsymbol{v}},\]
satisfying
\begin{equation}\label{ell2}
a_\mu^K({\Pi^a}{\boldsymbol{v}},\boldsymbol{p}) = a_\mu^K({\boldsymbol{v}},\boldsymbol{p}),\quad \boldsymbol{p} \in {(\mathbb{P}_k(K))^2}.
\end{equation}
Note that $\Pi^a\boldsymbol{v}$ is unique up to an additive vector in
\[K_0 :=RM = {\rm span} \left\{
\begin{bmatrix} 1 \\ 0 \end{bmatrix}, ~~
\begin{bmatrix} 0 \\ 1 \end{bmatrix},~~
\begin{bmatrix} -y \\ x \end{bmatrix}  \right\} .\]
For this reason, we can impose one of the following constraints:
\begin{itemize}
  \item[Choice 1:]
  \begin{equation}\label{ellv2}
\sum\limits_{i = 1}^{N_v} (\Pi^a\boldsymbol{v}(z_i),\boldsymbol{p}(z_i))  = \sum\limits_{i = 1}^{N_v} (\boldsymbol{v}(z_i),\boldsymbol{p}(z_i)),~~\boldsymbol{p} \in K_0,
\end{equation}
  \item[Choice 2:]
  \begin{equation}\label{linearcond2}
\int_{\partial K} (\Pi^a\boldsymbol{v})\cdot\boldsymbol{p} \mathrm{d}s = \int_{\partial K} \boldsymbol v\cdot \boldsymbol p \mathrm{d}s,~~\boldsymbol{p} \in K_0.
\end{equation}
  \item[Choice 3:]
    \begin{align}\label{constrainPiC}
  \int_K {\rm rot}~ \Pi_K^1\boldsymbol v \mathrm{d}x & = \int_K {\rm rot} \boldsymbol v \mathrm{d}x ,\nonumber\\
  \int_{\partial K} \Pi_K^1\boldsymbol v \mathrm{d}s  & = \int_{\partial K} \boldsymbol v \mathrm{d}s,
\end{align}
  where ${\rm rot} \bb{v} = \nabla  \times \bb{v} = \partial_1v_2 - \partial_2v_1$. An integration by parts gives
  \[\int_K {\rm rot} \bb v {\rm d}x =  \int_{\partial K} \bb v \cdot \bb{t}_K {\rm d}s, \]
where $\bb{t}_K = (-n_2, n_1)^T$ is the anti-clockwise tangential along $\partial K$.
\end{itemize}
It is evident that all the three choices can be computed by using the given d.o.f.s.

The discrete problem is the same as the one in \eqref{VEMLinearElasticity} with the local approximate bilinear form replaced by
\[a_h^K(\bb{v},\bb{w}) = 2\mu a_{\mu,h}^K(\bb{v},\bb{w}) + \lambda a_{\lambda,h}^K(\bb{v},\bb{w}), \]
where
\begin{align*}
& a_{\mu,h}^K(\bb{v},\bb{w}) = a_{\mu}^K(\Pi^a\bb{v},\Pi^a\bb{w}) + S^K(\bb{v}-\Pi^a\bb{v}, \bb{w}-\Pi^a\bb{w}), \\
& a_{\lambda,h}^K(\bb{v},\bb{w}) = (\Pi_{k-1}^0{\rm div}\bb{v}, \Pi_{k-1}^0{\rm div}\bb{w})_K.
\end{align*}

For the lowest-order case $k=1$, an integration by parts gives
\begin{align}
\bb{B}
& = a_\mu^K(\bb{m},\bb \phi^T)
       = \int_{\partial K} (\bb{\varepsilon}(\bb{m}) \cdot {\bb n}) \cdot {\bb \phi^T} \mathrm{d}s  \label{BLinearElas} \\
& = \Big[ \varepsilon_{11}(\boldsymbol{m})\int_{\partial K} \phi^Tn_x \mathrm{d}s  + \varepsilon_{12}(\boldsymbol{m})\int_{\partial K} \phi^Tn_y \mathrm{d}s, \nonumber \\
& \hspace{2cm} \varepsilon_{21}(\boldsymbol{m})\int_{\partial K} \phi^Tn_x \mathrm{d}s + \varepsilon_{22}(\boldsymbol{m})\int_{\partial K} \phi^Tn_y  \mathrm{d}s \Big] \nonumber
\end{align}
with
\begin{align*}
& \varepsilon_{11}(\boldsymbol{m}) = [0,\frac{1}{h_K},0,0,0,0]^T,  \qquad
  \varepsilon_{22}(\boldsymbol{m}) = [0,0,0,0,0,\frac{1}{h_K}]^T, \\
& \varepsilon_{12}(\boldsymbol{m}) = \varepsilon_{21}(\boldsymbol{m}) = [0,0,\frac{1}{2h_K},0,\frac{1}{2h_K},0]^T.
\end{align*}
Using the assembling technique for FEMs, the computation of $\bb{B}$ in MATLAB reads
\vspace{-0.8cm}
\begin{lstlisting}
elem1 = [v1(:), v2(:)];
C0 = zeros(1,2*Nv);
F = 1/2*[(1*Ne); (1*Ne)]; % [he*n1, he*n2]
C0(:) = accumarray([elem1(:);elem1(:)+Nv], F(:), [2*Nv 1]);
E = zeros(6,4);
E(2,1) = 1/hK;  E([3,5],[2,3]) = 1/(2*hK); E(6,4) = 1/hK;
B = [E(:,1)*C01x+E(:,2)*C01y, E(:,3)*C01x+E(:,4)*C01y];
\end{lstlisting}

Note that for the tensor type, we can consider the pure traction problem, i.e., $\Gamma_1 = \partial \Omega$. The numerical results are similar to that of the VEM of displacement type. The test script is \mcode{main\_elasticityVEM.m}. It is worth pointing out that we proposed in \cite{Huang-Lin-Yu-2021} a novel conforming locking-free method with a rigorous proof, where several benchmarks are tested and the test script can be found in \mcode{main\_elasticityVEM\_reducedIntegration.m}.

\subsection{Nonconforming VEMs}

Nonconforming VEMs for the linear elasticity problems are first introduced in \cite{Zhang-Zhao-Yang-Chen-2019}
for the pure displacement/traction formulation in two or three dimensions.
The proposed method is robust with respect to the Lam\'{e} constant for $k\ge2$.

For $k=1$, one easily finds that
\[\int_{\partial K} \phi^Tn_x {\rm d}s = [h_{e_1}n_{1,x}, \cdots ,  h_{e_{N_v}} n_{N_v,x} ], \]
where $h_{e_i}$ is the length of the edge $e_i$ and $\bb{n}_{e_i} = (n_{i,x}, n_{i,y})^T$. In this case, the computation of $\bb{B}$ in \eqref{BLinearElas} for the tensor type reads
\vspace{-0.8cm}
\begin{lstlisting}
C01x = Ne(:,1)';  C01y = Ne(:,2)';
B = [E(:,1)*C01x+E(:,2)*C01y, E(:,3)*C01x+E(:,4)*C01y];
\end{lstlisting}

\begin{figure}[!htb]
  \centering
  \subfigure[Displacement type]{\includegraphics[scale=0.35]{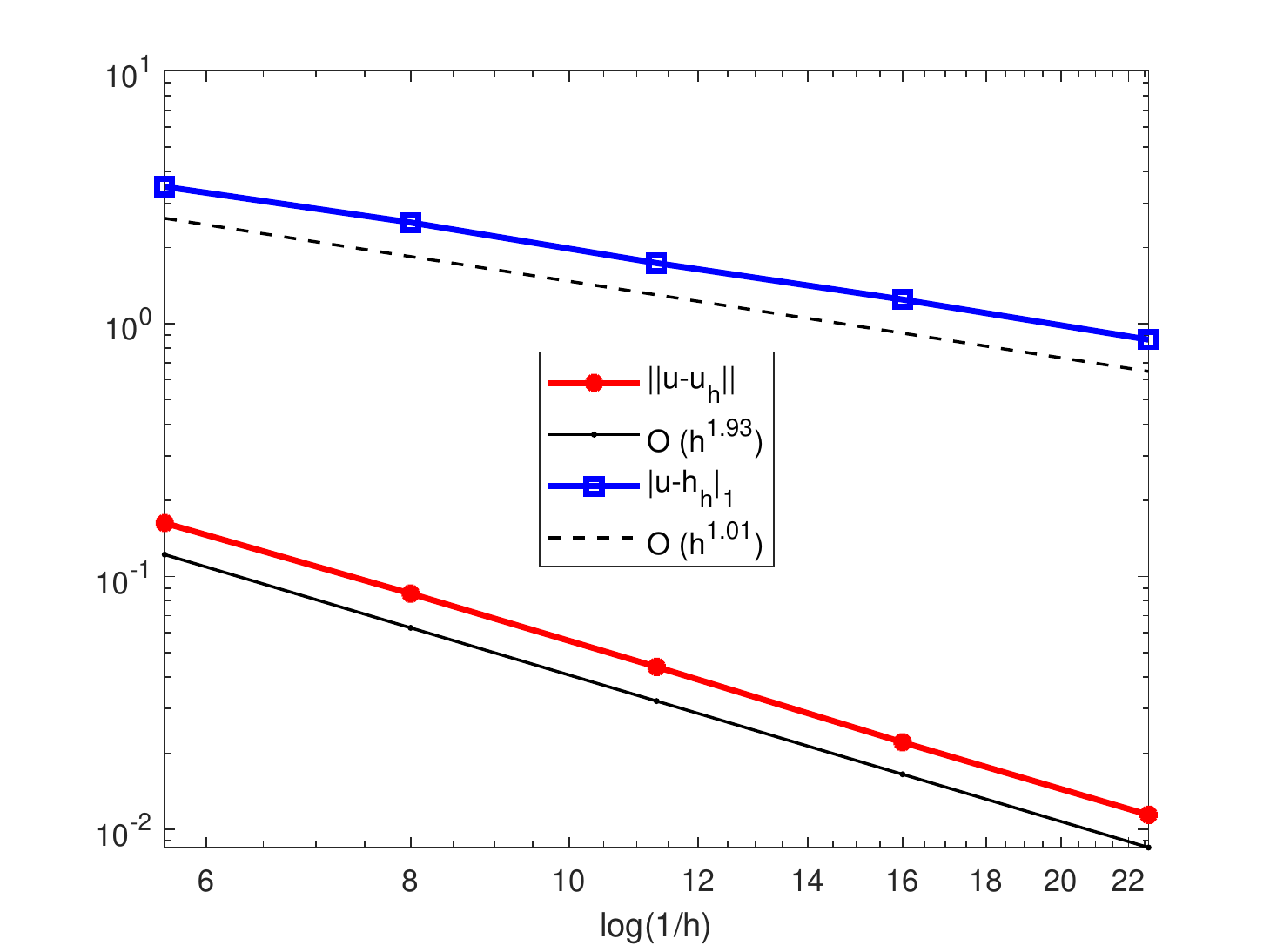}}
  \subfigure[Tensor type]{\includegraphics[scale=0.35]{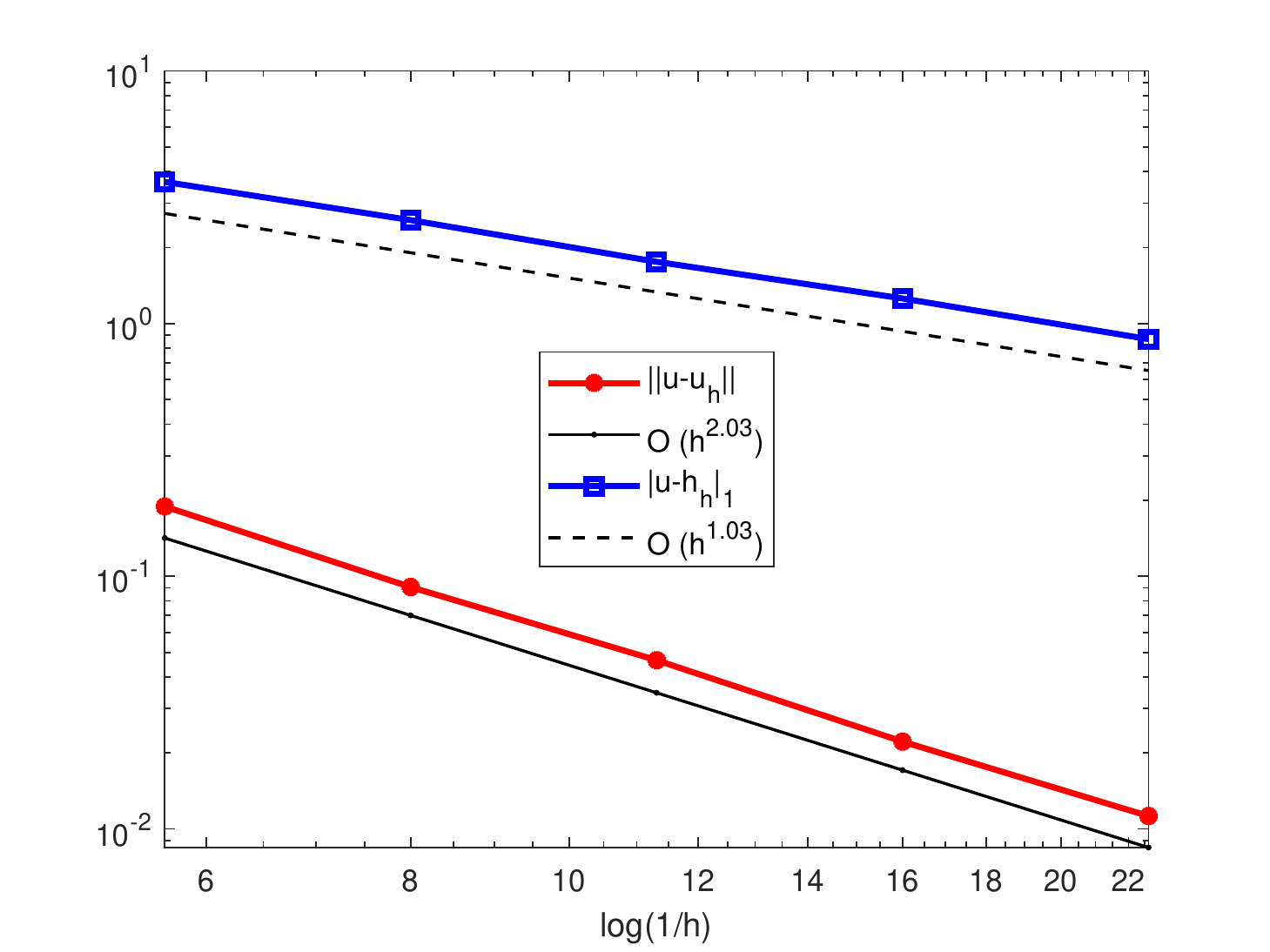}}
  \subfigure[Tensor type with the continuous treatment]{\includegraphics[scale=0.35]{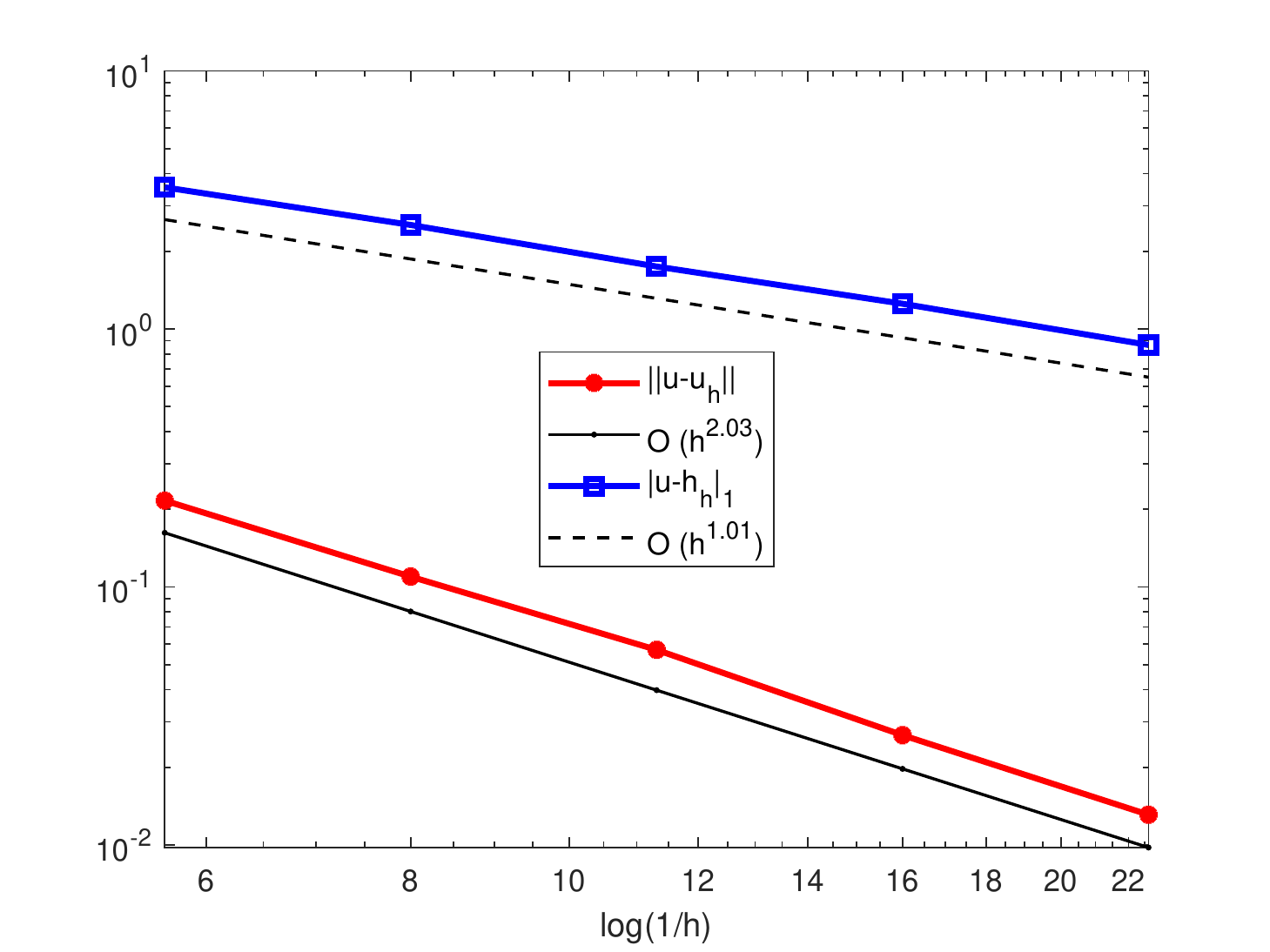}}\\
  \caption{Convergence rates of the nonconforming VEMs for the linear elasticity.}\label{fig:ElasNc}
\end{figure}

We provide the implementation of the displacement type and tensor type in the original nonconforming spaces with the results displayed in Fig.~\ref{fig:ElasNc}(a)(b). Similar to the Poisson equation, we also present the realization of the tensor type in the modified nonconforming spaces and the convergence rate is shown in Fig.~\ref{fig:ElasNc}(c). Here, we still consider the test in Example \ref{ex:linearElasticity} with $\lambda = 10^{10}$ and $\mu=1$. As again observed in Fig.~\ref{fig:ElasNc}, the optimal rates of convergence are obtained for all the three methods with general polygonal meshes applied in the nearly incompressible case although we cannot justify it or provide a counter example.

The test scripts are \mcode{main\_elasticityVEM\_NavierNC.m}, \mcode{main\_elasticityVEM\_NC.m} and \\ \mcode{main\_elasticityVEM\_NCb.m}, respectively.

\subsection{Several locking-free VEMs}

The authors in \cite{Kwak-Park-2022} present two kinds of lowest-order VEMs with consistent convergence, in which the first one is achieved by introducing a special stabilization term to ensure the discrete Korn's inequality, and the second one can be seen as an extension of the idea of Kouhia and Stenberg suggested in \cite{Kouhia-Stenberg-1995} to the virtual element method.
We provide the implementation of the second method in \mcode{elasticityVEM\_KouhiaStenberg.m}. In this approach, the local space is taken as
\[\bb{V}(K) = V_1^{nc}(K) \times V_1(K),\]
where $V_1^{nc}(K)$ and $V_1(K)$ are the lowest-order nonconforming and conforming virtual element spaces, respectively. In this case, the computation for $k=1$ reads
\vspace{-0.8cm}
\begin{lstlisting}
C01xe = Ne(:,1)';  C01ye = Ne(:,2)';
C01xv = 0.5*(Ne(p1,1)+Ne(p2,1))';
C01yv = 0.5*(Ne(p1,2)+Ne(p2,2))';
B = [E(:,1)*C01xe+E(:,2)*C01ye, E(:,3)*C01xv+E(:,4)*C01yv];
\end{lstlisting}

We also develop a lowest-order nonconforming virtual element method for planar linear elasticity in \cite{NonconformingLockingfree}, which can be viewed as an extension of the idea in \cite{Falk-1991} to the virtual element method, with the family of polygonal meshes satisfying a very general geometric assumption. The method is shown to be uniformly convergent for the nearly incompressible case with optimal rates of convergence. In addition, we provide a unified locking-free scheme both for the conforming and nonconforming VEMs in the lowest order case. The implementation and the numerical test can be found in \cite{NonconformingLockingfree}. The test scripts are \mcode{main\_elasticityVEM\_NCreducedIntegration.m} and \mcode{main\_elasticityVEM\_NCUniformReducedIntegration.m}, respectively.

\section{Plate bending problems}

In this section we focus on the plate bending problem in the form of
\[\begin{cases}
- \partial_{ij}M_{ij}(w) = f \quad & \mbox{in}~~\Omega \subset \mathbb{R}^2, \\
w = \partial_{\bb n}w = 0 \quad & \mbox{on} ~~ \partial \Omega,
\end{cases}\]
where
\[M_{ij} = D\left( (1 - \nu )K_{ij} + \nu K_{kk}\delta_{ij} \right),~~~i,j = 1,2,~~k \in \{ 1,2\} ,\]
\[K_{ij} =  - \partial_{ij}w,~~~i,j = 1,2, \quad D = \frac{Et^3}{12(1 - \nu^2)}.\]
Note that we have used the summation convention whereby summation is implied when an index is repeated exactly
two times. Please refer to \cite{FengShi1996,Brezzi-Marini-2013,Chinosi-Marini-2016} for details. The variational problem is to find $w\in V: = H_0^2(\Omega )$ such that
\[a(w,v) = \ell(v),\quad v \in V,\]
where
\[a(w,v) = \int_\Omega  M_{ij}(w)K_{ij}(v) \mathrm{d}x , \qquad \ell(v) = \int_\Omega  fv \mathrm{d}x.\]

\subsection{$C^1$-continuous VEMs}

We first recall the $H^2$-conforming virtual element space $V_k^{2,c}(K)$ introduced in \cite{Brezzi-Marini-2013,Chinosi-Marini-2016}. For $k\ge 3$, define
\[{V_k^{2,c}}(K)= \{ v \in H^2(K): \Delta ^2v \in \mathbb{P}_{k - 4}(K),~ v|_e \in \mathbb{P}_k(e),~ \partial_{\bb n} v|_e \in {\mathbb{P}_{k - 1}}(e),~e \subset \partial K \},\]
while for the lowest order $k=2$, the space is modified as
\[V_2^{2,c}(K)= \Big\{ v \in H^2(K): \Delta ^2v = 0,~ v |_e \in \mathbb{P}_3(e),~ \partial_{\bb n}v |_e \in \mathbb{P}_1(e),~e \subset \partial K \Big\}.\]
The d.o.f.s are:
\begin{itemize}
  \item The values of $v(z)$ at the vertices of $K$.
  \item The values of $h_z\partial _1v(z)$ and $h_z\partial _2v(z)$ at the vertices of $K$, where $h_z$ is a characteristic length attached to each vertex $z$, for instance, the average of the diameters of the elements having $z$ as a vertex.
  \item The moments of $v$ on edges up to degree $k-4$,
  \[\chi _e(v) = | e |^{-1}(m_e,v)_e,\quad m_e \in \mathbb{M}_{k - 4}(e).\]
  \item The moments of $\partial_{\bb n}v$ on edges up to degree $k-3$,
  \[\chi _{n_e}(v) = (m_e, \partial_{\bb n}v)_e,\quad m_e \in \mathbb{M}_{k - 3}(e).\]
  \item The moments on element $K$ up to degree $k-4$,
  \[\chi _K(v) = | K |^{-1}(m_K, v)_K,\quad m_K \in \mathbb{M}_{k - 4}(K).\]
\end{itemize}

The elliptic projection can be defined in vector form as
\begin{equation}\label{H2C1Elliptic}
\begin{cases}
  a^K(m,\Pi^K\phi^T) = a^K(m,\phi^T),  \\
  P_0^1(\Pi^K\phi^T) = P_0^1(\phi^T),  \\
  P_0^2(\nabla \Pi^K\phi^T) = P_0^2(\nabla \phi^T),
\end{cases}
\end{equation}
where
\[P_0^1(v) = \frac{1}{N_v}\sum\limits_{i = 1}^{N_v} v(a_i), \quad P_0^2(v) = \int_{\partial K} v \mathrm{d}s.\]
Introduce the following notation
\[Q_{3i} = M_{ij,j},~~Q_{3\bb n} = Q_{3i}n_i, \quad
M_{\bb{nn}} = M_{ij}n_in_j,~~M_{\bb{tn}} = M_{ij}t_in_j,\]
with $\bb{n} = (n_1,n_2)^T$ and $\bb{t} = (-n_2,n_1)^T$. We then have the integration by parts formula
\begin{align*}
  a^K(p,v)
  & =  - \int_K  Q_{3i,i}(p)v \mathrm{d}x   \\
  & \quad + \int_{\partial K } \left( Q_{3\bb n}(p) + \partial_{\bb t}M_{\bb{tn}}(p) \right)v \mathrm{d}s- \int_{\partial K } M_{\bb{nn}}(p)\partial_{\bb n}v \mathrm{d}s  \\
  & \quad + \sum\limits_{i = 1}^{N_v} [M_{\bb{tn}}(p)](z_i)v(z_i),
\end{align*}
which implies the computability of the elliptic projection $\Pi^K v$ for any $v\in V_k^{2,c}(K)$, where the last term is a jump at $z_i$ along the boundary $\partial K$, i.e.,
\[[M_{\bb{tn}}(p)](z_i) = M_{\bb{tn}}(p) \Big|_{z_i^ - }^{z_i^ + }.\]

We now consider the computation of the elliptic projection matrices. In the lowest order case, the d.o.f.s are
\begin{itemize}
      \item The values of $v$ at the vertices of $K$,
      \[\chi_a(v) = v(z_i),\quad i = 1, \cdots ,N_v.\]
      \item The values of $\nabla v$ at the vertices of $K$,
      \[\chi_{a1}(v) = \partial_xv(z_i),~~\chi_{a2}(v) = \partial_yv(z_i),\quad i = 1, \cdots ,N_v.\]
\end{itemize}
In the implementation, they are arranged as
\begin{align*}
& \chi_i(v) = v(z_i), \quad i = 1, \cdots ,N_v, \\
& \chi_{N_v + i}(v) = h_\xi\partial_xv(z_i),  \quad i = 1, \cdots ,N_v, \\
& \chi_{2N_v + i}(v) = h_\xi\partial_yv(z_i), \quad i = 1, \cdots ,N_v,
\end{align*}
where $h_\xi$ is the characteristic length attached to each vertex.
To compute the transition matrix $D$, we first provide the characteristic lengths by using the data structure \mcode{node2elem}, given as
\vspace{-0.8cm}
\begin{lstlisting}
% characteristic length
hxi = cellfun(@(id) mean(diameter(id)), node2elem);
index = elem{iel};
hxiK = hxi(index);
\end{lstlisting}
Then the computation of the matrix $D$ reads
\vspace{-0.8cm}
\begin{lstlisting}
% ------- scaled monomials ---------
% m'
m = @(x,y) [1+0*x, (x-xK)/hK, (y-yK)/hK, (x-xK).^2/hK^2, ...
            (x-xK).*(y-yK)/hK^2, (y-yK).^2/hK^2]; % m1,...,m6
% Dx(m'), Dy(m')
Dxm = @(x,y) [0*x, 1/hK+0*x, 0*x, 2*(x-xK)/hK^2, (y-yK)/hK^2, 0*x];
Dym = @(x,y) [0*x, 0*x, 1/hK+0*x, 0*x, (x-xK)/hK^2, 2*(y-yK)/hK^2];

% ------ transition matrix ---------
D = zeros(Ndof,Nm);
D(1:Nv,:) = m(x,y);
D(Nv+1:2*Nv,:) = repmat(hxiK,1,Nm).*Dxm(x,y);
D(2*Nv+1:end,:) = repmat(hxiK,1,Nm).*Dym(x,y);
\end{lstlisting}

For $k=2$, one has
\begin{align}
  B_{\alpha j}
  & = {a^K}(m_\alpha ,\phi_j) =  - \sum\limits_{e \subset \partial K} \int_e M_{\bb{nn}}(m_\alpha)\partial_{\bb n}\phi_j \mathrm{d}s  + \sum\limits_{i = 1}^{N_v} [M_{\bb{tn}}(m_\alpha )](z_i)\phi_j(z_i)  \nonumber \\
 & = :J_1(\alpha ,j) + J_2(\alpha ,j). \label{BPlate}
\end{align}
Since $\partial_{\bb n}\phi_j|_e \in \mathbb{P}_1(e)$, the trapezoidal rule gives
\[J_1(\alpha ,j) =   - \sum\limits_{i = 1}^{{N_v}} {{{\left. {{M_{\bb{nn}}}({m_\alpha })} \right|}_{{e_i}}}\frac{{\left| {{e_i}} \right|}}{2}\left( {{\partial_{\bb n}}{\phi_j}({z_i}) + {\partial_{\bb n}}{\phi_j}({z_{i + 1}})} \right)} . \]
Noting that
\[{\chi_{a1,i}}(v) = {h_\xi }{\partial_x}v({z_i})\quad ~~ \Rightarrow \quad ~~{\partial_x}v({z_i}) = \frac{1}{{{h_\xi }}}{\chi_{a1,i}}(v),\quad i = 1,\cdots, Nv,\]
we then compute the values of $\partial_{\bb n}\phi_j$ at the vertices as follows.
\vspace{-0.8cm}
\begin{lstlisting}
% Dx(phi'), Dy(phi') at z1,...,zNv (each row)
Dxphi = zeros(Nv,Ndof); Dyphi = zeros(Nv,Ndof);
Dxphi(:,Nv+1:2*Nv) = eye(Nv)./repmat(hxiK,1,Nv);
Dyphi(:,2*Nv+1:end) = eye(Nv)./repmat(hxiK,1,Nv);
\end{lstlisting}

For $J_2$ with the entry given by
\[J_2(\alpha ,j) = \sum\limits_{i = 1}^{N_v} [M_{\bb{tn}}(m_\alpha )](z_i)\phi_j(z_i) ,\]
by the definition of the jump,
\[[{M_{\bb{tn}}}({m_\alpha })]({z_i}) = \left. {{M_{\bb{tn}}}({m_\alpha })} \right|_{z_i^ - }^{z_i^ + } = {M_{\bb{tn}}}({m_\alpha })(z_i^ + ) - {M_{\bb{tn}}}({m_\alpha })(z_i^ - ),\]
where $M_{\bb{tn}}(m_\alpha )(z_i^ + )$ is the evaluation on the edge $e_i$ to the right of $z_i$, $M_{\bb{tn}}(m_\alpha )(z_i^ - )$ is the evaluation on the edge $e_{i - 1}$, and
\[{\phi^T}({z_i}) = [{{\bb{e}}_i},{\bb{0}},{\bb{0}}],\]
where $\bb{e}_i$ is a zero vector with $i$-th entry being 1.

The above discussion is summarized in the following code.
\vspace{-0.8cm}
\begin{lstlisting}
% --------- elliptic projection -----------
% \partial_ij (m)
D11 = zeros(Nm,1); D11(4) = 2/hK^2;
D12 = zeros(Nm,1); D12(5) = 1/hK^2;
D22 = zeros(Nm,1); D22(6) = 2/hK^2;
% Mij(m)
M11 = -para.D*((1-para.nu)*D11 + para.nu*(D11+D22));
M12 = -para.D*(1-para.nu)*D12;
M22 = -para.D*((1-para.nu)*D22 + para.nu*(D11+D22));
% Mnn(m) on e1,...,eNv
n1 = ne(:,1);  n2 = ne(:,2);
Mnn = M11*(n1.*n1)' + M12*(n1.*n2+n2.*n1)' + M22*(n2.*n2)';
% Mtn(m) on e1,...,eNv
t1 = te(:,1); t2 = te(:,2);
Mtn = M11*(t1.*n1)' + M12*(t1.*n2+t2.*n1)' + M22*(t2.*n2)';
% Dx(phi'), Dy(phi') at z1,...,zNv (each row)
Dxphi = zeros(Nv,Ndof); Dyphi = zeros(Nv,Ndof);
Dxphi(:,Nv+1:2*Nv) = eye(Nv)./repmat(hxiK,1,Nv);
Dyphi(:,2*Nv+1:end) = eye(Nv)./repmat(hxiK,1,Nv);
% B, Bs, G, Gs
B = zeros(Nm,Ndof);
p1 = [Nv,1:Nv-1]; p2 = 1:Nv;
for j = 1:Nv % loop of edges
    % int[\partial_n (phi')] on ej
    Dnphi1 = Dxphi(v1(j),:)*n1(j) + Dyphi(v1(j),:)*n2(j); % zj
    Dnphi2 = Dxphi(v2(j),:)*n1(j) + Dyphi(v2(j),:)*n2(j); % z_{j+1}
    nphi = 0.5*he(j)*(Dnphi1+Dnphi2);
    % Jump(m) at zj
    Jump = Mtn(:,p2(j))-Mtn(:,p1(j));
    % phi' at zj
    phi = zeros(1,Ndof);  phi(j) = 1;
    % B1 on e and at zj
    B = B - Mnn(:,j)*nphi + Jump*phi;
end
\end{lstlisting}

The first constraint in \eqref{H2C1Elliptic} is
\[\tilde B(1,:) = P_0^1(\phi^T) = \frac{1}{N_v}[\bb{1},\bb{0},\bb{0}].\]
The second constraint can be computed by using the integration by parts as
\begin{align*}
  \int_{\partial K} \nabla v \mathrm{d}s
  & = \sum\limits_{e \subset \partial K} \int_e \partial_{\bb{n}_e}v\bb{n}_e \mathrm{d}s + \int_e \partial_{\bb{t}_e}v \bb{t}_e \mathrm{d}s   \\
  & = \sum\limits_{e \subset \partial K} \bb{n}_e\int_e \partial_{\bb n}v \mathrm{d}s  + \sum\limits_{i = 1}^{N_v} \bb{t}_{e_i}(v(z_{i + 1}) - v(z_i)) ,
\end{align*}
where
\[ \int_{\partial K} \nabla \phi_j \mathrm{d}s = \bb{t}_{e_{j - 1}} - \bb{t}_{e_j}.\]
The computation reads
\vspace{-0.8cm}
\begin{lstlisting}
Bs = B;
% first constraint
Bs(1,1:Nv) = 1/Nv;
% second constraint
for j = 1:Nv % loop of edges
    Bs(2:3,1:Nv) = te([Nv,1:Nv-1],:)' - te';
    Dnphi1 = Dxphi(v1(j),Nv+1:end)*n1(j) + Dyphi(v1(j),Nv+1:end)*n2(j); % zj
    Dnphi2 = Dxphi(v2(j),Nv+1:end)*n1(j) + Dyphi(v2(j),Nv+1:end)*n2(j); % z_{j+1}
    Nphi = 0.5*Ne(j,:)'*(Dnphi1+Dnphi2); % scaled
    Bs(2:3,Nv+1:end) = Bs(2:3,Nv+1:end) + Nphi;
end
\end{lstlisting}

\begin{example}\label{ex:Plate}
The exact solution is chosen as $u = \sin(2\pi x) \cos(2 \pi y)$ with the parameters $t = 0.1$, $E = 10920$ and $\nu = 0.3$.
\end{example}

\begin{figure}[!htb]
  \centering
  \includegraphics[scale=0.45]{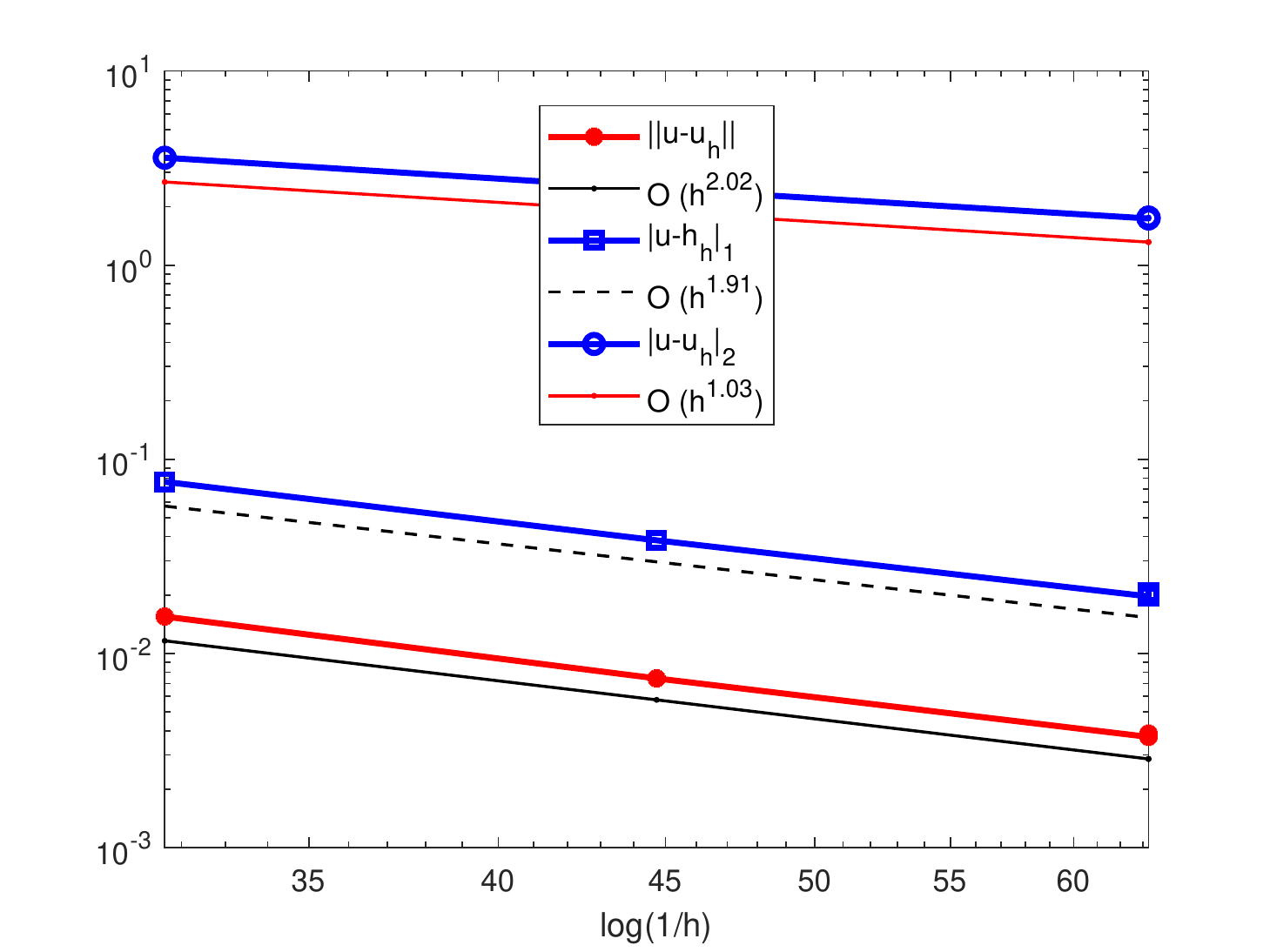}\\
  \caption{Convergence rates of the plate bending problem for the $C^1$-continuous virtual elements}\label{fig:PlateC1}
\end{figure}

We consider a sequence of meshes, which is a Centroidal Voronoi Tessellation of the unit square in 32, 64, 128, 256 and 512 polygons. The results are shown in Fig.~\ref{fig:PlateC1}. The optimal rate of convergence of the discrete $H^2$-norm (1st order), $H^1$-norm (2nd order) and $L^2$-norm (2nd order) is observed for the lowest order $k = 2$ when meshes are fine enough. The test script is \mcode{main\_PlateBending\_C1VEM.m}.

\subsection{$C^0$-continuous VEMs}

The Ref.~\cite{Zhao-Chen-Zhang-2016} gives the $C^0$-continuous nonconforming virtual element method for the plate bending problem, with the local virtual element space defined as
\[V_k^{2,0}(K) = \{ v \in H^2(K): \Delta^2v \in \mathbb{P}_{k - 4}(K),~~v|_{\partial K} \in \mathbb{B}_k(\partial K) \},\]
where $K$ is a convex polygon and
\[\mathbb{B}_k(\partial K) = \{ v \in C^0(\partial K): v|_e \in \mathbb{P}_k(e),~~\Delta v|_e \in \mathbb{P}_{k-2}(e),~~e \subset \partial K \}.\]
A function in $V_k^{2,0}(K)$ is uniquely identified by the following degrees of freedom:
\begin{itemize}
  \item The values at the vertices of $K$,
  \[\chi_i(v) = v(z_i),\quad z_i~\mbox{is a vertex of $K$}.\]
  \item The moments of $v$ on edges up to degree $k-2$,
  \[\chi_e(v) = |e|^{-1} (m_e,v)_e,\quad v \in \mathbb{M}_{k - 2}(e).\]
  \item The moments of $\partial_{\bb n} v$ on edges up to degree $k-2$,
  \[\chi_{n_e}(v) = (m_e, \partial_{\bb n}v)_e,\quad m \in \mathbb{M}_{k - 2}(e).\]
  \item The moments of $v$ on $K$ up to degree $k-4$,
  \[\chi_K(v) = |K|^{-1} (m_K,v)_K,\quad v \in \mathbb{M}_{k - 4}(K).\]
\end{itemize}
\begin{figure}[!htb]
  \centering
  \includegraphics[scale=0.25]{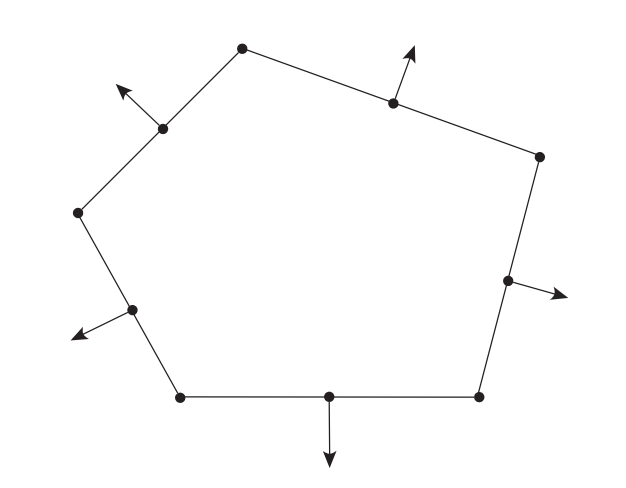}\\
 \caption{Local degrees of freedom for the $C^0$-continuous virtual elements~($k=2$)} \label{fig:C0Morley}
\end{figure}
We only consider the lowest order case $k=2$. The d.o.f.s contain the values at the vertices and the following moments
\[\chi_e(v) = \frac{1}{| e |}\int_e v \mathrm{d}s,\quad \chi_{n_e}(v) = \int_e \partial_{\bb n}v \mathrm{d}s,\]
as shown in Fig.~\ref{fig:C0Morley}. Note that the first-type moments can also be replaced by the midpoint values on edges and the second-type moment has different signs when restricted to the left and right elements of an interior edge. Let $N_v$ be the number of vertices. Then there are $N_k=3N_v$ d.o.f.s on each element, arranged as
\begin{itemize}
  \item The values at the vertices of $K$,
  \[\chi_i(v) = v(z_i),\quad i = 1, \cdots ,N_v.\]
  \item The mid-point values on each edge of $K$,
  \[\chi_{i+N_v}(v) = v(m_i),\quad i = 1, \cdots ,N_v.\]
  \item The moments
  \[\chi_{2N_v+i}(v) = \int_{e_i} \partial_{\bb n}v \mathrm{d}s,\quad i = 1, \cdots ,N_v.\]
\end{itemize}

Since
\[\int_e \partial_{\bb n}\phi^T \mathrm{d}s = [\bb{0},\bb{0},\bb{e}_j], \quad  \phi^T(z_i) = [\bb{e}_i,\bb{0},\bb{0}],\]
the matrix $B$ for the elliptic projection can be realized as
\vspace{-0.8cm}
\begin{lstlisting}
B = zeros(Nm,Ndof);
p1 = [Nv,1:Nv-1]; p2 = 1:Nv;
for j = 1:Nv % loop of edges
    % nphi on ej
    nphi = zeros(1,Ndof); nphi(2*Nv+j) = 1;
    % Jump at zj
    Jump = Mtn(:,p2(j))-Mtn(:,p1(j));
    % phi at zj
    phi = zeros(1,Ndof); phi(j) = 1;
    % B1
    B = B - Mnn(:,j)*nphi + Jump*phi;
end
\end{lstlisting}
Noting that
\begin{align*}
 \int_{\partial K} \nabla \phi_j \mathrm{d}s = \bb{t}_{e_{j - 1}} - \bb{t}_{e_j},  \quad
 \int_{\partial K} \nabla \phi_{N_v + j} \mathrm{d}s = \boldsymbol{0}, \quad
 \int_{\partial K} \nabla \phi_{2N_v + j} \mathrm{d}s = \bb{n}_{e_j},
\end{align*}
we can compute the constraints as follows.
\vspace{-0.8cm}
\begin{lstlisting}
Bs = B;
% first constraint
Bs(1,1:Nv) = 1/Nv;
% second constraint
Bs(2:3,1:Nv) = te([Nv,1:Nv-1],:)' - te';
Bs(2:3,2*Nv+1:end) = ne';
\end{lstlisting}

\begin{figure}[!htb]
  \centering
  \includegraphics[scale=0.45]{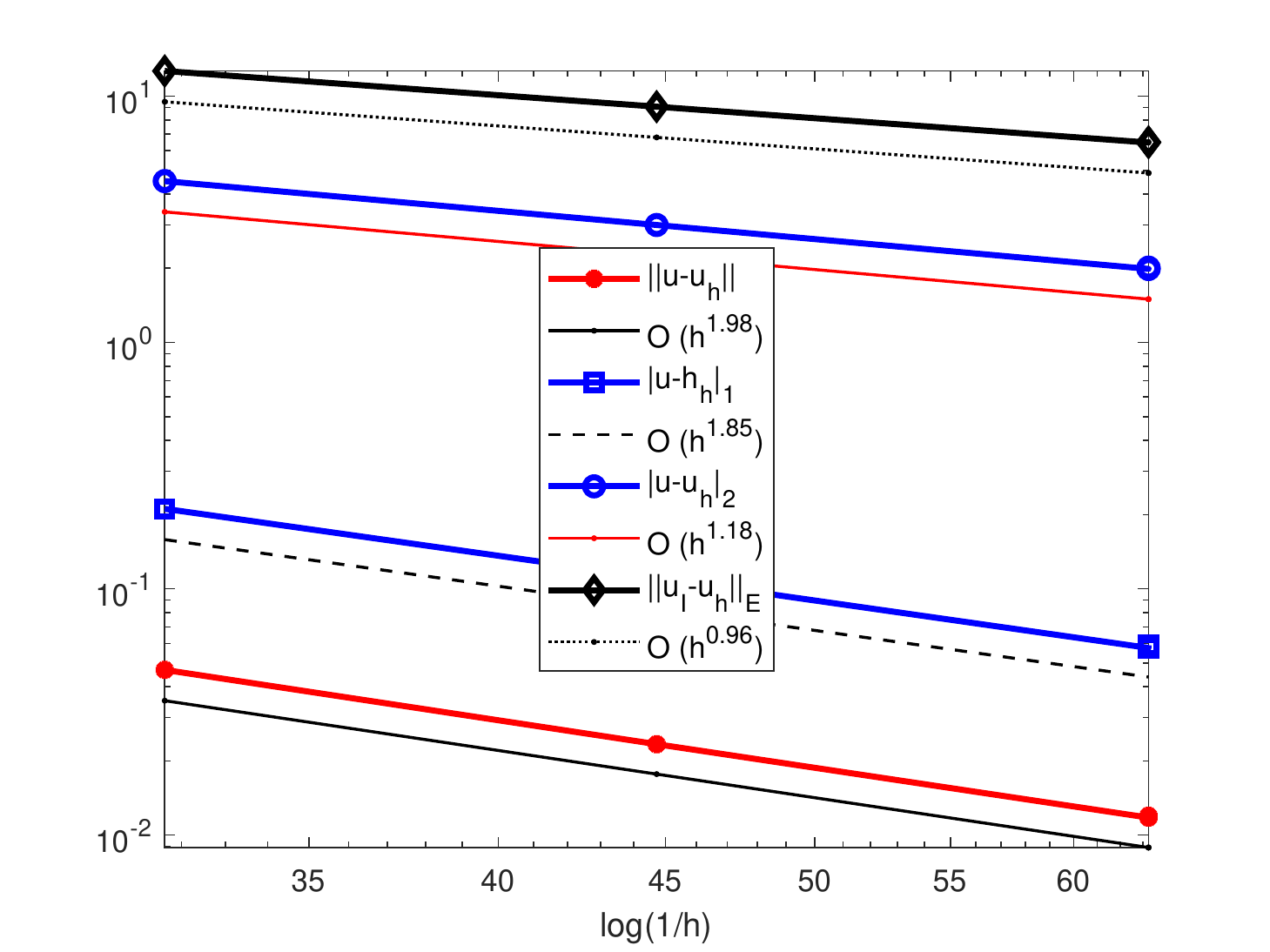}\\
  \caption{Convergence rates of the plate bending problem for the $C^0$-continuous virtual elements}\label{fig:PlateC0}
\end{figure}

We repeat the test in Example \ref{ex:Plate} and display the result in Fig.~\ref{fig:PlateC0}. In this case, we still observe the optimal convergence rates. The test script is \mcode{main\_PlateBending\_C0VEM.m}.

\subsection{Morley-type VEMs}

The fully $H^2$-nonconforming virtual element method was proposed in \cite{Antonietti-Manzini-Verani-2018} for biharmonic problem.
The local virtual element space is defined by
\begin{align*}
  V_k^2(K)
  & = \{ v\in H^2(K):  \Delta^2v\in \mathbb{P}_{k-4}(K), \\
  &  \hspace{3cm} M_{\boldsymbol{nn}}(v)|_e\in \mathbb{P}_{k-2}(e),~~Q_{\boldsymbol n}(v)|_e\in \mathbb{P}_{k-3}(e),~~e\subset \partial K \}.
\end{align*}
The corresponding degrees of freedom can be chosen as:
\begin{itemize}
  \item The values at the vertices of $K$.
  \item The moments of $v$ on edges up to degree $k-3$,
  \[\chi _e(v) = |e|^{-1}(m_e,v)_e,\quad m_e \in \mathbb{M}_{k - 3}(e).\]
  \item The moments of $\partial_{\boldsymbol n}v$ on edges up to degree $k-2$,
  \[\chi _{n_e}(v) = (m_e, \partial_{\boldsymbol n}v)_e,\quad m_e \in \mathbb{M}_{k - 2}(e).\]
  \item The moments on element $K$ up to degree $k-4$,
  \[\chi _K(v) = |K|^{-1}(m_K, v)_K,\quad {m_K} \in \mathbb{M}_{k - 4}(K).\]
\end{itemize}
If $k=2$ and $K$ is a triangle, one can prove that $V_k^2(K)$ becomes the well-known Morley element. We remark that the above Morley-type virtual element is also proposed in \cite{Zhao-Zhang-Chen-2018} with the same d.o.f.s but different local space, where the enhancement technique in \cite{Ahmad-Alsaedi-Brezzi-2013} is utilized to modify the $C^0$-continuous spaces.

\begin{figure}[!htb]
  \centering
  \includegraphics[scale=0.45]{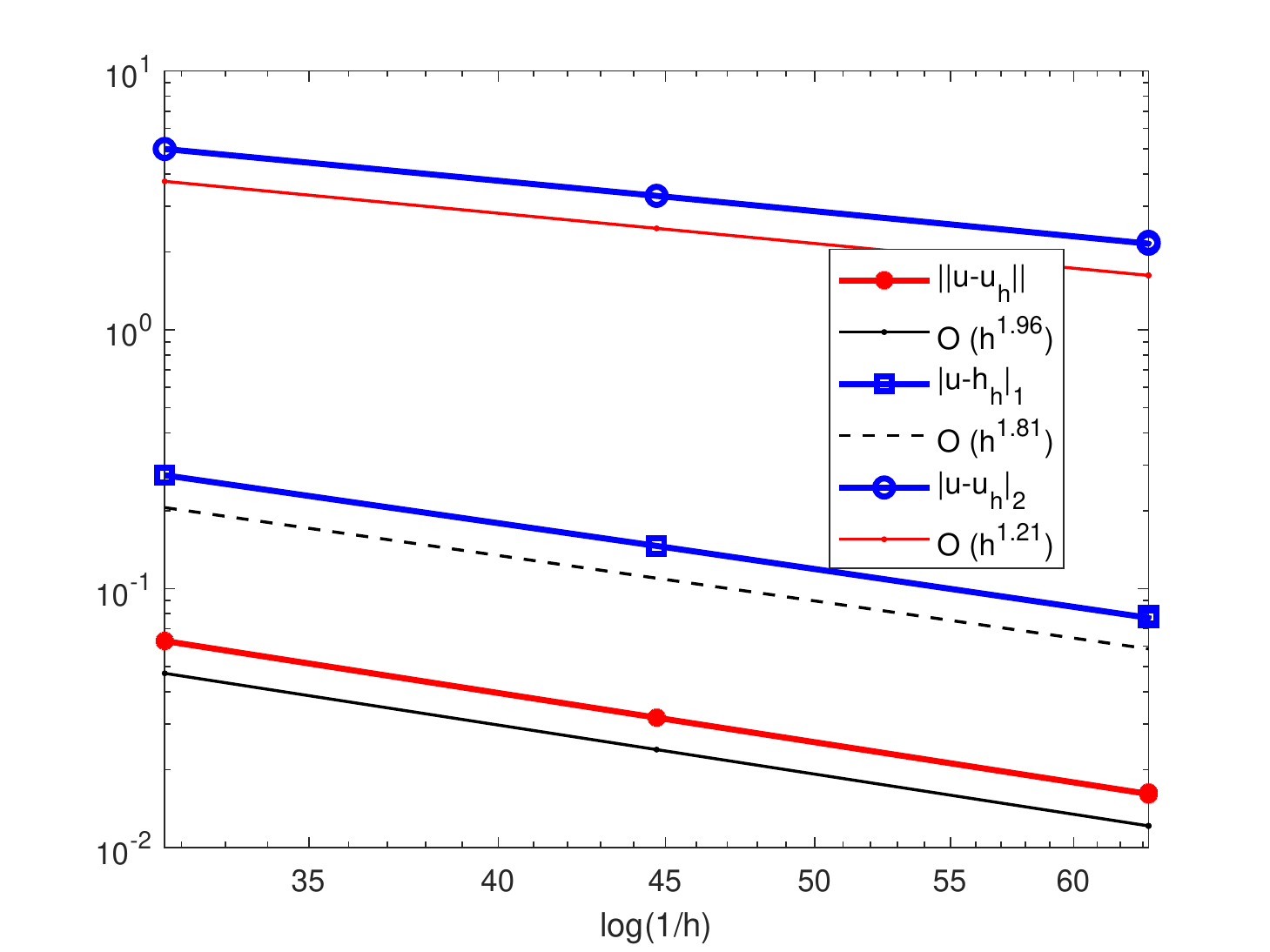}\\
  \caption{Convergence rates of the plate bending problem for the fully nonconforming virtual elements}\label{fig:PlateMorley}
\end{figure}

The test script is \mcode{main\_PlateBending\_MorleyVEM.m}. We repeat the test in Example \ref{ex:Plate} and display the result in Fig.~\ref{fig:PlateMorley}, from which we again observe the optimal rate of convergence in all the discrete norms.

\section{Stokes problem}

The Stokes problem with homogeneous Dirichlet boundary conditions is to find $(\bb{u}, p)$ such that
\[
\begin{cases}
-\nu \Delta \bb{u} - \nabla p = \bb{f} \quad & \mbox{in}~~\Omega, \\
{\rm div} \bb{u} = 0 \quad & \mbox{in}~~\Omega, \\
\bb{u} = \bb{0} \quad & \mbox{on}~~\partial \Omega.
\end{cases}
\]
Define
\[\bb{V} = \bb{H}_0^1(\Omega), \quad  Q = L_0^2(\Omega).\]
The mixed variational problem is: Find $(\bb{u}, p) \in \bb{V} \times Q$ such that
\[
\begin{cases}
a(\bb{u},\bb{v}) + b(\bb{v}, p) & = (\bb{f}, \bb{v}), \quad \bb{v}\in \bb{V}, \\
b(\bb{u},q) & = 0, \quad q\in Q,
\end{cases}
\]
where
\[a(\bb{u},\bb{v}) = (\nu \nabla \bb{u}, \nabla \bb{v}), \qquad b(\bb{v}, q) = ({\rm div}\bb{v}, q).\]

\subsection{The mixed VEMs}

In this subsection, we review the divergence-free virtual elements proposed in \cite{Beirao-Lovadina-Vacca-2017} for the Stokes problem.

\subsubsection{Virtual element space and elliptic projection}

The local virtual element space associated with $\bb{V}$ for $k\ge 2$ is
\begin{align*}
  V_k(K)
  & = \Big \{\bb{v}\in \bb{H}^1(K): \bb{v}|_{\partial K} \in [\mathbb{B}_k (\partial K)]^2, \\
  &  \hspace{3cm}
  \begin{cases}
  -\nu \Delta \bb{v} - \nabla s \in \mathcal{G}_{k-2}(K)^\bot, \\
  {\rm div} \bb{v} \in \mathbb{P}_{k-1}(K),
  \end{cases}
  \quad
  \mbox{for some $s \in L^2(K)$}
  \Big\},
\end{align*}
where
\[\mathcal{G}_{k-2}(K) = \nabla (\mathbb{P}_{k-1}(K)) \subset (\mathbb{P}_{k-2}(K))^2.\]
The finite dimensional space for $Q$ is taken as $Q_k(K) = \mathbb{P}_{k-1}(K)$.

The d.o.f.s for $V_k(K)$ and $Q_k(K)$ can be chosen as
\begin{itemize}
  \item $\mathbb{D_V}1$: the values at the vertices of $K$.
  \item $\mathbb{D_V}2$: $k-1$ interior points on each edge $e$.
  \item $\mathbb{D_V}3$: the moments
  \[\int_K \bb{v}\cdot \bb{g}_{k-2}^\bot {\rm d}x, \quad \bb{g}_{k-2}^\bot \in \mathcal{G}_{k-2}(K)^\bot. \]
  \item $\mathbb{D_V}4$: the moments
  \[\int_K ({\rm div} \bb{v})q_{k-1} {\rm d}x, \quad q_{k-1} \in \mathbb{P}_{k-1}(K)/\mathbb{R}. \]
  \item $\mathbb{D_Q}$: the moments
  \[\int_K q p_{k-1} {\rm d}x, \quad p_{k-1} \in \mathbb{P}_{k-1}(K). \]
\end{itemize}

As usual, we define the elliptic projection $\Pi^K: V_k(K)\to (\mathbb{P}_k(K))^2$, $\bb{v}\mapsto \Pi^K\bb{v}$ satisfying
\begin{equation}\label{ellipticprojStokes}
\begin{cases}
a^K(\Pi^K\bb{v}, \bb{q}) = a^K(\bb{v}, \bb{q}), \quad \bb{q}\in (\mathbb{P}_k(K))^2, \\
P_0^K(\Pi^K\bb{v}) = P_0^K(\bb{v}),
\end{cases}
\end{equation}
where $P_0^K$ is the constant projector on $K$. One can check that the elliptic projection $\Pi^K\bb{v}$ is computable by using the given d.o.f.s. In fact, the integration by parts gives

\[a^K(\bb{v}, \bb{q})  = \int_K \nu \nabla \bb{q} \cdot \nabla\bb{v}{\rm d}x
= -\int_K \nu \Delta \bb{q} \cdot \bb{v}{\rm d}x +  \int_{\partial K} (\nu \nabla \bb{q} \bb{n}) \cdot \bb{v}{\rm d}s,\]
where $\bb{n} = (n_1,n_2)^T$ and
\[\nabla \bb{q} = \begin{bmatrix} \nabla q_1 \\ \nabla q_2 \end{bmatrix} =
\begin{bmatrix} q_{1,x} & q_{1,y} \\ q_{2,x} & q_{2,y} \end{bmatrix}.\]
Since $\nu \Delta \bb{q}\in (\mathbb{P}_{k-2}(K))^2$, there exists $q_{k-1}\in \mathbb{P}_{k-1}(K)/\mathbb{R}$ and $\bb{g}_{k-2}^\bot \in \mathcal{G}_{k-2}(K)^\bot$ such that
\[\nu \Delta \bb{q} = \nabla q_{k-1} + \bb{g}_{k-2}^\bot.\]
Then
\begin{align}
a^K(\bb{v}, \bb{q})
& = -\int_K \nabla q_{k-1} \cdot \bb{v}{\rm d}x - \int_K \bb{g}_{k-2}^\bot \cdot \bb{v}{\rm d}x +  \int_{\partial K} (\nu \nabla \bb{q} \bb{n}) \cdot \bb{v}{\rm d}s \nonumber\\
& = \int_K q_{k-1}{\rm div}\bb{v} {\rm d}x - \int_K \bb{g}_{k-2}^\bot \cdot \bb{v}{\rm d}x +  \int_{\partial K} (\nu \nabla \bb{q} \bb{n}-q_{k-1}\bb{n}) \cdot \bb{v}{\rm d}s. \label{integrationstokes}
\end{align}
Obviously, the first term and the second term are determined by $\mathbb{D_V}4$ and $\mathbb{D_V}3$, respectively, while the boundary term can be computed by $\mathbb{D_V}1$ and $\mathbb{D_V}2$. Note that $\bb{g}_{k-2}^\bot = \bb{0}$ in the lowest order case $k=2$.

\subsubsection{The discrete problem}

The bilinear form $a^K(\bb{u},\bb{v})$ is approximated by
\[a_h^K(\bb{u},\bb{v}) = a^K({\Pi}^K\bb{u},{\Pi}^K\bb{v}) + S^K(\bb{u}-{\Pi}^K\bb{u},\bb{v}-{\Pi}^K\bb{v}),\]
with the stabilization term $S^K$ given by the $l^2$-inner product of the d.o.f. vectors.

In what follows we use $V_h$ to denote the global virtual element space of $V_k(K)$. The discrete space $Q_h$ of $Q$ is
\[Q_h:= \{q\in Q: q|_K\in \mathbb{P}_{k-1}(K),~~K\in \mathcal{T}_h.\}\]
The discrete mixed problem is: Find $(\boldsymbol{u}_h, p_h)\in V_h^g\times Q_h$ such that
\begin{equation}\label{discretemixedStokes}
\begin{cases}
a_h(\bb{u}_h,\bb{v}_h) + b(\bb{v}_h, p_h) & = (\bb{f}_h, \bb{v}_h), \quad \bb{v}_h\in V_h, \\
b(\bb{u}_h,q_h) & = 0, \quad q_h\in Q_h.
\end{cases}
\end{equation}
The constraint $\int_\Omega p_h{\rm d}x$ is not naturally imposed in the above system. To this end, we introduce a Lagrange multiplier and consider the augmented variational formulation: Find $((\boldsymbol{u}_h, p_h),\lambda)\in V_h^g\times Q_h \times \mathbb{R}$ such that
\begin{equation}\label{augdiscretemixedStokes}
\begin{cases}
a_h(\bb{u}_h,\bb{v}_h) + b(\bb{v}_h, p_h) & = (\bb{f}_h, \bb{v}_h), \quad \bb{v}_h\in V_h, \\
\displaystyle  b(\bb{u}_h,q_h) + \lambda \int_\Omega q_h{\rm d}x & = 0, \quad q_h\in Q_h,\\
\displaystyle \mu\int_\Omega p_h{\rm d}x & = 0, \quad \mu\in \mathbb{R}.
\end{cases}
\end{equation}

Let $\bb{\phi}_i$, $i=1,\cdots, N$ be the basis functions of $V_h$. Then
\[\bb{u} = \sum\limits_{i=1}^N \chi_i(\bb{u})\bb{\phi}_i =: \bb{\phi}^T\bb{\chi}(\bb{u}).\]
Similarly, the basis functions of $Q_h$ are denote by $\psi_l$, $l=1,\cdots,M$, with
\[p_h = \sum\limits_{l=1}^Mp_l\psi_l.\]
Plugging these expansions in \eqref{discretemixedDarcy} and taking $\bb{v}_h = \bb{\phi}_j$ and  $q_h = \psi_l$, we obtain
\[
\begin{cases}
\sum\limits_{i=1}^Na_h(\bb{\phi}_i,\bb{\phi}_j)\chi_i + \sum\limits_{l=1}^M b(\bb{\phi}_j, \psi_l)p_l & = (\bb{f}_h, \bb{\phi}_j), \quad j=1,\cdots,N, \\
\displaystyle \sum\limits_{i=1}^Nb(\bb{\phi}_i,\psi_l)\chi_i + \lambda \int_\Omega \psi_l {\rm d}x & = 0, \quad l=1,\cdots,M,\\
\displaystyle \sum\limits_{l=1}^M \int_\Omega \psi_l {\rm d}x p_l & = 0.
\end{cases}
\]
Let
\[d_l = \int_\Omega \psi_l {\rm d}x, \quad \bb{d} = [d_1,\cdots, d_M]^T.\]
The linear system can be written in matrix form as
\begin{equation}\label{linearsystemStokes}
\begin{bmatrix}
A   & B & \bb{0}\\
B^T & O & \bb{d}\\
\bb{0}^T & \bb{d}^T  & 0
\end{bmatrix}
\begin{bmatrix}
\bb{\chi} \\
\bb{p} \\
\lambda
\end{bmatrix}
=\begin{bmatrix}
\bb{f} \\
\bb{0} \\
0
\end{bmatrix},
\end{equation}
where
\[A=(a_h(\bb{\phi}_j,\bb{\phi}_i)), \quad  B = (b(\bb{\phi}_j, \psi_l)),\quad \bb{f} = ((\bb{f}_h, \bb{\phi}_j)).\]

We only consider the lowest order case $k=2$. In this case, the third-type d.o.f.s vanish and the boundary part corresponds to a tensor-product space. We arrange the d.o.f.s in the following order:
\begin{align*}
&\chi_i (\bb{v}) = \bb{v}_1(z_i),   \quad i = 1,\cdots, N_v,   \\
&\chi_{N_v+i} (\bb{v}) = \bb{v}_1(m_i),   \quad i = 1,\cdots, N_v,   \\
&\chi_{2*N_v+i} (\bb{v}) = \bb{v}_2(z_i),   \quad i = 1,\cdots, N_v,   \\
&\chi_{3N_v+i} (\bb{v}) = \bb{v}_2(m_i),   \quad i = 1,\cdots, N_v,   \\
&\chi_{4N_v+1} (\bb{v}) = \int_K {\rm div}\bb{v} m_2(x,y){\rm d}x, \quad m_2(x,y) = \frac{x-x_K}{h_K}, \\
&\chi_{4N_v+2} (\bb{v}) = \int_K {\rm div}\bb{v} m_3(x,y){\rm d}x, \quad m_3(x,y) = \frac{y-y_K}{h_K}.  \\
\end{align*}

Denote the basis functions of $\mathbb{B}_k(\partial K)$ by $\phi_1, \cdots, \phi_{N_v}; \phi_{N_v+1}, \cdots, \phi_{2N_v}$. Then the tensor-product space $(\mathbb{B}_k(\partial K))^2$ has the basis functions:
\[\overline{\phi}_1, \cdots, \overline{\phi}_{2N_v},  \underline{\phi}_{1}, \cdots, \underline{\phi}_{2N_v},\]
which correspond to the first $4N_v$ d.o.f.s. For convenience, we denote these functions by $\bb{\phi}_1, \bb{\phi}_2, \cdots, \bb{\phi}_{4N_v}$. The basis functions associated with the last two d.o.f.s are then denoted by $\bb{\varphi}_1 = \bb{\phi}_{4N_v+1}, \bb{\varphi}_2 = \bb{\phi}_{4N_v+2}$.

In the following, we only provide the details of computing the elliptic projection.

\subsection{Computation of the elliptic projection}

\subsubsection{Transition matrix}

We rewrite the basis of $V_k(K)$ in a compact form as
\[\bb{\phi}^T = (\bb{\phi}_1, \bb{\phi}_2, \cdots, \bb{\phi}_{N_k}),\]
where $N_k = 4N_v+2$ is the number of the d.o.f.s.  The basis of $(\mathbb{P}_k(K))^2$ can be denote by
\[\bb{m}^T = (\bb{m}_1, \bb{m}_2, \cdots, \bb{m}_{N_p})= (\overline{m}_1, \cdots, \overline{m}_6, \underline{m}_1, \cdots, \underline{m}_6).\]
Since $(\mathbb{P}_k(K))^2 \subset V_k(K)$, one has $\bb{m}^T = \bb{\phi}^T\bb{D}$, where $\bb{D}$ is referred to as the transition matrix from $(\mathbb{P}_k(K))^2$ to $V_k(K)$. Let
\[\boldsymbol{m}^T = [{\overline m_1},{\overline m_2},{\overline m_3},{\underline m_1},{\underline m_2},{\underline m_3}] = :[\overline m^T,\underline m^T],\]
\[{\boldsymbol \phi^T} = [{\overline \phi_1}, \cdots ,{\overline \phi_{{N}}},{\underline \phi_1}, \cdots ,{\underline \phi_{{N}}}, \bb{\varphi}_1, \bb{\varphi}_2] = :[\overline \phi^T ,\underline \phi^T , \bb{\varphi}_1, \bb{\varphi}_2].\]
Blocking the matrix $\bb{D}$ as $\bb{D} = \begin{bmatrix} \bb{D}_1 \\ \bb{D}_2 \end{bmatrix}$, one has
\begin{align*}
\bb{m} = [\overline m^T,\underline m^T]
& = [\overline \phi^T ,\underline \phi^T , 0, 0] \bb{D} + [\bb{0}^T ,\bb{0}^T  , \bb{\varphi}_1, \bb{\varphi}_2] \bb{D} \\
& = [\overline \phi^T ,\underline \phi^T , 0, 0] \begin{bmatrix} \bb{D}_1 \\ \bb{O} \end{bmatrix}
   + [\bb{0}^T ,\bb{0}^T  , \bb{\varphi}_1, \bb{\varphi}_2] \begin{bmatrix} \bb{O} \\ \bb{D}_2 \end{bmatrix}.
\end{align*}
Let $m^T = \phi^TD$. One easily finds that
\[\bb{D}_1 = \begin{bmatrix} D & \\ & D \end{bmatrix}, \qquad
\bb{D}_2 =
\begin{bmatrix}
\chi_{4N_v+1}(\overline{m}^T) & \chi_{4N_v+1}(\underline{m}^T)\\
\chi_{4N_v+2}(\overline{m}^T) & \chi_{4N_v+2}(\underline{m}^T)
\end{bmatrix}\]

We first provide some necessary information.
\vspace{-0.8cm}
\begin{lstlisting}
    % ------- element information ----------
    index = elem{iel};     Nv = length(index);
    xK = centroid(iel,1); yK = centroid(iel,2); hK = diameter(iel);
    x = node(index,1); y = node(index,2);
    v1 = 1:Nv; v2 = [2:Nv,1]; % loop index for vertices or edges
    xe = (x(v1)+x(v2))/2;  ye = (y(v1)+y(v2))/2; % mid-edge points
    Ne = [y(v2)-y(v1), x(v1)-x(v2)]; % he*ne
    nodeT = [node(index,:);centroid(iel,:)];
    elemT = [(Nv+1)*ones(Nv,1),(1:Nv)',[2:Nv,1]'];

    % --------------- scaled monomials -----------------
    m1 = @(x,y) 1+0*x;                  gradm1 = @(x,y) [0+0*x, 0+0*x];
    m2 = @(x,y) (x-xK)./hK;             gradm2 = @(x,y) [1/hK+0*x, 0+0*x];
    m3 = @(x,y) (y-yK)./hK;             gradm3 = @(x,y) [0+0*x, 1/hK+0*x];
    m4 = @(x,y) (x-xK).^2/hK^2;         gradm4 = @(x,y) [2*(x-xK)./hK^2, 0+0*x];
    m5 = @(x,y) (x-xK).*(y-yK)./hK^2;   gradm5 = @(x,y) [(y-yK)./hK^2, (x-xK)./hK^2];
    m6 = @(x,y) (y-yK).^2./hK^2;        gradm6 = @(x,y) [0+0*x, 2*(y-yK)./hK^2];

    m = @(x,y) [m1(x,y), m2(x,y), m3(x,y), m4(x,y), m5(x,y), m6(x,y)];
    Gradm = {gradm1, gradm2, gradm3, gradm4, gradm5, gradm6};
    divmm = @(x,y) [0+0*x, 1/hK+0*x, 0+0*x, 2*(x-xK)./hK^2, (y-yK)./hK^2, 0+0*x, ...
        0+0*x, 0+0*x, 1/hK+0*x, 0+0*x, (x-xK)./hK^2, 2*(y-yK)./hK^2];
\end{lstlisting}
Then the transition matrix can be realized as
\vspace{-0.8cm}
\begin{lstlisting}
    % -------- transition matrix ----------
    NdofBd = 2*Nv; NdofA = 2*NdofBd+2;
    divmm2 = @(x,y) divmm(x,y).*repmat(m2(x,y),1,Nmm);
    divmm3 = @(x,y) divmm(x,y).*repmat(m3(x,y),1,Nmm);
    D = zeros(NdofA, Nmm);
    Dbd = [m(x,y); m(xe,ye)];
    D(1:4*Nv, :) = blkdiag(Dbd, Dbd);
    D(end-1,:) = integralTri(divmm2,4,nodeT,elemT);
    D(end,:) = integralTri(divmm3,4,nodeT,elemT);
\end{lstlisting}

\subsubsection{Elliptic projection matrices}

The elliptic projection satisfies
\[
\begin{cases}
\bb{G}\bb{\Pi}_*^K = \bb{B}, \\
P_0^K(\bb{m}^T)\bb{\Pi}_*^K = P_0^K(\bb{\phi}^T)
\end{cases} \quad \mbox{or} \quad
\tilde{\bb{G}}\bb{\Pi}_*^K = \tilde{\bb{B}},
\]
where
\[{\bb{G}} = a^K({\bb{m}}, {\bb{m}}^T), \quad {\bb{B}} = a^K({\bb{m}}, \bb{\phi}^T).\]

For $k=2$,
\begin{align*}
{\bb{B}}_{\alpha i}
& = a^K({\bb{m}}_\alpha, \bb{\phi}_i) = \int_K q_\alpha{\rm div}\bb{\phi}_i {\rm d}x  +  \int_{\partial K} (\nu \nabla {\bb{m}}_\alpha \bb{n}-q_\alpha\bb{n}) \cdot \bb{\phi}_i{\rm d}s \\
& =:  I_1(\alpha,i) +  I_2(\alpha,i),
\end{align*}
where
\[I_1(\alpha,i) = \int_K q_\alpha{\rm div}\bb{\phi}_i {\rm d}x, \quad I_2(\alpha,i) = \int_{\partial K} (\nu \nabla {\bb{m}}_\alpha \bb{n}-q_\alpha\bb{n}) \cdot \bb{\phi}_i{\rm d}s,\]
and
\[\nu \Delta {\bb{m}}_\alpha = \nabla q_\alpha + 0, \quad q_\alpha \in \mathbb{P}_1(K)/\mathbb{P}_0(K).\]

First consider $I_1$. By definition, let $q_\alpha = \nu h_K(c_{2,\alpha} m_2 + c_{3,\alpha} m_3)$. Then
\[\begin{bmatrix}c_{2,\alpha} \\ c_{3,\alpha} \end{bmatrix}  \longleftarrow  \begin{bmatrix} \Delta m^T & \bb{0}^T \\ \bb{0}^T & \Delta m^T\end{bmatrix}.\]
The Kronecher's property gives
\begin{align*}
I_1(\alpha,i)
& = \nu h_K \Big(c_{2,\alpha}\int_K  m_2{\rm div}\bb{\phi}_i {\rm d}x+ c_{3,\alpha} \int_K m_3 {\rm div}\bb{\phi}_i {\rm d}x\Big) \\
& = \nu h_K \Big(c_{2,\alpha}\delta_{i, (4N_v+1)}+ c_{3,\alpha} \delta_{i, (4N_v+2)}\Big).
\end{align*}
The computation of $I_1$ reads
\vspace{-0.8cm}
\begin{lstlisting}
    % --- first term ---
    Lapm = [0, 0, 0, 2/hK^2, 0, 2/hK^2];
    B = zeros(Nmm, NdofA);
    B(1:Nm, end-1) = pde.nu*hK*Lapm;
    B(Nm+1:end, end) = pde.nu*hK*Lapm;
\end{lstlisting}

For the second term $I_2$, let $(g_1^\alpha,g_2^\alpha)^T = \nu \nabla {\bb{m}}_\alpha \bb{n}-q_\alpha\bb{n}$, Then
\[I_2(\alpha,i)
 = \int_{\partial K} g_1^\alpha \bb{\phi}_{1,i} {\rm d}s + \int_{\partial K} g_2^\alpha \bb{\phi}_{2,i} {\rm d}s=:J_1(\alpha,i)+J_2(\alpha,i).\]
Let $\phi$ be the basis of $\mathbb{B}_k(\partial K)$. One has
\[J_1(\alpha,:)^T =
\begin{bmatrix}
( g_1^\alpha, \phi )_{\partial K} \\
( g_1^\alpha, \bb{0} )_{\partial K} \\
( g_1^\alpha, \bb{\varphi}_{1,1} )_{\partial K}\\
( g_1^\alpha, \bb{\varphi}_{2,1} )_{\partial K}
\end{bmatrix}
= \begin{bmatrix}
( g_1^\alpha, \phi )_{\partial K} \\
\bb{0} \\
0\\
0
\end{bmatrix},
\]
\[J_2(\alpha,:)^T =
\begin{bmatrix}
( g_2^\alpha, \bb{0} )_{\partial K} \\
( g_2^\alpha, \phi )_{\partial K} \\
( g_2^\alpha, \bb{\varphi}_{1,2} )_{\partial K}\\
( g_2^\alpha, \bb{\varphi}_{2,2} )_{\partial K}
\end{bmatrix}
= \begin{bmatrix}
\bb{0}\\
( g_2^\alpha, \phi )_{\partial K} \\
0\\
0
\end{bmatrix},
\]
where $( g_1^\alpha, \phi )_{\partial K}$ and $( g_2^\alpha, \phi )_{\partial K}$ can be computed using the  assembling technique for FEMs. For example, for $( g_1^\alpha, \phi )_{\partial K}$ there exist three integrals on $e$:
\[F_i = \int_e g_1^\alpha \phi_i {\rm d}s, \quad i = 1,2,3,\]
where $e=[a_e,m_e,b_e]$ is an edge with $a_e$ and $b_e$ being the endpoints and $m_e$ the midpoint. By the Simpson's formula,
\[F = \begin{bmatrix} F_1 \\ F_2 \\ F_3 \end{bmatrix}
    = \frac{h_e}{6}\begin{bmatrix} g_1^\alpha(a_e) \\ g_1^\alpha(b_e) \\ 4g_1^\alpha(m_e)\end{bmatrix}.\]
The above discussion can be realized as follows.
\vspace{-0.8cm}
\begin{lstlisting}
    % --- second term ---
    elem1 = [v1(:), v2(:), v1(:)+Nv]; % elem2dof for [ae, be, me]
    Gradmm = cell(2,Nmm);
    Gradmm(1,1:Nm) = Gradm;
    Gradmm(2,Nm+1:end) = Gradm;
    for im = 1:Nm
        Gradmm{1,im+Nm} = @(x,y) [0+0*x, 0+0*x];
        Gradmm{2,im} = @(x,y) [0+0*x, 0+0*x];
    end
    qmm = cell(1,Nmm);  % q = nu*hK(c2*m2 + c3*m3)
    c2 = [Lapm, zeros(1,Nm)];
    c3 = [zeros(1,Nm), Lapm];
    for im = 1:Nmm
        qmm{im} = @(x,y) pde.nu*hK*(c2(im)*m2(x,y) + c3(im)*m3(x,y));
    end
    for s = 1:2
        id = (1:NdofBd) + (s-1)*NdofBd;
        for im = 1:Nmm
            pm = @(x,y) pde.nu*Gradmm{s,im}(x,y);
            qa = @(x,y) qmm{im}(x,y);
            F1 = 1/6*(sum(pm(x(v1),y(v1)).*Ne, 2) - qa(x(v1),y(v1)).*Ne(:,s));
            F2 = 1/6*(sum(pm(x(v2),y(v2)).*Ne, 2) - qa(x(v2),y(v2)).*Ne(:,s));
            F3 = 4/6*(sum(pm(xe,ye).*Ne, 2) - qa(xe,ye).*Ne(:,s));
            B(im, id) = accumarray(elem1(:), [F1; F2; F3], [NdofBd, 1]);
        end
    end
\end{lstlisting}

We finally consider the implementation of the constraint. At first glance the $L^2$ projection is not computable since there is no zero-order moment on $K$. In fact, the computability can be obtained using the decomposition of polynomial spaces.
Let $\bb{\phi}_i = [\bb{\phi}_{1,i},\bb{\phi}_{1,i}]^T$. Then
\[\bb{\phi}_{1,i} = \bb{\phi}_i \cdot \begin{bmatrix} 1 \\ 0 \end{bmatrix}, \quad
\bb{\phi}_{2,i} = \bb{\phi}_i \cdot \begin{bmatrix} 0 \\ 1 \end{bmatrix}.\]
It is easy to get
\[\begin{bmatrix} 1 \\ 0 \end{bmatrix}
= \nabla p_{k-1} + \bb{g}_{k-2}^\bot, \quad
p_{k-1} = h_K  m_2 , \quad \bb{g}_{k-2}^\bot = \bb{0},\]
\[\begin{bmatrix} 0 \\ 1 \end{bmatrix}
= \nabla q_{k-1} + \bb{g}_{k-2}^\bot, \quad
q_{k-1} = h_K  m_3 , \quad \bb{g}_{k-2}^\bot = \bb{0},\]
which yield
\begin{align}
& P_0^K(\bb{\phi}_{1,i})
 = |K|^{-1} \int_K \bb{\phi}_i \cdot \nabla p_{k-1} {\rm d}x
 =  |K|^{-1}h_K\Big(- \int_K {\rm div} \bb{\phi}_i  m_2 {\rm d}x + \int_{\partial K} m_2 \bb{\phi}_i \cdot \bb{n} {\rm d}s\Big), \nonumber\\
& P_0^K(\bb{\phi}_{2,i})
 = |K|^{-1} \int_K \bb{\phi}_i \cdot \nabla q_{k-1} {\rm d}x
 = |K|^{-1}h_K\Big(- \int_K {\rm div} \bb{\phi}_i  m_3 {\rm d}x + \int_{\partial K} m_3\bb{\phi}_i \cdot \bb{n} {\rm d}s\Big). \label{constraintStokes}
\end{align}
As you can see, their computation is similar to the previous one for $B_{\alpha, i}$ with $\alpha$ fixed. The resulting two row vectors will replace the first and seventh rows of $B$.
\vspace{-0.8cm}
\begin{lstlisting}
    % constraint
    P0K = zeros(2,NdofA);
    P0K(1, end-1) = -1; P0K(2, end) = -1;
    m23 = {m2, m3};
    for s = 1:2
        mc = m23{s};
        F1 = 1/6*(mc(x(v1),y(v1)).*Ne);  % [n1, n2]
        F2 = 1/6*(mc(x(v2),y(v2)).*Ne);
        F3 = 4/6*(mc(xe,ye).*Ne);
        F = [F1; F2; F3];
        P0K(s, 1:NdofBd) = accumarray(elem1(:), F(:,1), [NdofBd 1]);
        P0K(s, NdofBd+1:2*NdofBd) = accumarray(elem1(:), F(:,2), [NdofBd 1]);
    end
    P0K = 1/area(iel)*hK*P0K;
    % Bs, G, Gs
    Bs = B;  Bs([1,7], :) = P0K;
\end{lstlisting}

\begin{example}
Let $\Omega = (0,1)^2$. We choose the load term $\bb{f}$ in such a way that the analytical solution is
\[\bb{u}(x,y) =
\begin{bmatrix}
-\frac{1}{2}\cos(x)^2 \cos(y) \sin(y)\\
\frac{1}{2} \cos(y)^2 \cos(x) \sin(x)
 \end{bmatrix},
 \quad
 p(x,y) = \sin(x) - \sin(y).\]
\end{example}

\begin{figure}[!htb]
  \centering
  \includegraphics[scale=0.45]{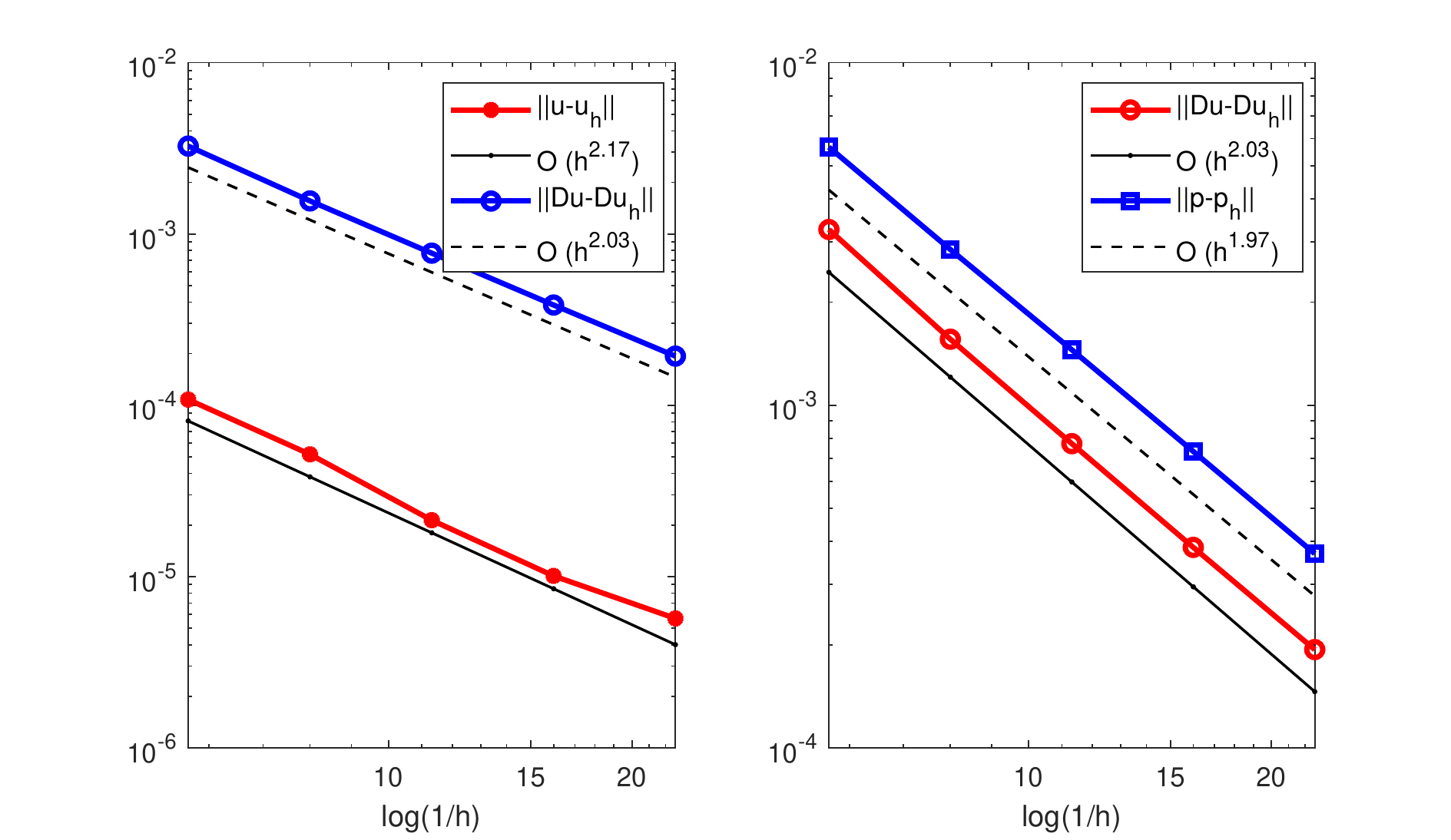}\\
  \caption{Convergence rates of the divergence free VEM for the Stokes problem}\label{fig:Stokes}
\end{figure}

The results are displayed in Fig.~\ref{fig:Stokes} and Tab.~\ref{tab:Stokes}, from which we observe the optimal rates of convergence both for $u$ and $p$.

\begin{table}[!htb]
 \centering
 \caption{The discrete errors for the Stokes problem} \label{tab:Stokes}
\begin{tabular}{rcccccccccc}
  \hline
  $\sharp$Dof   &    $h$    &    $\|u-u_h\|$ &   $|u-u_h|_1$  &  $ \|p-p_h\|$ \\
  \hline
 390  & 1.768e-01  & 1.07945e-04   & 3.25802e-03  & 5.66765e-03  \\
 774  & 1.250e-01  & 5.16780e-05   & 1.55844e-03  & 2.84772e-03  \\
1534  & 8.839e-02  & 2.13132e-05   & 7.72427e-04  & 1.45175e-03  \\
3042  & 6.250e-02  & 1.00841e-05   & 3.84674e-04  & 7.33126e-04  \\
6090  & 4.419e-02  & 5.69480e-06   & 1.93534e-04  & 3.68775e-04  \\
  \hline
\end{tabular}
\end{table}

\section{Adaptive virtual element methods}

Due to the large flexibility of the meshes, researchers have focused on the a posterior error analysis of the VEMs and made some progress in recent years \cite{Beirao-Manzini-2015,Cangiani-Georgoulis-Pryer-Sutton-2017,Berrone-Borio-2017a,
Chi-Beirao-Paulino-2019,Beirao-Manzini-Mascotto-2019,Huang-Lin-2021}.

We consider the adaptive VEMs for the Poisson equation. The computable error estimator is taken from \cite{Berrone-Borio-2017a,Huang-Lin-2021}, given as
\[\eta(u_h) = \Big( \sum\limits_{K\in \mathcal{T}_h} \eta_K^2(u_h) \Big)^{1/2},\]
where
\[\eta_K^2(u_h) = \sum\limits_{i=1}^4 \eta_{i,K}^2(u_h),\]
with
\[\eta_{1,K}^2 = h_K^2 \|f-\Pi_0^K f\|_{0,K}^2, \quad
\eta_{2,K}^2 = h_K^2 \|\Pi_0^K f\|_{0,K}^2, \quad
\eta_{3,K}^2 = \|\bb{\chi}(u_h-\Pi_1^\nabla u_h)\|_{l^2}^2,\]
and
\[\eta_{4,K}^2 = \frac{1}{2} \sum\limits_{e\subset \partial K}h_e \| [\partial_{\bb{n}} \Pi_1^\nabla u_h ] \|_{0,e}^2.\]

Standard adaptive algorithms based on the local mesh refinement can be written as loops of the form
\[{\bf SOLVE} \to {\bf ESTIMATE} \to {\bf MARK} \to {\bf REFINE}.\]
Given an initial polygonal subdivision $\mathcal{T}_0$, to get $\mathcal{T}_{k+1}$ from $\mathcal{T}_k$ we
first solve the VEM problem
under consideration to get the numerical solution $u_k$ on $\mathcal{T}_k$. The error is then estimated
by using $u_k$, $\mathcal{T}_k$ and the a posteriori error bound. The local error bound is used to mark a subset
of elements in $\mathcal{T}_k$ for refinement. The marked polygons and possible more neighboring elements are refined
in such a way that the subdivision meets certain conditions, for example, the resulting polygonal mesh is still shape
regular. In the implementation, it is usually time-consuming to write a mesh refinement function since we need to
carefully design the rule for dividing the marked elements to get a refined mesh of high quality.
We have present such an implementation of the mesh refinement for polygonal meshes in \cite{PolyMeshRefine}.

\begin{figure}[!htb]
  \centering
  \subfigure[]{\includegraphics[scale=0.45]{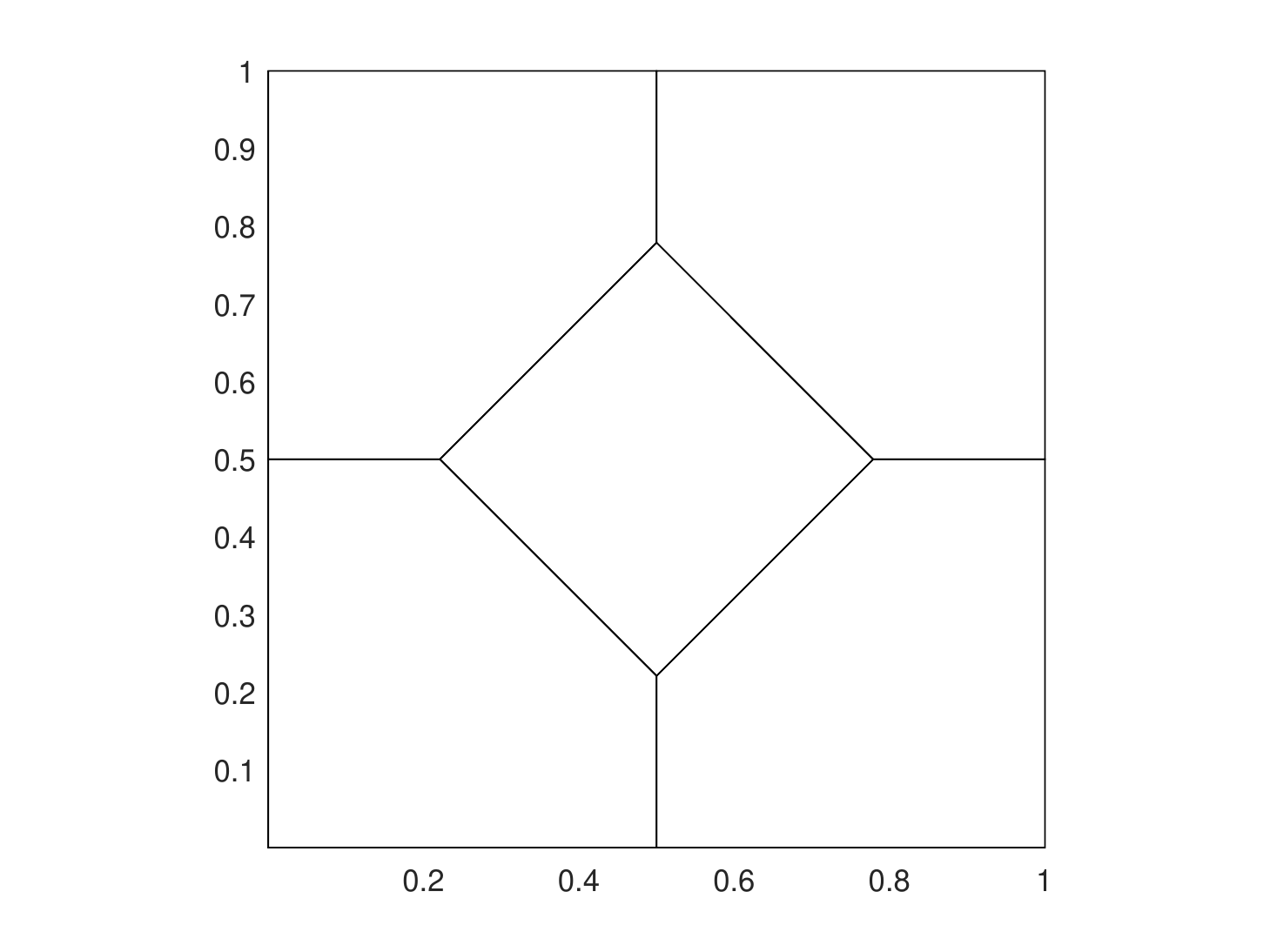}}
  \subfigure[]{\includegraphics[scale=0.45]{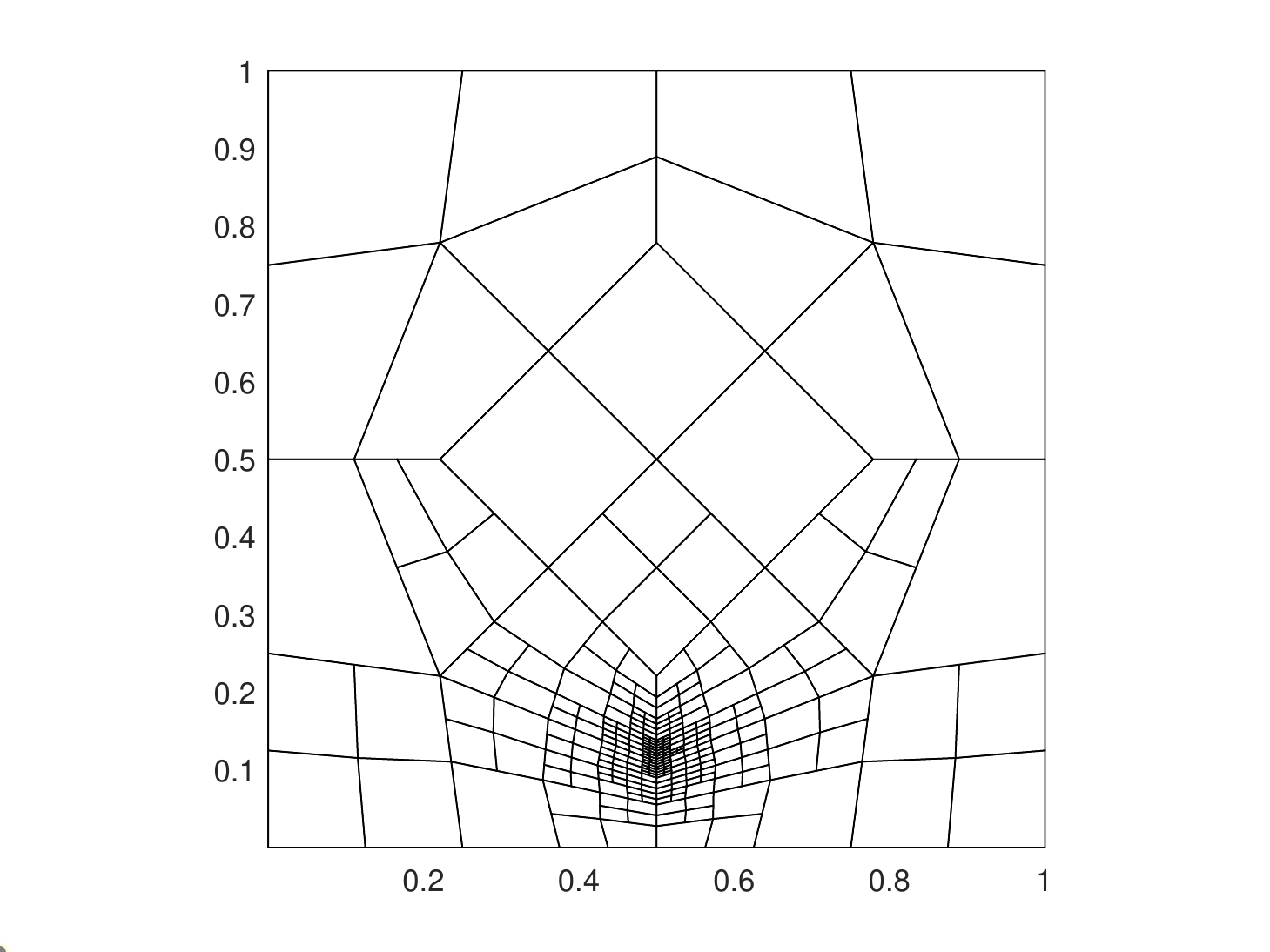}}\\
  \subfigure[]{\includegraphics[scale=0.45]{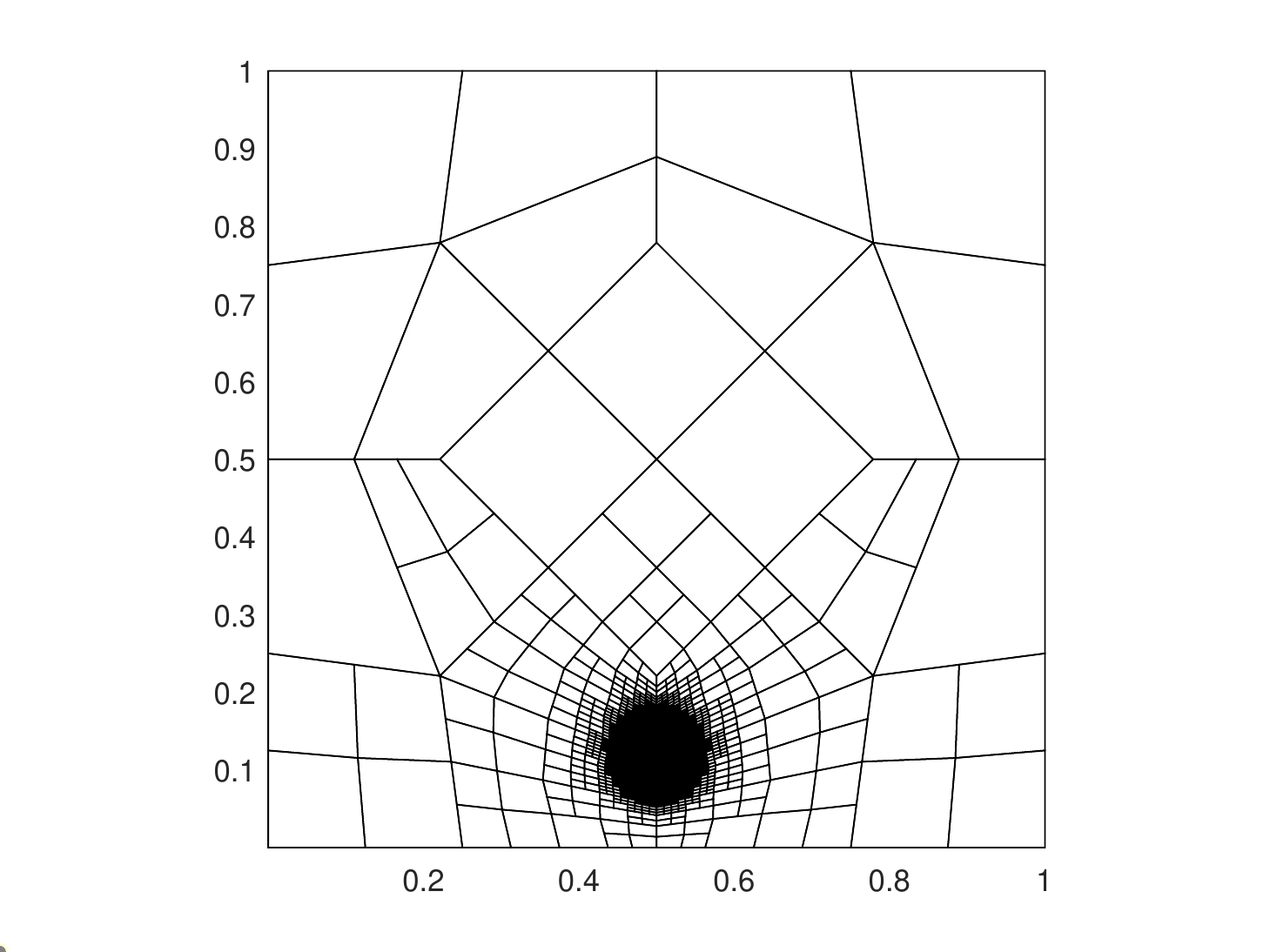}}
  \subfigure[]{\includegraphics[scale=0.45]{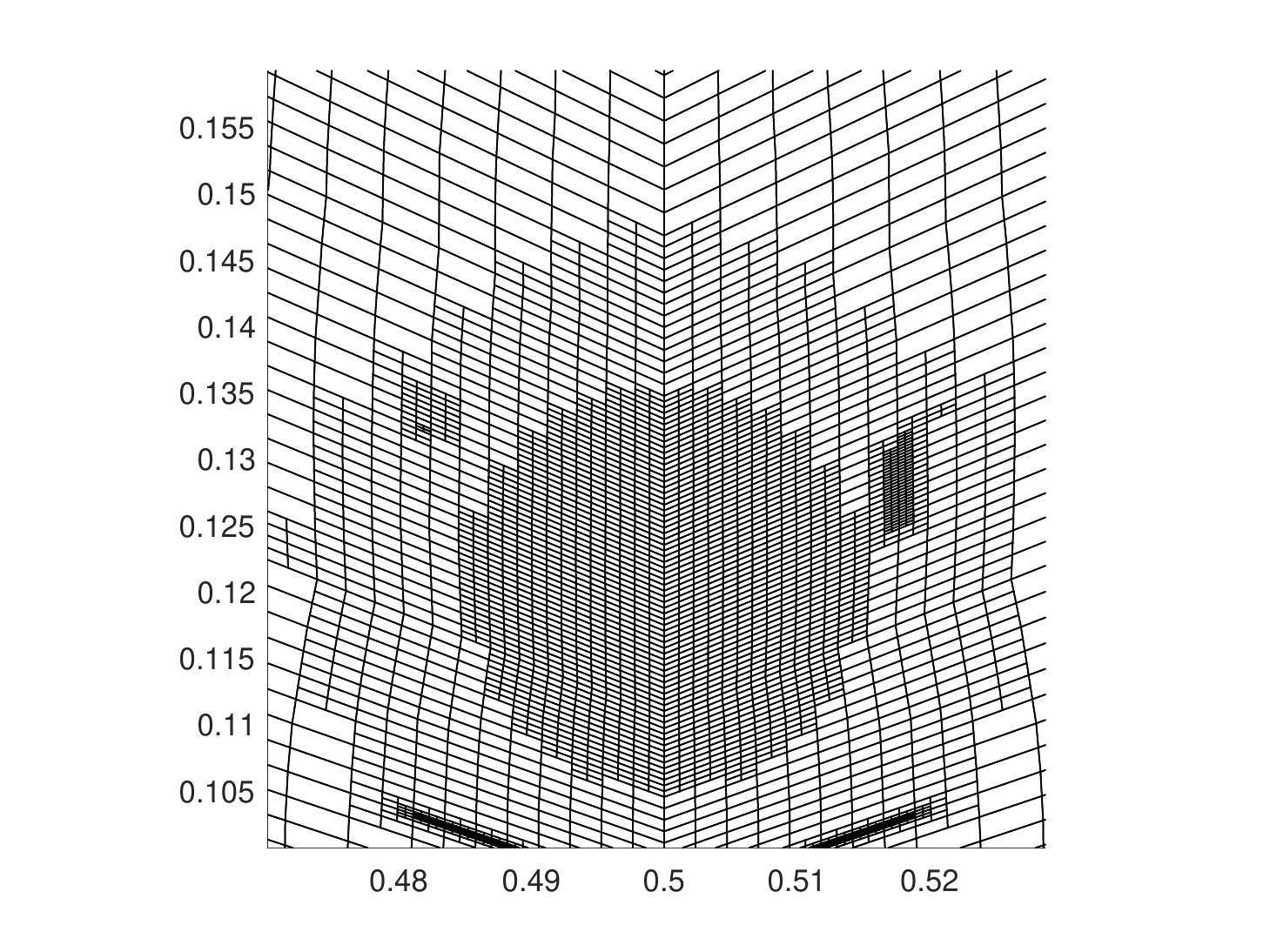}}\\
  \caption{The initial and the final adapted meshes. (a) The initial mesh;
  (b) After 20 refinement steps; (c) After 30 refinement steps; (d) The zoomed mesh in (c)}\label{refineVEMmesh}
\end{figure}

\begin{figure}[!htb]
  \centering
  \subfigure[Exact]{\includegraphics[scale=0.45]{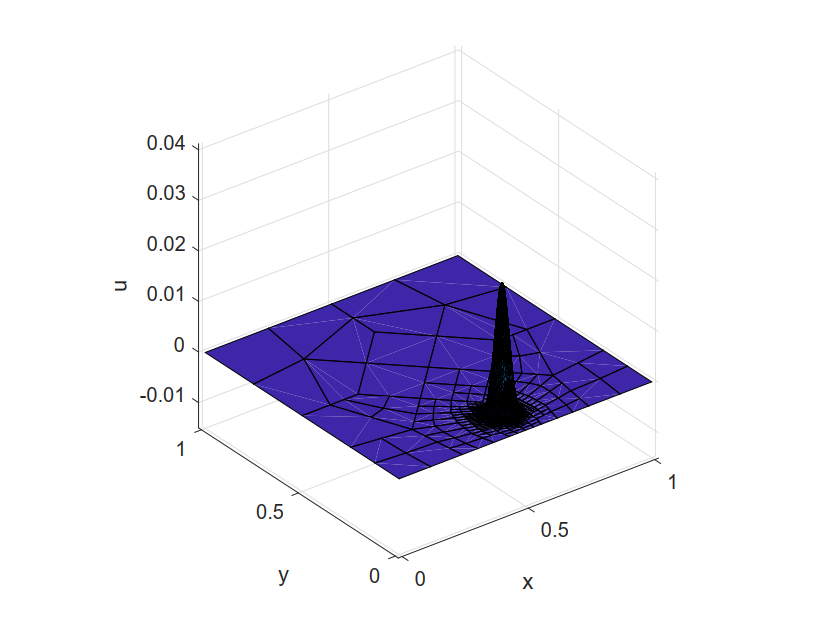}}
  \subfigure[Numerical]{\includegraphics[scale=0.45]{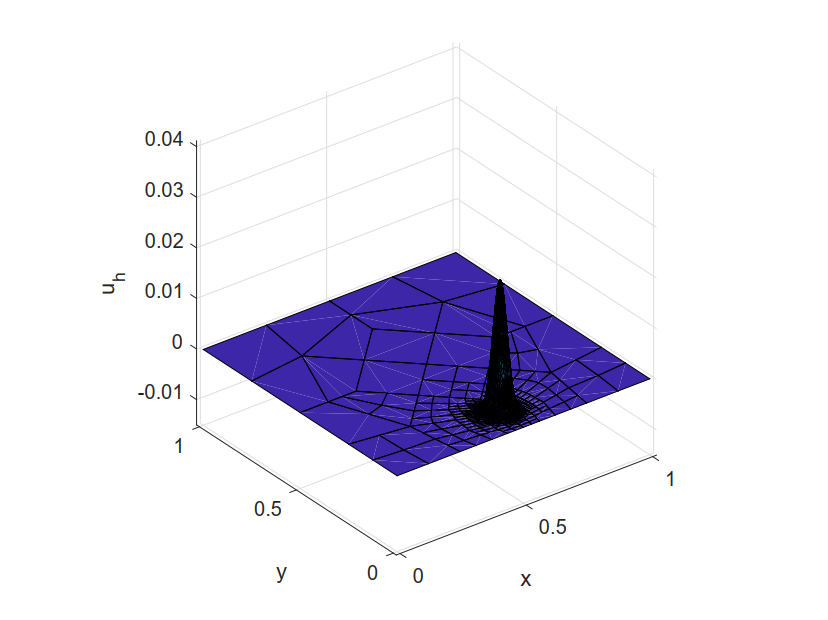}}\\
  \caption{The exact and numerical solutions}\label{refineVEMsol}
\end{figure}

\begin{figure}[!htb]
  \centering
  \includegraphics[scale=0.45]{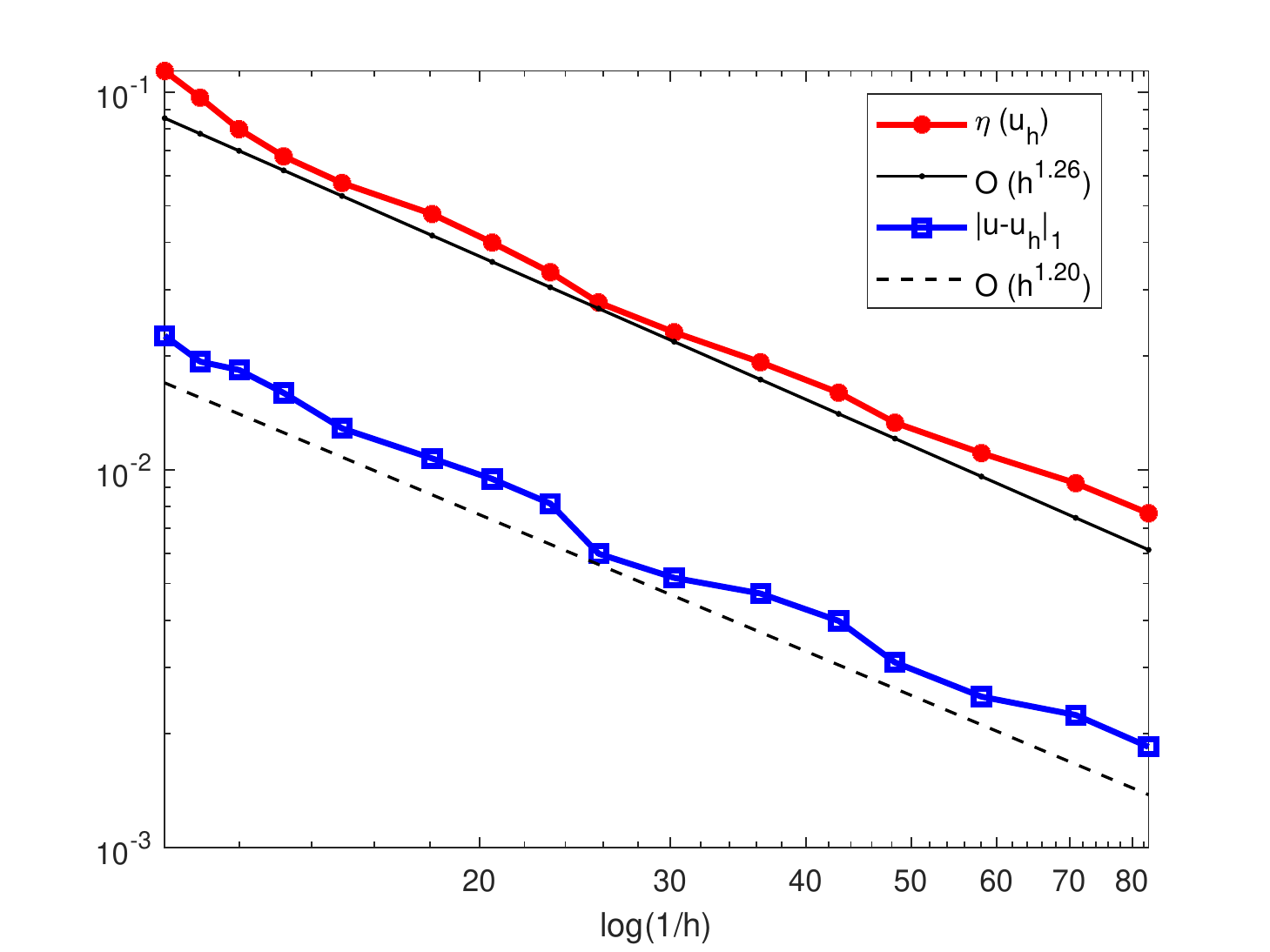}\\
  \caption{Convergence rates of the error $|u - \Pi_1^\nabla u_h|_1$ and the error estimator $\eta(u_h)$}\label{fig:adaptiveRat}
\end{figure}

Consider the Poisson equation with Dirichlet boundary condition on the unit square. The exact solution is given by
\[u(x,y) = xy(1 - x)(1 - y){\text{exp}}\left(  - 1000((x - 0.5)^2 + (y - 0.117)^2) \right).\]
We employ the VEM in the lowest order case and use the D\"{o}rfler marking strategy with parameter $\theta = 0.4$ to select the subset of elements for refinement.
The initial mesh and the final adapted meshes after 20 and 30 refinement steps are presented in Fig.~\ref{refineVEMmesh}~(a-c),
respectively. The detail of the last mesh is shown in Fig.~\ref{refineVEMmesh}~(d). Clearly, no small edges are observed. We also plot the adaptive approximation in Fig.~\ref{refineVEMsol}, which almost coincides with the exact solution. The full code is available from mVEM package. The subroutine \mcode{PoissonVEM\_indicator.m} is used to compute the local indicator and the test script is \mcode{main\_Poisson\_avem.m}. As shown in Fig.~\ref{fig:adaptiveRat}, we see the adaptive strategy correctly refines the mesh in a neighborhood of the singular point and there is a good level of agreement between the $H^1$ error and error estimator.

\section{Variational inequalities} \label{sec:VI}

We now focus on the the virtual element method to solve a simplified friction problem, which is a typical elliptic variational inequality of the second kind.

\subsection{The simplified friction problem}

Let $\Omega \subset \mathbb{R}^2$ be a bounded domain with a Lipschitz boundary $\Gamma = \partial \Omega$ that is divided into two parts $\Gamma_C$ and $\Gamma_D$. The problem is
\begin{equation}\label{VI0}
\begin{cases}
 - \Delta u + \alpha u = f \quad & \mbox{in} ~~ \Omega, \\
 \partial_{\bb{n}} u \le g, \quad
  u \partial_{\bb{n}} u + g |u| = 0, \quad & \mbox{on}~~\Gamma_C, \\
 u = 0  \quad \mbox{on}~~ & \Gamma_D,
\end{cases}
\end{equation}
where $\alpha>0$ is a constant, $f\in L^2(\Omega)$, $g\in L^2(\Gamma_C)$, $\Gamma_C$ is the frictional boundary part and $\Gamma_D$ is the Dirichlet boundary part.

Define
\[V = \{v\in H^1(\Omega): v|_{\Gamma_D} = 0\}.\]
The variational inequality is \cite{Wang-Wei-2018}: Find $u\in V$ such that
\begin{equation}\label{VI}
a(u,v-u) + j(v)-j(u) \ge \ell(v-u), \quad v\in V,
\end{equation}
where
\[a(u,v) = \int_{\Omega} (\nabla u \cdot \nabla v + \alpha u v ) {\rm d}x, \quad \ell(v) = \int_{\Omega} f v {\rm d}x, \quad
 j(v) = \int_{\Gamma_C} g |v| {\rm d}s.  \]

\subsection{The VEM discretization}

We consider the lowest-order virtual element space. The local space is taken as the enhanced virtual element space $W_1(K)$ defined in \eqref{second-VEMspace}. Let $V_h$ be the global space. The virtual element method for solving the simplified friction problem is: Find $u_h\in V_h$ such that
\begin{equation}\label{VIvem}
a_h(u_h, v_h-u_h) + j(v_h)-j(u_h)  \ge \ell_h(v_h-u_h), \quad v_h \in V_h,
\end{equation}
where
\[a_h^K(v,w) = (\nabla \Pi_1^\nabla v, \nabla \Pi_1^\nabla w)_K + \alpha (\Pi_1^0 v, \Pi_1^0 w)_K +
S^K(v-\Pi_1^\nabla v, w-\Pi_1^\nabla w),\]
\[S^K(v, w) := (1+\alpha h_K^2) \bb{\chi} (v)\cdot \bb{\chi} (w), \quad
\ell_h(v_h)  = \sum\limits_{K\in\mathcal{T}_h} (f, \Pi_0^0v_h)_K.\]

By introducing a Lagrangian multiplier
\[\lambda \in \Lambda = \{\lambda \in L^\infty(\Gamma_C): |\lambda|\le 1 \quad \mbox{a.e.~~on} ~~\Gamma_C\}, \]
the inequality problem \eqref{VI} can be rewritten as
\[\begin{cases}
a(u,v) + \displaystyle \int_{\Gamma_C} g \lambda v {\rm d}s = \ell(v), \quad & v \in V, \\
\lambda u = |u| \quad  \mbox{a.e.~~on} ~~\Gamma_C.
\end{cases}\]
For this reason, the discrete problem \eqref{VIvem} can be recast as
\[\begin{cases}
a_h(u_h,v_h) + \displaystyle \int_{\Gamma_C} g \lambda_h v_h {\rm d}s = \ell_h(v_h), \quad & v_h \in V_h, \\
\lambda_h u_h = |u_h| \quad  \mbox{a.e.~~on} ~~\Gamma_C,
\end{cases}\]
where $\lambda_h \in L^\infty(\Gamma_C)$ and $|\lambda_h| \le 1$. Then the Uzawa algorithm for solving the
above problem is \cite{Wang-Zhao-2021,Wu-Wang-Han-2022}: given any $\lambda_h^{(0)} \in \Lambda$, for $n\ge 1$, find $u_h^{(n)}$ and $\lambda_h^{(n)}$ by solving
\begin{equation}\label{Uzawa}
a_h(u_h^{(n)},v_h) = \ell_h(v_h) - \int_{\Gamma_C} g \lambda_h^{(n-1)} v_h {\rm d}s
\end{equation}
and
\[\lambda_h^{(n)} = P_{\Lambda}(\lambda_h^{(n-1)} + \rho g u_h^{(n)}),\]
where $P_{\Lambda}(\mu) = \sup \{-1, \inf \{1,\mu\}\}$ and $\rho$ is a constant parameter.

\subsection{Numerical example}

Let $\Omega = (0,1)^2$ and suppose that the frictional boundary condition is imposed on $y=0$. The function $g$ can be simply chosen as $\sup_{\Gamma_C} |\partial_{\bb{n}} u|$. The right-hand function $f$ is chosen such that the exact solution is $u = (\sin(x)-x\sin(1))\sin(2\pi y)$.

\begin{figure}[!htb]
  \centering
  \includegraphics[scale=0.45,trim = 20 80 20 80,clip]{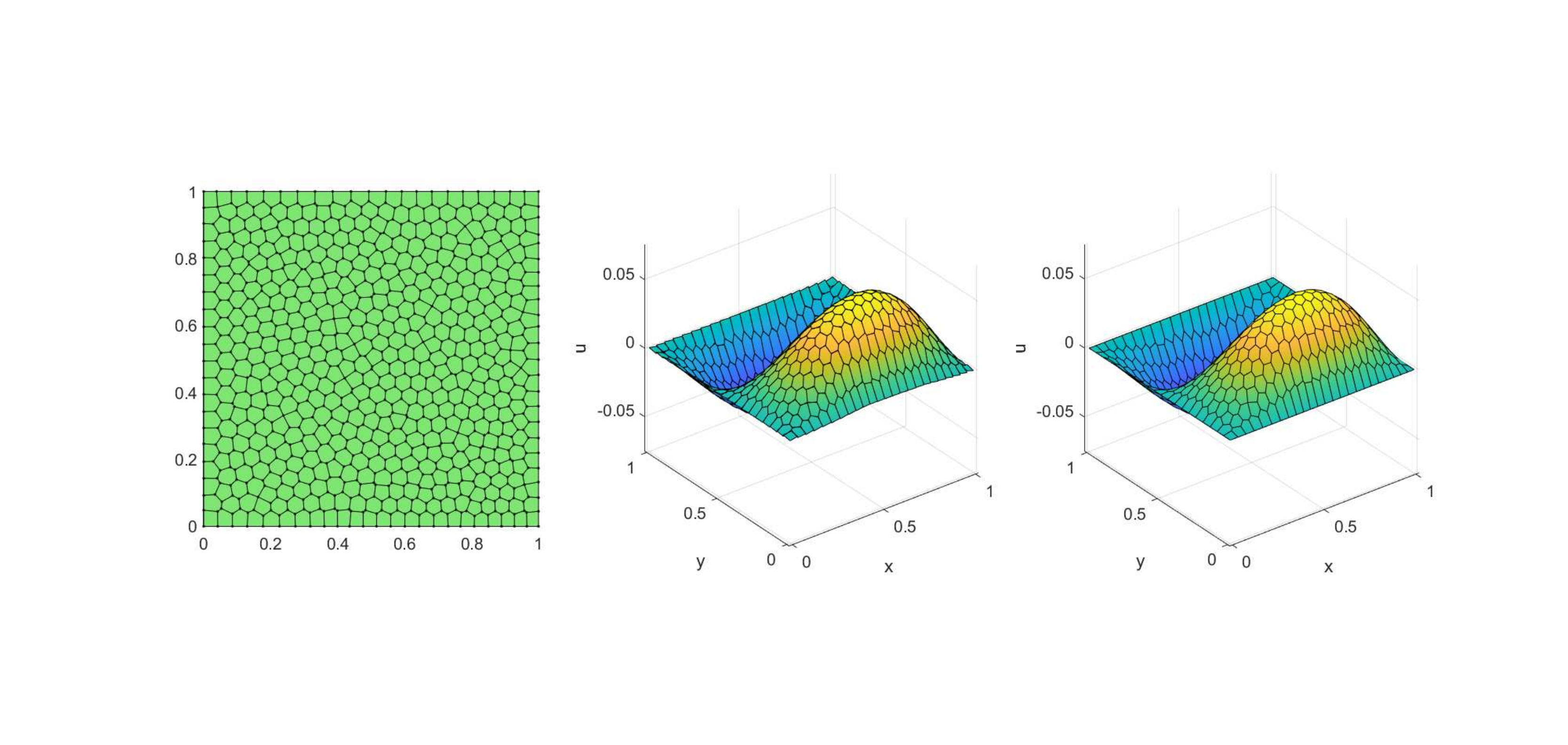}\\
  \caption{Numerical and exact solutions for the simplified friction problem ($\alpha = 10^4$)}\label{fig:VI}
\end{figure}

\begin{figure}[!htb]
  \centering
  \includegraphics[scale=0.45]{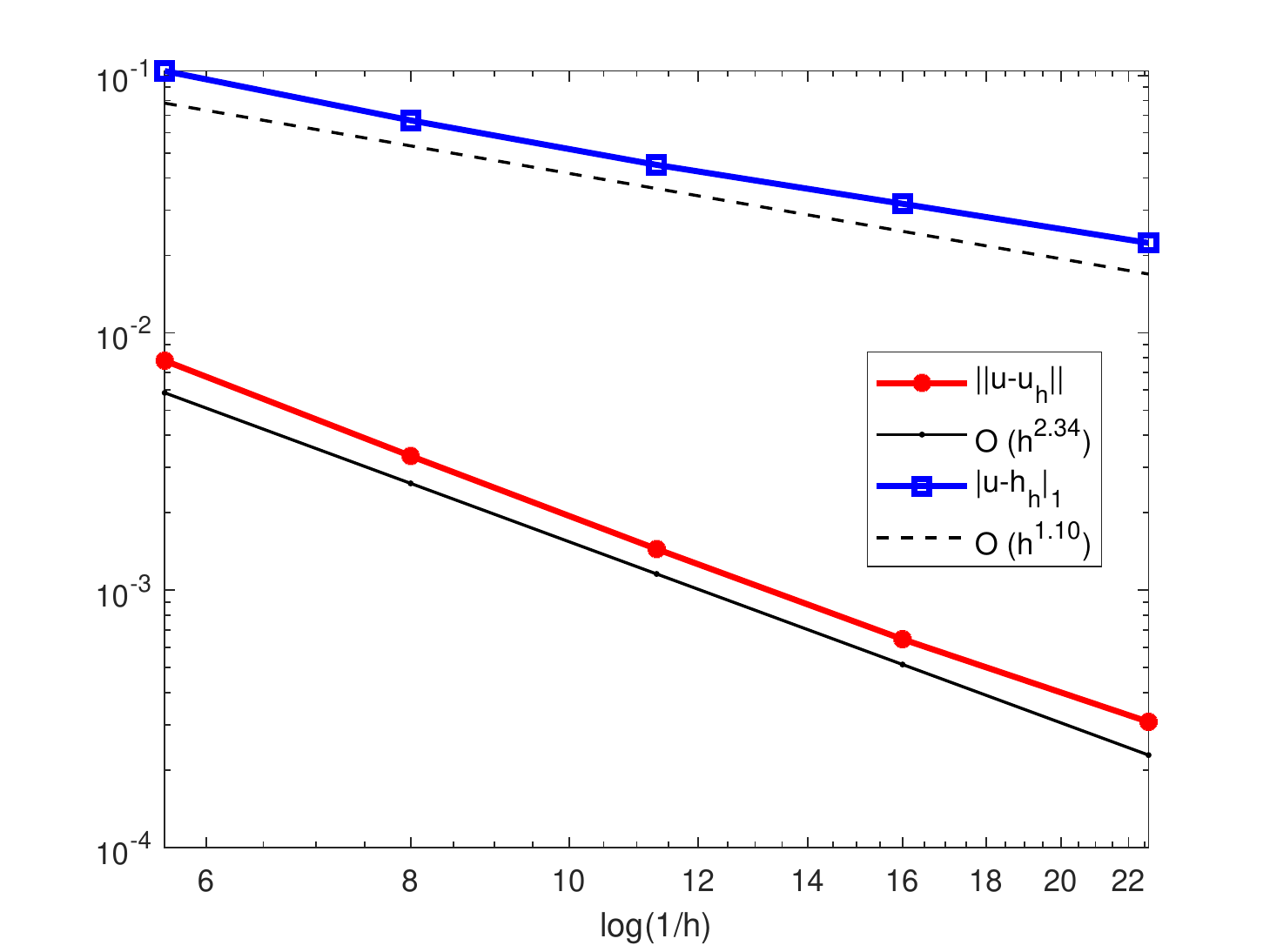}\\
  \caption{Convergence rates of the VEM for the simplified friction problem ($\alpha = 10^4$)}\label{fig:VIrate}
\end{figure}

The Uzawa iteration stops when $\|\bb{\chi}(u_h^{n+1}-u_h^n)\|_{l^2} \le \mbox{tol}$ or $n\ge \mbox{maxIt}$. It is evident that the problem \eqref{Uzawa} is exactly the VEM discretization for the reaction-diffusion problems, with the Neumann boundary data replaced by $g \lambda_h^{(n-1)}$. In addition, we only need to assemble the integral on $\Gamma_C$ in each iteration. Because of this, we will not give the implementation details. We set $\mbox{tol}=10^{-8}$, $\mbox{maxIt}=500$ and $\rho = 10$. The results are shown in Figs.~\ref{fig:VI} and  \ref{fig:VIrate}, from which we see that the lowest-order VEM achieves the linear convergence order in the discrete $H^1$ norm, which is optimal according to the a priori error estimate in \cite{Wang-Wei-2018}. The test script is \mcode{main\_PoissonVEM\_VI\_Uzawa.m}.

\section{Three-dimensional problems} \label{sec:3D}

In this section we are concerned with the implementation of 3-D VEMs proposed in \cite{Beirao-Dassi-Russo-2017} for the Poisson equation in the lowest-order case. Considering the length of the article, we omit the details of the introduction to the function spaces.

\subsection{Virtual element methods for the 3-D Poisson equation}

Let $\Omega \subset \mathbb{R}^3$ be a polyhedral domain and let $\Gamma$ denote a subset of its boundary consisting of some faces. We consider the following model problem
\begin{equation} \label{Poisson3D}
\begin{cases}
- \Delta u  = r \quad & \mbox{in}~~\Omega, \\
u = g_D  \quad & \mbox{on} ~~\Gamma, \\
\partial_n u = g_N \quad & \mbox{on} ~~ \Gamma' = \partial \Omega \backslash \Gamma,
\end{cases}
\end{equation}
where $r\in L^2(\Omega)$ and $g_N\in L^2(\Gamma')$ are the applied load and Neumann boundary data, respectively, and $g_D\in H^{1/2}(\Gamma)$ is the Dirichlet boundary data function.

In what follows, we use $K$ to represent the generic polyhedral element with $f \subset \partial K$ being its generic face. The vertices of a face $f$ are in a counterclockwise order when viewed from the inside.
The virtual element method proposed in \cite{Beirao-Dassi-Russo-2017} for \eqref{Poisson3D} is to find $u_h \in V_\Gamma^k$ such that
\[a_h(u_h,v_h) = \ell_h(v_h), \quad v_h \in V_0^k,\]
where
\[a_h(u_h,v_h) = \sum\limits_{K\in \mathcal{T}_h} a_h^K(u_h,v_h), \quad
\ell_h(v_h) = \int_\Omega r_h v_h {\rm d}x + \int_{\Gamma'} g_h v_h {\rm d}s.\]
The local bilinear form is split into two parts:
\[a_h^K(v,w) = (\nabla \Pi_1^\nabla v, \nabla \Pi_1^\nabla w)_K  + h_KS^K(v-\Pi_1^\nabla v, w - \Pi_1^\nabla w),\]
where $\Pi_1^\nabla : V^1(K) \to \mathbb{P}_1(K)$ is the elliptic projector and $S^K$ is the stabilization term given as
\[S^K(v,w) = \sum\limits_{i=1}^{N^K} \chi_i(v) \chi_i(w),\]
where $\chi_i(v) = v(p_i)$ and $p_i$ is the $i$-th vertex of $K$ for $i=1,2,\cdots, N_K$. The local linear form of the right-hand side will be approximated as
\[\ell_h^K(v_h) = \int_K r \Pi_1^\nabla v_h {\rm d}x + \sum\limits_{f \subset \Gamma' \cap \partial K}\int_f g_N \Pi_{1,f}^\nabla v_h {\rm d}s,\]
where $\Pi_{1,f}^\nabla: V^1(f) \to \mathbb{P}_1(f)$ is the elliptic projector defined on the face $f$.

For the detailed introduction of the virtual element spaces, please refer to Section 2 in \cite{Beirao-Dassi-Russo-2017}. In this paper, we only consider the lowest-order case, but note that the hidden ideas can be directly generalized to higher order cases.

\subsection{Data structure and test script}

\begin{figure}[!htb]
  \centering
  \subfigure[Polyhedral mesh]{\includegraphics[height=6cm]{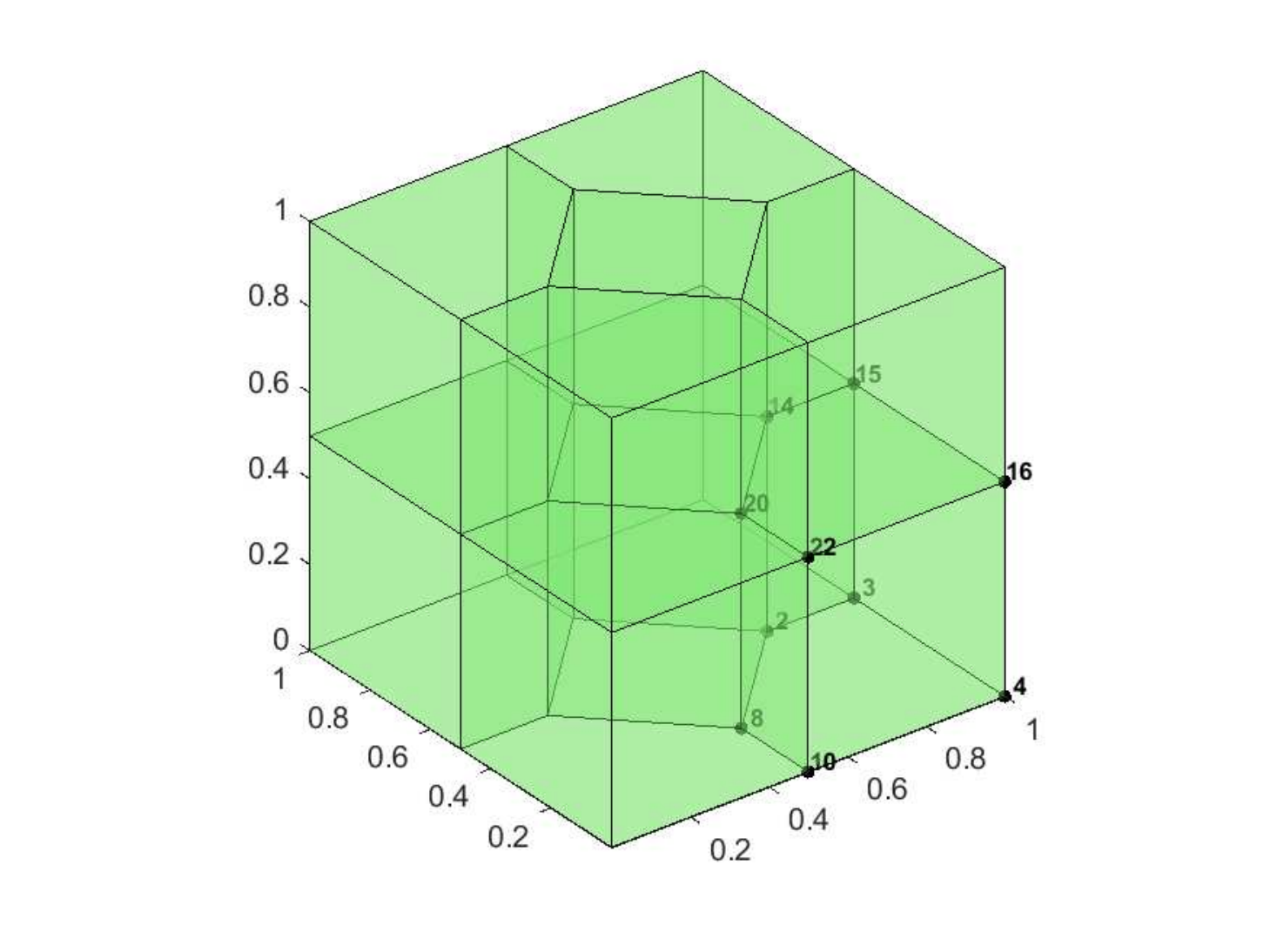}} \hspace{1cm}
  \subfigure[Representation of the first element]{\includegraphics[height=5cm]{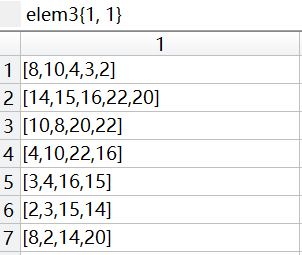}}\\
  \caption{Polyhedral mesh of a cube}\label{fig:mesh3ex1}
\end{figure}

We first discuss the data structure to represent polyhedral meshes. In the implementation, the mesh is represented by \mcode{node3} and \mcode{elem3}. The $\mcode{N} \times 3$ matrix \mcode{node3} stores the coordinates of all vertices in the mesh. \mcode{elem3} is a cell array with each entry storing the face connectivity, for example, the first entry \mcode{elemf = elem3\{1\}} for the mesh given in Fig.~\ref{fig:mesh3ex1}(a) is shown in Fig.~\ref{fig:mesh3ex1}(b), which is still represented by a cell array since the faces may have different numbers of vertices.

All faces including the repeated internal ones can be gathered in a cell array as

\vspace{0.5em}

\hspace{2cm} \mcode{allFace = vertcat(elem3\{:\});  \% cell }

\vspace{0.5em}

\noindent By padding the vacancies and using the \mcode{sort} and \mcode{unique} functions to rows, we obtain the face set \mcode{face}. The cell array \mcode{elem2face} then establishes the map of local index of faces in each polyhedron to its global index in face set \mcode{face}. The above structural information is summarized in the subroutine \mcode{auxstructure3.m}. The geometric quantities such as the diameter \mcode{diameter3}, the barycenter \mcode{centroid3} and the volume \mcode{volume} are computed by \mcode{auxgeometry3.m}. We remark that these two subroutines may be needed to add more information when dealing with higher order VEMs.

The test script is \mcode{main_PoissonVEM3.m} listed as follows. In the \mcode{for} loop, we first load the pre-defined mesh data, which immediately returns the matrix \mcode{node3} and the cell array \mcode{elem3} to the MATLAB workspace. Then we set up the Neumann boundary conditions to get the structural information of the boundary faces. The subroutine \mcode{PoissonVEM3.m} is the function file containing all source code to implement the 3-D VEM. When obtaining the numerical solutions, we calculate
the discrete $L^2$ errors and $H^1$ errors defined as
\[\mbox{ErrL2} = \sum\limits_{K\in \mathcal{T}_h} \| u - \Pi_1^\nabla u_h \|_{0,K}, \quad
  \mbox{ErrH1} = \sum\limits_{K\in \mathcal{T}_h} | u - \Pi_1^\nabla u_h |_{1,K}\]
by using respectively the subroutine \mcode{getError3.m}. The procedure is completed by verifying the rate of convergence through \mcode{showrateErr.m}.
\vspace{-0.8cm}
\begin{lstlisting}
%% Parameters
maxIt = 5;
h = zeros(maxIt,1);      N = zeros(maxIt,1);
ErrL2 = zeros(maxIt,1);  ErrH1 = zeros(maxIt,1);

%% PDE data
pde = Poisson3data();
bdNeumann = 'x==0'; % string for Neumann

%% Virtual element method
for k = 1:maxIt
    % load mesh
    fprintf('Mesh %d: \n', k);
    load( ['SimpleMesh3data', num2str(k), '.mat'] ); % polyhedral mesh
    %load( ['mesh3data', num2str(k), '.mat'] ); % polyhedral mesh
    %[node3,~,elem3] = cubemesh([0 1 0 1 0 1], 1/(2*k)); % tetrahedral mesh
    % get boundary information
    bdStruct = setboundary3(node3,elem3,bdNeumann);
    % solve
    [uh,info] = PoissonVEM3(node3,elem3,pde,bdStruct);
    % record
    N(k) = length(uh);  h(k) = (1/size(elem3,1))^(1/3);
    % compute errors in discrete L2, H1 and energy norms
    kOrder = 1;
    [ErrH1(k),ErrL2(k)] = getError3(node3,elem3,uh,info,pde,kOrder);
end

%% Plot convergence rates and display error table
figure, showrateErr(h,ErrL2,ErrH1);

fprintf('\n');
disp('Table: Error')
colname = {'#Dof','h','||u-u_h||','|u-u_h|_1'};
disptable(colname,N,[],h,'%0.3e',ErrL2,'%0.5e',ErrH1,'%0.5e');
\end{lstlisting}

In the following subsections, we shall go into the details of the implementation of the 3-D VEM in \mcode{PoissonVEM3.m}.

\subsection{Elliptic projection on polygonal faces}

Let $f $ be a face of $K$ or a polygon embedded in $\mathbb{R}^3$. In the VEM computing, we have to get all elliptic projections $\Pi_{1,f}^\nabla \phi_f^T$ ready in advance, where $\phi_f$ is the nodal basis of the enhanced virtual element space $V^1(f)$ (see Subsection 2.3 in \cite{Beirao-Dassi-Russo-2017}).  To this end, it may be necessary to establish local coordinates $(s,t)$ on the face $f$.

\begin{figure}[!htb]
  \centering
  \includegraphics[scale=0.5]{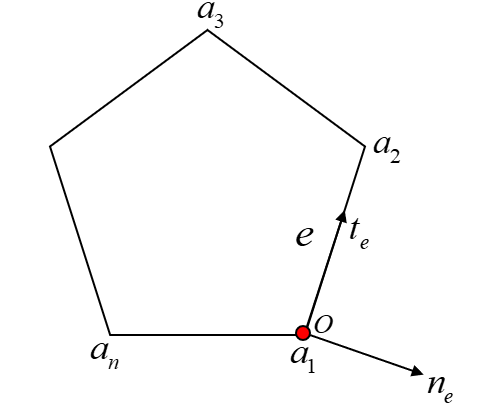}\\
  \caption{Local coordinate system of a face or polygon embedded in $\mathbb{R}^3$}\label{fig:localCord}
\end{figure}

As shown in Fig.~\ref{fig:localCord}, the boundary of the polygon is oriented in a counterclockwise order as $a_1, \cdots, a_n$. Let $e = a_1a_2$ be the first edge, and $n_e$ and $t_e$ be the normal vector and tangential vector, respectively. Then we can define a local coordinate system with $a_1$ being the original point by using these two vectors. Let $n_e = (n_1,n_2,n_3)$ and $t_e = (t_1,t_2, t_3)$. For any $a = (x,y,z) \in f$, its local coordinate $(s,t)$ is related by
\[\overrightarrow{Oa} = s\cdot n_e + t\cdot t_e, \quad \mbox{or} \quad
(x-x_1,y-y_1,z-z_1) = s \cdot (n_1,n_2,n_3) + t \cdot (t_1,t_2,t_3),\]
which gives
\[(s,t) = (x-x_1,y-y_1,z-z_1) \begin{bmatrix} n_1  &  n_2  & n_3 \\ t_1 & t_2  & t_2  \end{bmatrix}^{-1}, \]
with the inverse understood in the least squares sense.
When converting to the local coordinate system, we can compute all the matrices of elliptic projection in the same way for the Poisson equation in two-dimensional cases. For completeness, we briefly recall the implementation.
In what follows, we use the subscript ``$f$'' to indicate the locally defined symbols.

Let $\phi_f^T = [\phi_{f,1}, \cdots, \phi_{f,n}]$ be the basis functions of $V^1(f)$ and $m_f^T =[m_{f,1}, m_{f,2}, m_{f,3}]$ the scaled monomials on $f$ given as
\[m_{f,1} = 1, \quad m_{f,2} = \frac{s - s_f}{h_f}, \quad m_{f,3} = \frac{t - t_f}{h_f},\]
where $(s_f, t_f)$ and $h_f$ are the barycenter and the diameter of $f$, respectively. The vector form of the elliptic projector $\Pi_{1,f}^\nabla$ can be represented as
\begin{equation} \label{ellipticf}
\begin{cases}
(\nabla_f  m_f, \nabla_f \Pi_{1,f}^\nabla \phi_f^T)_f = (\nabla_f  m_f,  \nabla_f \phi_f^T )_f,  \\
P_0(\Pi_{1,f}^\nabla \phi_f^T) = P_0(\phi_f^T).
\end{cases}
\end{equation}
where
\[P_0(v) = \frac{1}{n}\sum\limits_{i=1}^n v (a_i).\]

Since $\mathbb{P}_1(f) \subset V^1(f)$, we can write
\[m_f^T = \phi_f^T \bb{D}_f, \qquad \bb{D}_f = (D_{i\alpha }), \quad D_{i\alpha } = \chi_{f,i} (m_{f,\alpha}),\]
where $\chi_{f,i}$ is the $i$-th d.o.f associated with $a_i$, and $\bb{D}_f$ is referred to as the transition matrix. We further introduce the following expansions
\[\Pi_{1,f}^\nabla{\phi_f^T} = \phi_f^T \bb{\Pi}_{1,f}^\nabla, \qquad
\Pi_{1,f}^\nabla{\phi_f^T} = m_f^T \bb{\Pi}_{1*,f}^\nabla.\]
One easily finds that
\[\bb{\Pi}_{1,f}^\nabla = \bb{D}_f\bb{\Pi}_{1*,f}^\nabla,\]
and \eqref{ellipticf} can be rewritten in matrix form as
\begin{align*}
\begin{cases}
\bb{G}_f \bb{\Pi}_{1*,f}^\nabla = \bb{B}_f,  \\
P_0(m_f^T) \bb{\Pi}_{1*,f}^\nabla = P_0(\phi_f^T)
\end{cases}, \quad \mbox{or denoted by} \quad
\tilde{\bb{G}}_f \bb{\Pi}_{1*,f}^\nabla = \tilde{\bb{B}}_f,
\end{align*}
where
\[\bb{G}_f = (\nabla_f  m_f, \nabla_f m_f^T)_f, \quad
\bb{B}_f = (\nabla_f  m_f, \nabla_f \phi_f^T)_f.\]
Note that the following consistency relation holds
\[\bb{G}_f = \bb{B}_f \bb{D}_f, \quad \tilde{\bb{G}}_f = \tilde{\bb{B}}_f \bb{D}_f.\]

Let \mcode{face} be the face set with internal faces repeated once. Then using the local coordinates we are able to derive all elliptic projections $\Pi_{1,f}^\nabla \phi_f^T$ as in 2-D cases. It is not recommended to carry out the calculation element by element in view of the repeated cost for the internal faces.

The above discussion is summarized in a subroutine with input and output as

\vspace{0.5em}

\hspace{2cm} \mcode{Pifs = faceEllipticProjection(P)},

\vspace{0.5em}

\noindent where \mcode{P} is the coordinates of the face $f$ and \mcode{Pifs} is the matrix representation $\bb{\Pi}_{1*,f}^\nabla$ of $\Pi_{1,f}^\nabla \phi_f^T$ in the basis $m_f^T$. One can derive all matrices by looping over the face set \mcode{face}:
\vspace{-0.8cm}
\begin{lstlisting}
%% Derive elliptic projections of all faces
faceProj = cell(NF,1);
for s = 1:NF
    % Pifs
    faces = face{s};  P = node3(faces,:);
    Pifs = faceEllipticProjection(P);
    % sort the columns
    [~,idx] = sort(faces);
    faceProj{s} = Pifs(:,idx);
end
\end{lstlisting}
Note that in the last step we sort the columns of $\bb{\Pi}_{1*,f}^\nabla$ in ascending order according to the numbers of the vertices. In this way we can easily find the correct correspondence on each element (see Lines 32-35 in the code of the next subsection).

The face integral is then given by
\begin{equation}\label{integralProj}
\int_f \Pi_{1,f}^\nabla \phi_f^T \mathrm{d}\sigma =  \int_f m_f^T \mathrm{d}\sigma \bb{\Pi}_{1*,f}^\nabla
=  ( |f|, 0, 0)\bb{\Pi}_{1*,f}^\nabla,
\end{equation}
where $|f|$ is the area of $f$ and the definition of the barycenter is used.

\subsection{Elliptic projection on polyhedral elements}

The 3-D scaled monomials $m^T = [m_1,m_2,m_3,m_4]$ are
\[m_1 = 1, \quad  m_2 = \frac{x-x_K}{h_K} , \quad  m_3 = \frac{y-y_K}{h_K} , \quad  m_4 = \frac{z-z_K}{h_K} ,\]
where $(x_K,y_K,z_K)$ is the centroid of $K$ and $h_K$ is the diameter, and the geometric quantities are computed by the subroutine \mcode{auxgeometry3.m}.  Similar to the 2-D case, we have the symbols $\bb{D}, \bb{G}, \tilde{\bb{G}}, \bb{B}$ and $ \tilde{\bb{B}}$. For example, the transition matrix is given by
\[\bb{D} = (D_{i \alpha}), \quad  D_{i \alpha} = \chi_i( m_\alpha) = m_\alpha (p_i).\]
The most involved step is to compute the matrix
\begin{align*}
\bb{B}
& = \int_K {\nabla m}  \cdot \nabla {\phi^T}\mathrm{d}x =  - \int_K {\Delta m}  \cdot {\phi^T}\mathrm{d}x + \sum\limits_{f \subset \partial K} {\int_f {(\nabla m \cdot {{\boldsymbol n}_f}){\phi^T}} } \mathrm{d}\sigma \\
& = \sum\limits_{f \subset \partial K} \int_f (\nabla m \cdot \bb{n}_f) \phi^T \mathrm{d}\sigma,
\end{align*}
where $\phi^T = [\phi_1,\phi_2,\cdots, \phi_{N_K}]$ are the basis functions with $\phi_i$ associated with the vertex $p_i$ of $K$.
According to the definition of $V^1(f)$, one has
\begin{align*}
\int_f (\nabla m \cdot \bb{n}_f) \phi^T \mathrm{d}\sigma
& = (\nabla m \cdot \bb{n}_f) \int_f \phi^T \mathrm{d}\sigma
  = (\nabla m \cdot \bb{n}_f) \int_f \Pi_{1,f}^\nabla \phi^T \mathrm{d}\sigma,
\end{align*}
and the last term is available from \eqref{integralProj}. Obviously, for the vertex $p_i$ away from the face $f$ there holds $\Pi_{1,f}^\nabla \phi_i = 0$. In the following code, \mcode{indexEdge} gives the row index in the face set \mcode{face} for each face of \mcode{elemf}, and \mcode{iel} is the index for looping over the elements.
\vspace{-0.8cm}
\begin{lstlisting}
    % ------- element information --------
    % faces
    elemf = elem3{iel};   indexFace = elem2face{iel};
    % global index of vertices and local index of elemf
    [~,index3,~,elemfLocal] = faceTriangulation(elemf);
    % centroid and diameter
    Nv = length(index3);  Ndof = Nv;
    V = node3(index3,:);
    xK = centroid3(iel,1); yK = centroid3(iel,2); zK = centroid3(iel,3);
    hK = diameter3(iel);
    x = V(:,1);  y = V(:,2);  z = V(:,3);

    % ------- scaled monomials ----------
    m1 = @(x,y,z) 1+0*x;
    m2 = @(x,y,z) (x-xK)/hK;
    m3 = @(x,y,z) (y-yK)/hK;
    m4 = @(x,y,z) (z-zK)/hK;
    m = @(x,y,z) [m1(x,y,z),m2(x,y,z),m3(x,y,z),m4(x,y,z)]; % m1,m2,m3,m4
    mc = {m1,m2,m3,m4};
    gradmMat = [0 0 0; 1/hK 0 0; 0 1/hK 0; 0 0 1/hK];

    % -------- transition matrix ----------
    D = m(x,y,z);

    % ----------- elliptic projection -------------
    B = zeros(4,Ndof);
    for s = 1:size(elemf,1)
        % --- information of current face
        % vertices of face
        faces = elemf{s};  P = node3(faces,:);
        % elliptic projection on the face
        idFace = indexFace(s);
        Pifs = faceProj{idFace}; % the order may be not correct
        [~,~,idx] = unique(faces);
        Pifs = Pifs(:,idx);
        % normal vector
        e1 = P(2,:)-P(1,:);  en = P(1,:)-P(end,:);
        nf = cross(e1,en); nf = nf./norm(nf);
        % area
        areaf = polyarea3(P);
        % --- integral of Pifs
        intFace = [areaf,0,0]*Pifs;  % local
        intProj = zeros(1,Ndof);  % global adjustment
        faceLocal = elemfLocal{s};
        intProj(faceLocal) = intFace;
        % add grad(m)*nf
        Bf = dot(gradmMat, repmat(nf,4,1), 2)*intProj;
        B = B + Bf;
    end
    % constraint
    Bs = B;  Bs(1,:) = 1/Ndof;
    % consistency relation
    G = B*D;  Gs = Bs*D;
\end{lstlisting}

\subsection{Computation of the right hand side and assembly of the linear system}

The right-hand side is approximated as
\[F_K = \int_K f \Pi_1^\nabla \phi \mathrm{d} x = (\bb{\Pi}_{1*}^\nabla)^T \int_K f m \mathrm{d}x,\]
where $\Pi_1^\nabla$ is the elliptic projector on the element $K$ and $\bb{\Pi}_{1*}^\nabla$ is the matrix representation in the basis $m^T$. The integral $\int_K f m \mathrm{d}x$ can be approximated by
\[\int_K f m \mathrm{d}x = |K| f(\bb{x}_K) m(\bb{x}_K) = |K| f(\bb{x}_K) [1, 0, 0, 0]^T, \quad \bb{x}_K = (x_K, y_K, z_K). \]
One can also divide the element $K$ as a union of some tetrahedrons and compute the integral using the Gaussian rule. Please refer to the subroutine \mcode{integralPolyhedron.m} for illustration.

One easily finds that the stiffness matrix for the bilinear form is
\[\bb{A}_K = (\bb{\Pi}_{1*}^\nabla )^T \bb{G}\bb{\Pi}_{1*}^\nabla + h_K (\bb{I} - \bb{\Pi}_1^\nabla )^T(\bb{I} - \bb{\Pi}_1^\nabla ).\]
We compute the elliptic projections in the previous section and provide the assembly index element by element. Then the linear system can be assembled using the MATLAB function \mcode{sparse} as follows.
\vspace{-0.8cm}
\begin{lstlisting}
for iel = 1:NT

    ...

    % --------- local stiffness matrix ---------
    Pis = Gs\Bs;   Pi  = D*Pis;   I = eye(size(Pi));
    AK  = Pis'*G*Pis + hK*(I-Pi)'*(I-Pi);
    AK = reshape(AK,1,[]);  % straighten

    % --------- load vector -----------
    %fK = Pis'*[pde.f(centroid3(iel,:))*volume(iel);0;0;0];
    fun = @(x,y,z) repmat(pde.f([x,y,z]),1,4).*m(x,y,z);
    fK = integralPolyhedron(fun,3,node3,elemf);
    fK = Pis'*fK(:);

    % --------- assembly index for ellptic projection -----------
    indexDof = index3;
    ii(ia+1:ia+Ndof^2) = reshape(repmat(indexDof, Ndof, 1), [], 1);
    jj(ia+1:ia+Ndof^2) = repmat(indexDof(:), Ndof, 1);
    ss(ia+1:ia+Ndof^2) = AK(:);
    ia = ia + Ndof^2;

    % --------- assembly index for right hand side -----------
    elemb(ib+1:ib+Ndof) = indexDof(:);
    Fb(ib+1:ib+Ndof) = fK(:);
    ib = ib + Ndof;

    % --------- matrix for L2 and H1 error evaluation  ---------
    Ph{iel} = Pis;
    elem2dof{iel} = indexDof;
end
kk = sparse(ii,jj,ss,N,N);
ff = accumarray(elemb,Fb,[N 1]);
\end{lstlisting}

Note that we have stored the matrix representation $\bb{\Pi}_{1*}^\nabla$ and the assembly index \mcode{elem2dof} in the M-file so as to compute the discrete $L^2$ and $H^1$ errors.

\subsection{Applying the boundary conditions}

We first consider the Neumann boundary conditions. Let $f$ be a boundary face with $n$ vertices. The local load vector is
\[F_f = \int_f g_N \Pi_{1,f}^\nabla \phi_f \mathrm{d} \sigma
 = (\bb{\Pi}_{1*,f}^\nabla)^T \int_f g_N  m_f \mathrm{d} \sigma ,\]
where $g_N = \partial_{\boldsymbol n_f}u = \nabla u \cdot \bb{n}_f$. For simplicity, we provide the gradient $g_N = \nabla u$ in the PDE data instead and compute the true $g_N$ in the M-file. Note that the above integral can be transformed to a 2-D problem by using the local coordinate system as done in the following code, where \mcode{localPolygon3.m} realizes the transformation and returns some useful information, and \mcode{integralPolygon.m} calculates the integral on a 2-D polygon.
\vspace{-0.8cm}
\begin{lstlisting}
%% Assemble Neumann boundary conditions
bdFaceN = bdStruct.bdFaceN;  bdFaceIdxN = bdStruct.bdFaceIdxN;
if ~isempty(bdFaceN)
    Du = pde.Du;
    faceLen = cellfun('length',bdFaceN);
    nnz = sum(faceLen);
    elemb = zeros(nnz,1); FN = zeros(nnz,1);
    ib = 0;
    for s = 1:size(bdFaceN,1)
        % vertices of face
        faces = bdFaceN{s};   nv = length(faces);
        P = node3(faces,:);
        % elliptic projection on the face
        idFace = bdFaceIdxN(s);
        Pifs = faceProj{idFace}; % the order may be not correct
        [~,idx] = sort(faces);
        Pifs = Pifs(:,idx);
        % 3-D polygon -> 2-D polygon
        poly = localPolygon3(P);
        nodef = poly.nodef;  % local coordinates
        nf = poly.nf;   % outer normal vector of face
        centroidf = poly.centroidf;
        sc = centroidf(1); tc = centroidf(2);
        hf = poly.diameterf;
        Coord = poly.Coord; % (s,t) ---> (x,y,z)
        % g_N
        g_N = @(s,t) dot(Du(Coord(s,t)), nf);
        fun = @(s,t) g_N(s,t)*[1+0*s, (s-sc)/hf, (t-tc)/hf];
        Ff = integralPolygon(fun,3,nodef);
        Ff = Pifs'*Ff(:);
        % assembly index
        elemb(ib+1:ib+nv) = faces(:);
        FN(ib+1:ib+nv) = Ff(:);
        ib = ib + nv;
    end
    ff = ff + accumarray(elemb(:), FN(:),[N 1]);
end
\end{lstlisting}

The Dirichlet boundary conditions are easy to handle. The code is given as follows.
\vspace{-0.8cm}
\begin{lstlisting}
%% Apply Dirichlet boundary conditions
g_D = pde.g_D;  bdNodeIdxD = bdStruct.bdNodeIdxD;
isBdNode = false(N,1); isBdNode(bdNodeIdxD) = true;
bdDof = (isBdNode); freeDof = (~isBdNode);
nodeD = node3(bdDof,:);
u = zeros(N,1); uD = g_D(nodeD); u(bdDof) = uD(:);
ff = ff - kk*u;
\end{lstlisting}

In the above codes, \mcode{bdStruct} stores all necessary information of boundary faces. We finally derive the linear system
\mcode{kk*uh = ff}, where \mcode{kk} is the resulting coefficient matrix and \mcode{ff} is the right-hand side.
For small scale linear system, we directly solve it using the backslash command in MATLAB, while for large systems the algebraic multigrid method is used instead.
\vspace{-0.8cm}
\begin{lstlisting}
%% Set solver
solver = 'amg';
if N < 2e3, solver = 'direct'; end
% solve
switch solver
    case 'direct'
        u(freeDof) = kk(freeDof,freeDof)\ff(freeDof);
    case 'amg'
        option.solver = 'CG';
        u(freeDof) = amg(kk(freeDof,freeDof),ff(freeDof),option);
end
\end{lstlisting}
Here, the subroutine \mcode{amg.m} can be found in $i$FEM --- a MATLAB software package for the finite element methods \cite{ChenL-iFEM-2009}.

The complete M-file is \mcode{PoissonVEM3.m}. The overall structure of a virtual element method implementation will be much the same as for a standard finite element method, as outlined in Algorithm \ref{alg:vem3}.
\begin{algorithm}[!htb]
  \caption{An overall structure of the implementation of a 3-D virtual element method \label{alg:vem3}}
\textbf{Input}: Mesh data and PDE data
  \begin{enumerate}
    \item Get auxiliary data of the mesh, including some data structures and geometric quantities;
    \item Derive elliptic projections of all faces;
    \item Compute and assemble the linear system by looping over the elements;
    \item Assemble Neumann boundary conditions and apply Dirichlet boundary conditions;
    \item Set solver and store information for computing errors.
  \end{enumerate}
\textbf{Output}: The numerical DoFs
\end{algorithm}

\subsection{Numerical examples} \label{sec:numerical}

It should be pointed out that all examples in this subsection are implemented in MATLAB R2019b. The domain $\Omega$ is always taken as the unit cube $(0,1)^3$ with the Neumann boundary condition imposed on $x=0$. The implementation can be adapted to the reaction-diffusion equation $-\Delta u + \alpha u = r$ with the details omitted in this article, where $\alpha$ is a nonnegative constant.

\begin{figure}[!htb]
  \centering
  \subfigure[Tetrahedral mesh]{\includegraphics[scale=0.4,trim = 50 20 50 20,clip]{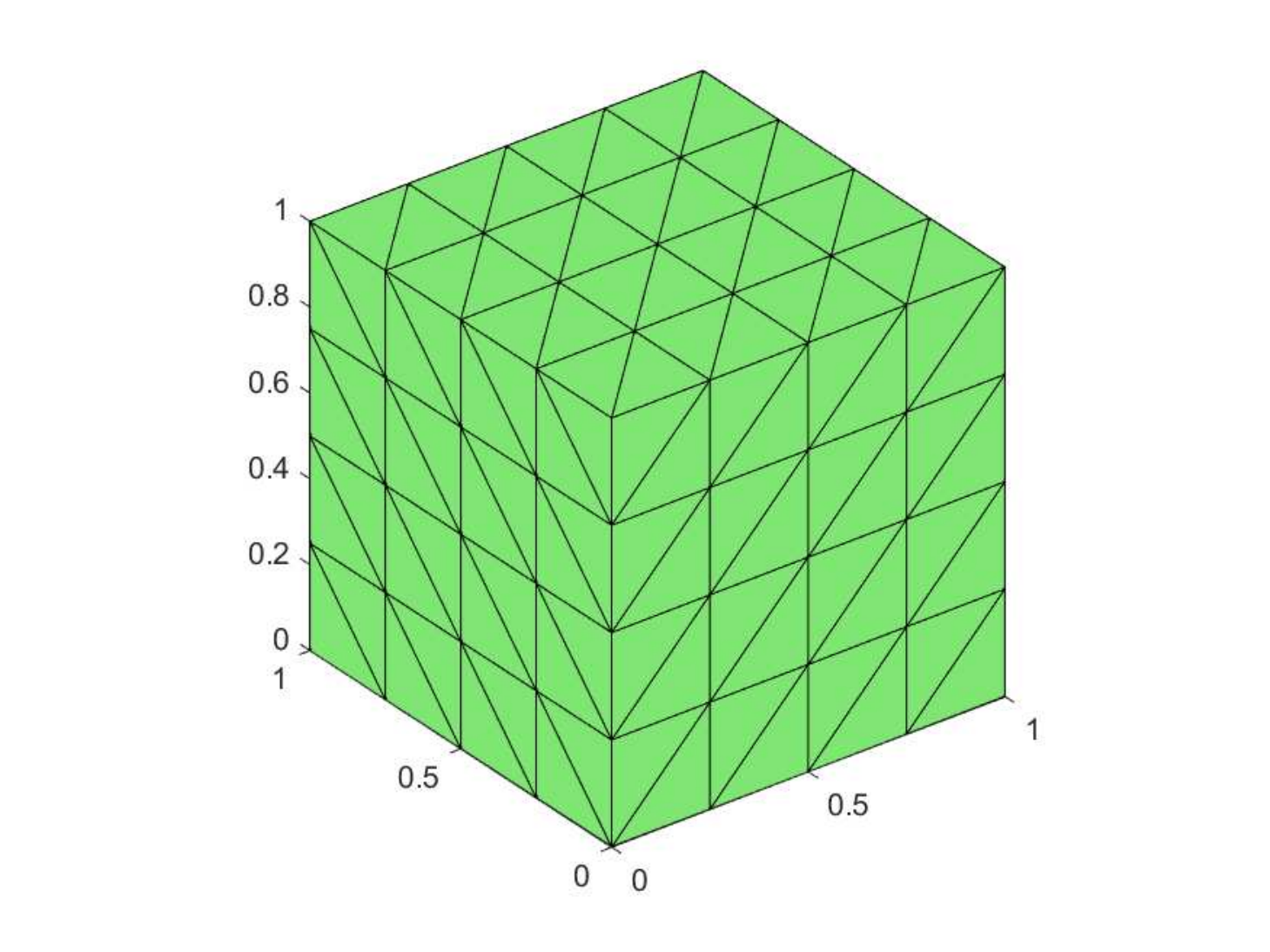}}
  \subfigure[Structured mesh]{\includegraphics[scale=0.4,trim = 50 20 50 20,clip]{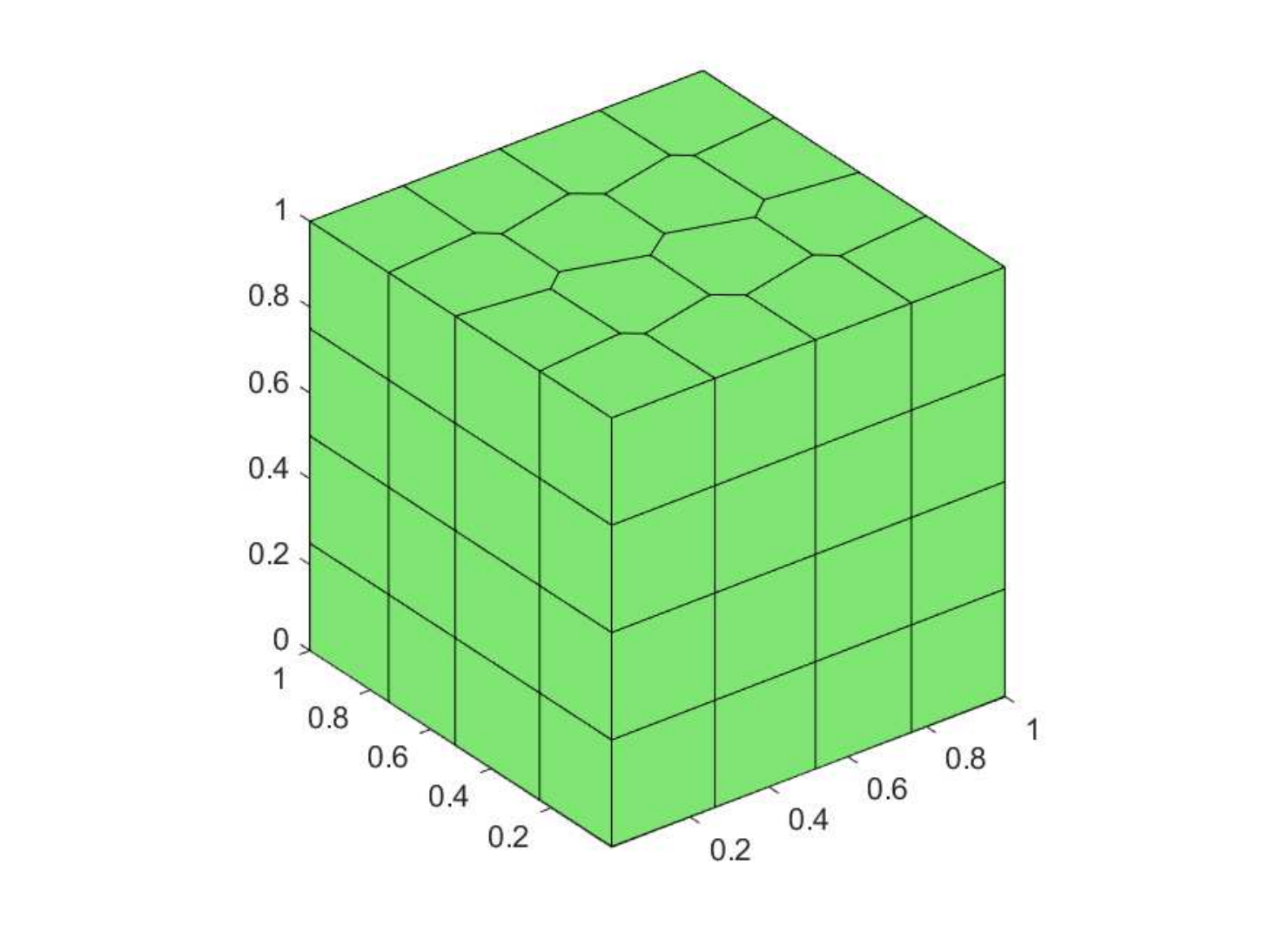}}
  \subfigure[CVT mesh]{\includegraphics[scale=0.4,trim = 50 20 50 20,clip]{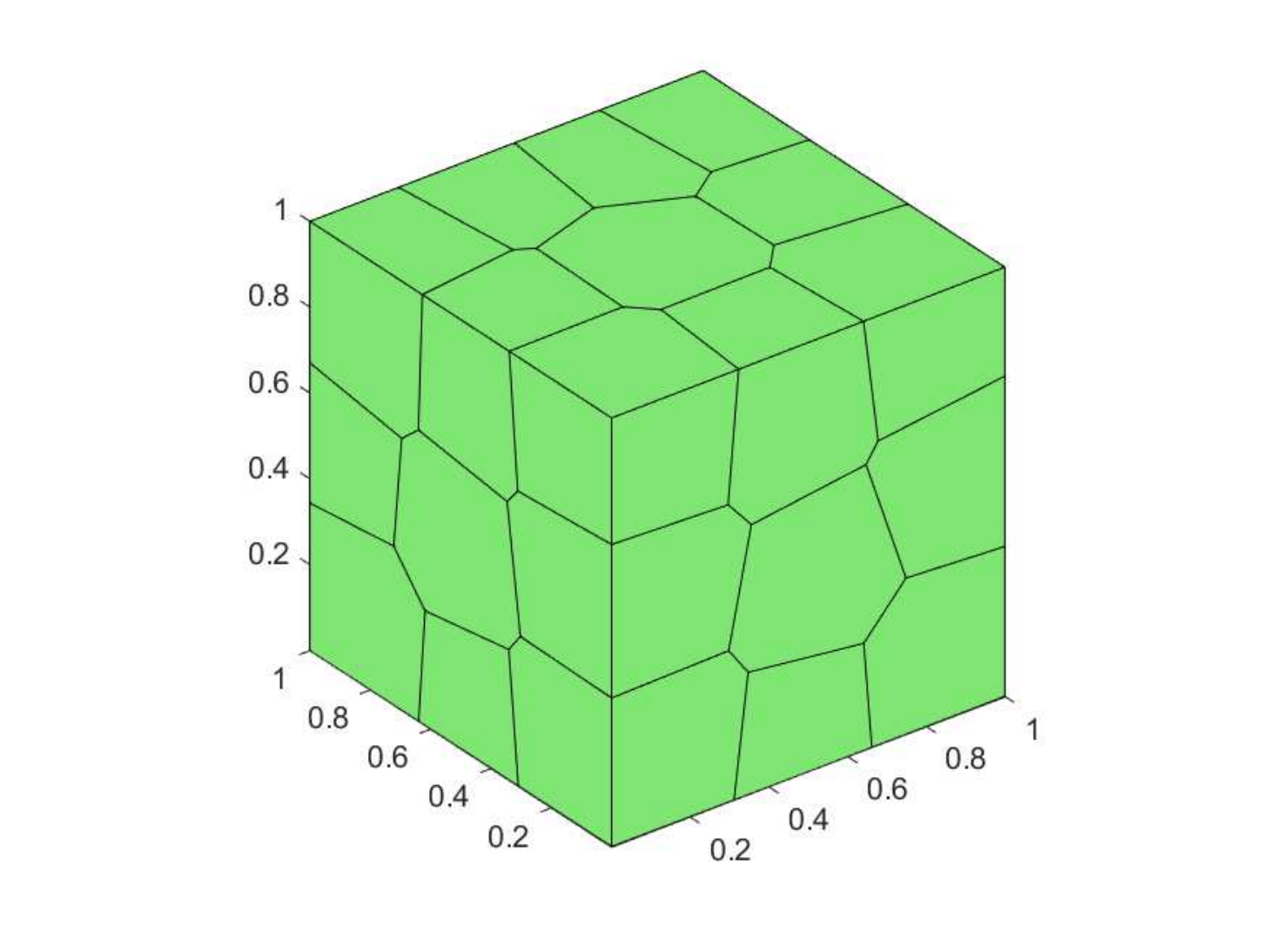}}\\
  \caption{Three types of discretizations}\label{fig:mesh3}
\end{figure}

 We solve the problem on three different kinds of meshes. One is the uniform triangulation shown in Fig.~\ref{fig:mesh3}(a) and the others are the polyhedral meshes displayed in Fig.~\ref{fig:mesh3}(b) and Fig.~\ref{fig:mesh3}(c). The structured polyhedral mesh in Fig.~\ref{fig:mesh3}(b) is formed by translating a two-dimensional polygonal mesh along the $z$-axis and connecting the corresponding vertices, hence all sides are quadrilaterals. The CVT meshes refer to the Centroidal Voronoi Tessellations, which are obtained from a set of seeds that coincide with barycenters of the resulting Voronoi cells. We generate such meshes via a standard Lloyd algorithm by extending the idea in the MATLAB toolbox - PolyMesher introduced in \cite{Talischi-Paulino-Pereira-2012} to a cuboid.

\begin{example}\label{ex1:Poisson3D}
Let $\alpha=1$. The right-hand side $f$ and the boundary conditions are chosen in such a way that the exact solution is
\[u(x,y,z) = \sin(2xy) \cos(z).\]
\end{example}

\begin{figure}[!htb]
  \centering
  \subfigure[Tetrahedral mesh]{\includegraphics[scale=0.35]{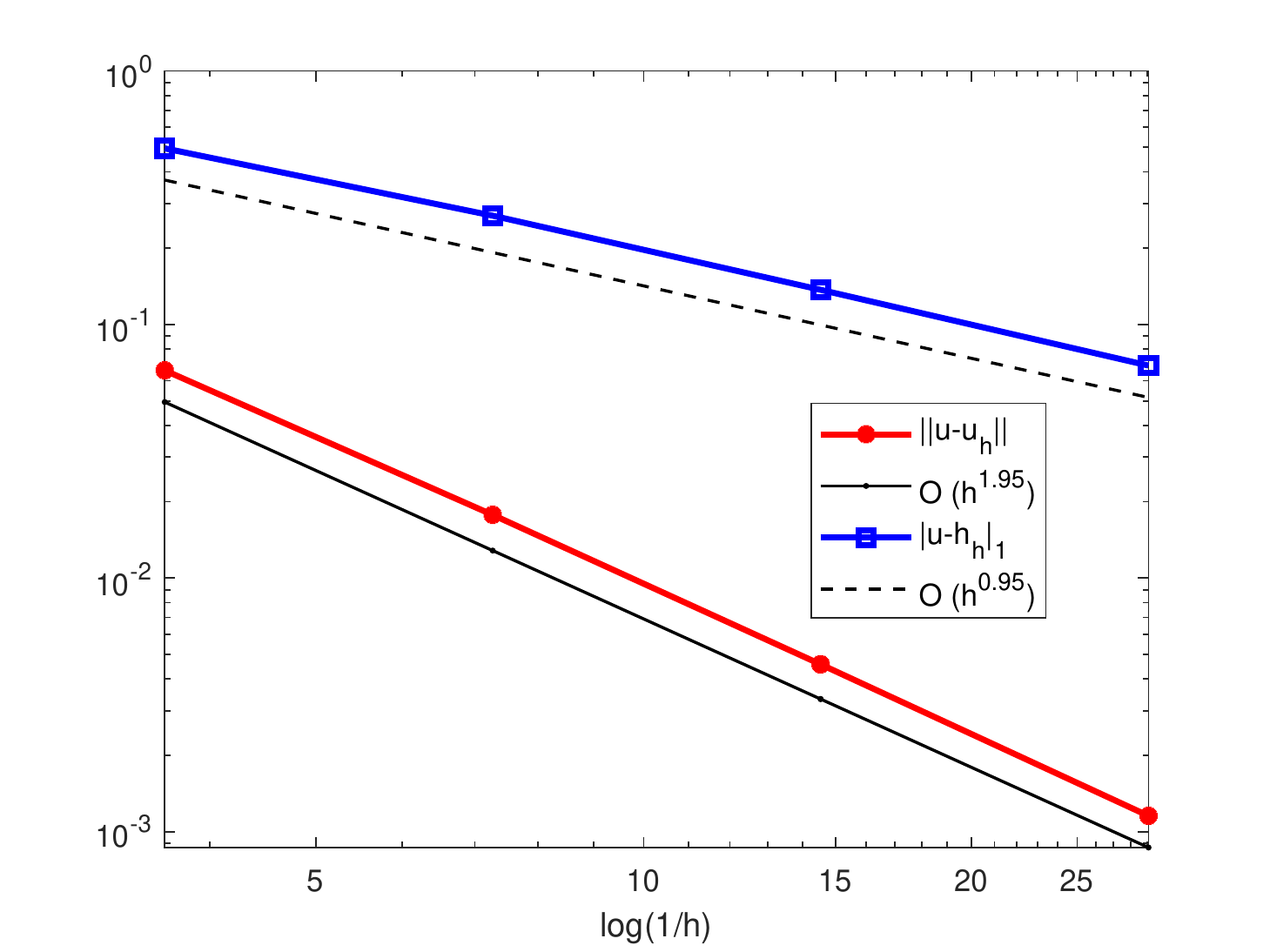}}
  \subfigure[Structured mesh]{\includegraphics[scale=0.35]{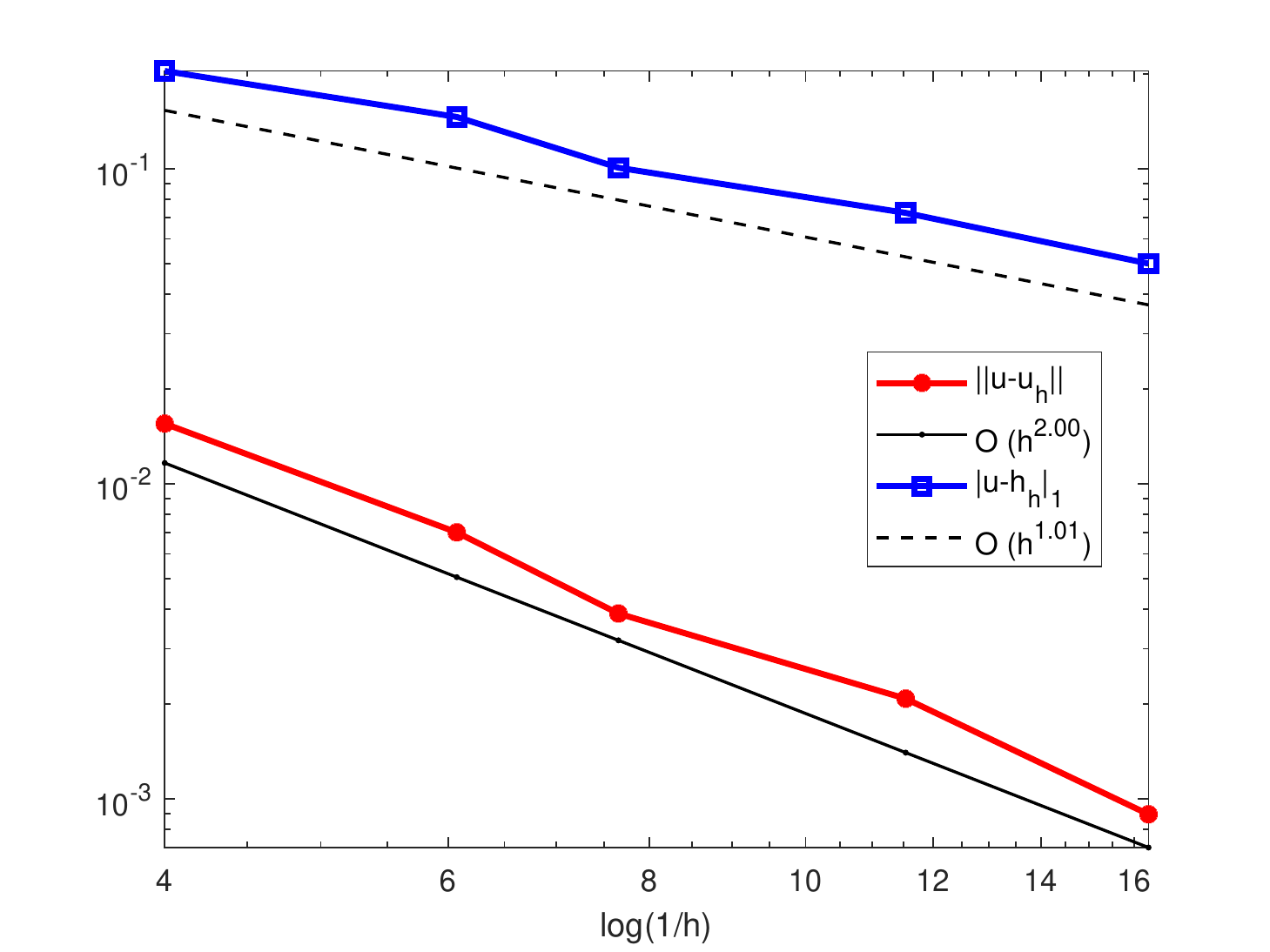}}
  \subfigure[CVT mesh]{\includegraphics[scale=0.35]{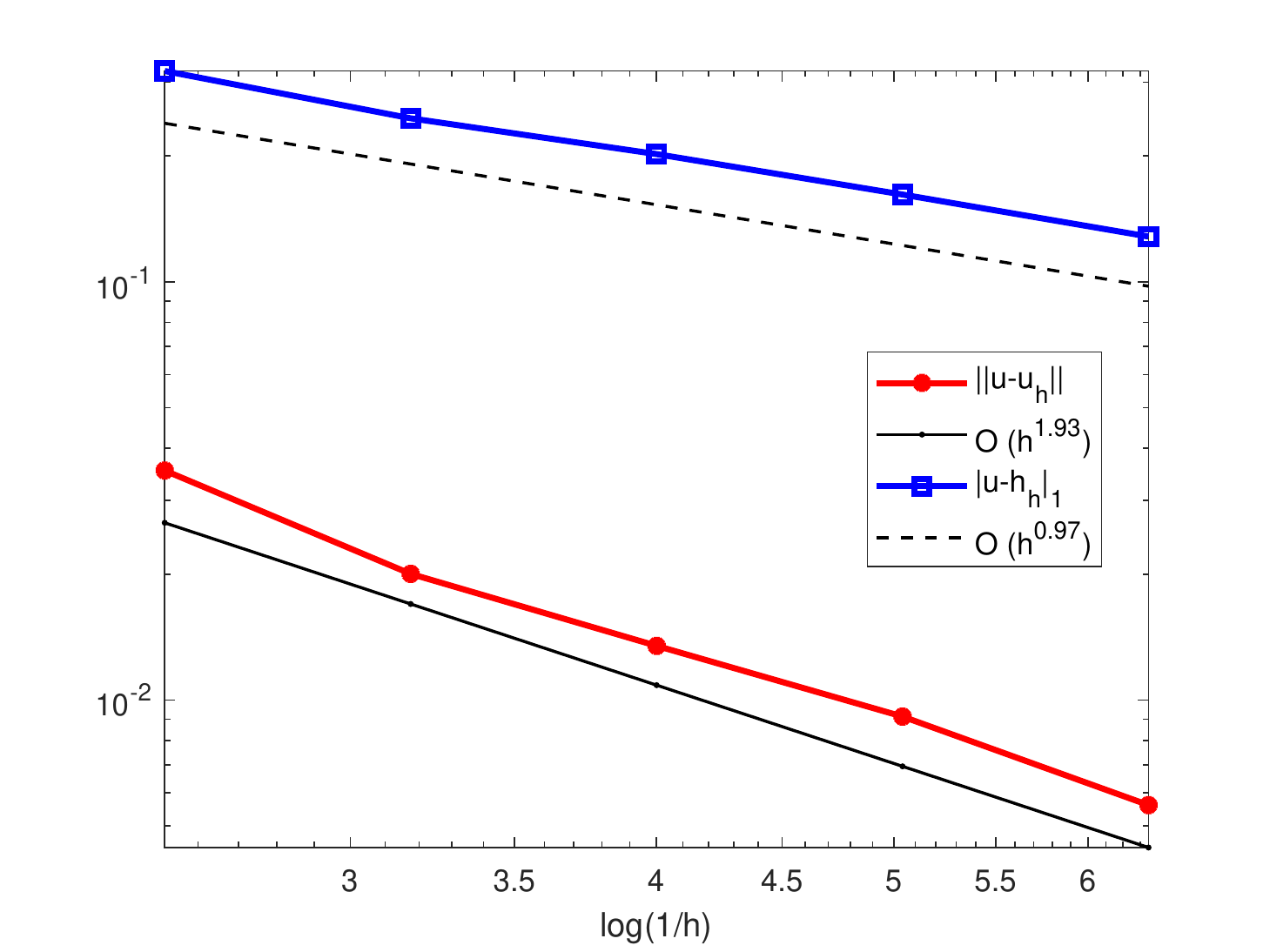}}\\
  \caption{Convergence rates in $L^2$ and $H^1$ norms for Example \ref{ex1:Poisson3D}}\label{fig:PoissonVEM3RatEx1}
\end{figure}

The results are displayed in Fig.~\ref{fig:PoissonVEM3RatEx1}, from which we observe that the optimal rates of convergence are achieved for all the three types of discretizations in the $H^1$ and $L^2$ norms. Note that the uniform triangulation is generated by \mcode{cubemesh.m} in $i$FEM.

\begin{example}\label{ex2:Poisson3D}
In this example, the exact solution is chosen as
\[u(x,y,z) = \sin(\pi x) \cos(\pi y) \cos(\pi z).\]
\end{example}

\begin{figure}[!htb]
  \centering
  \subfigure[Tetrahedral mesh]{\includegraphics[scale=0.45]{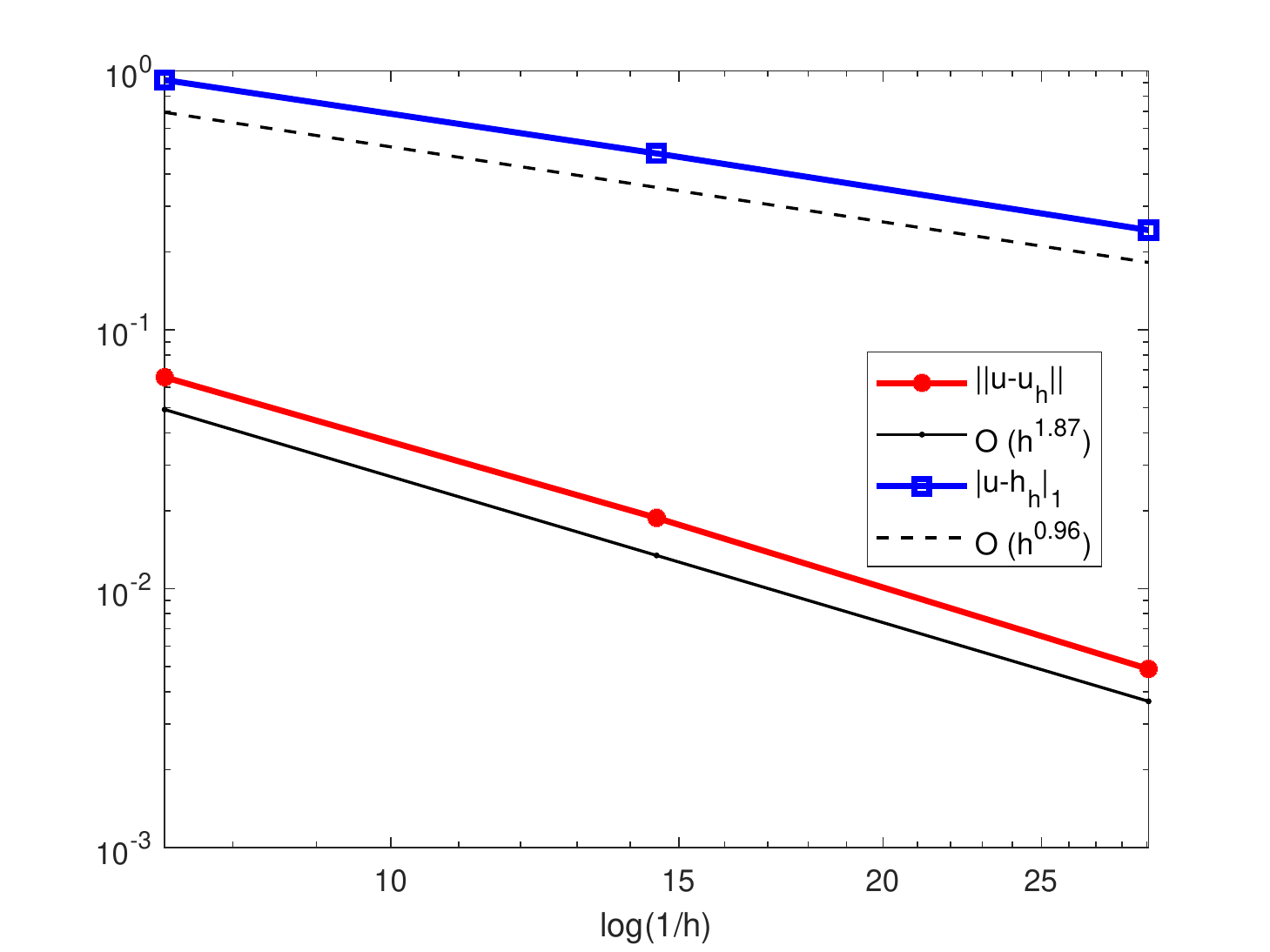}}
  \subfigure[Structured mesh]{\includegraphics[scale=0.45]{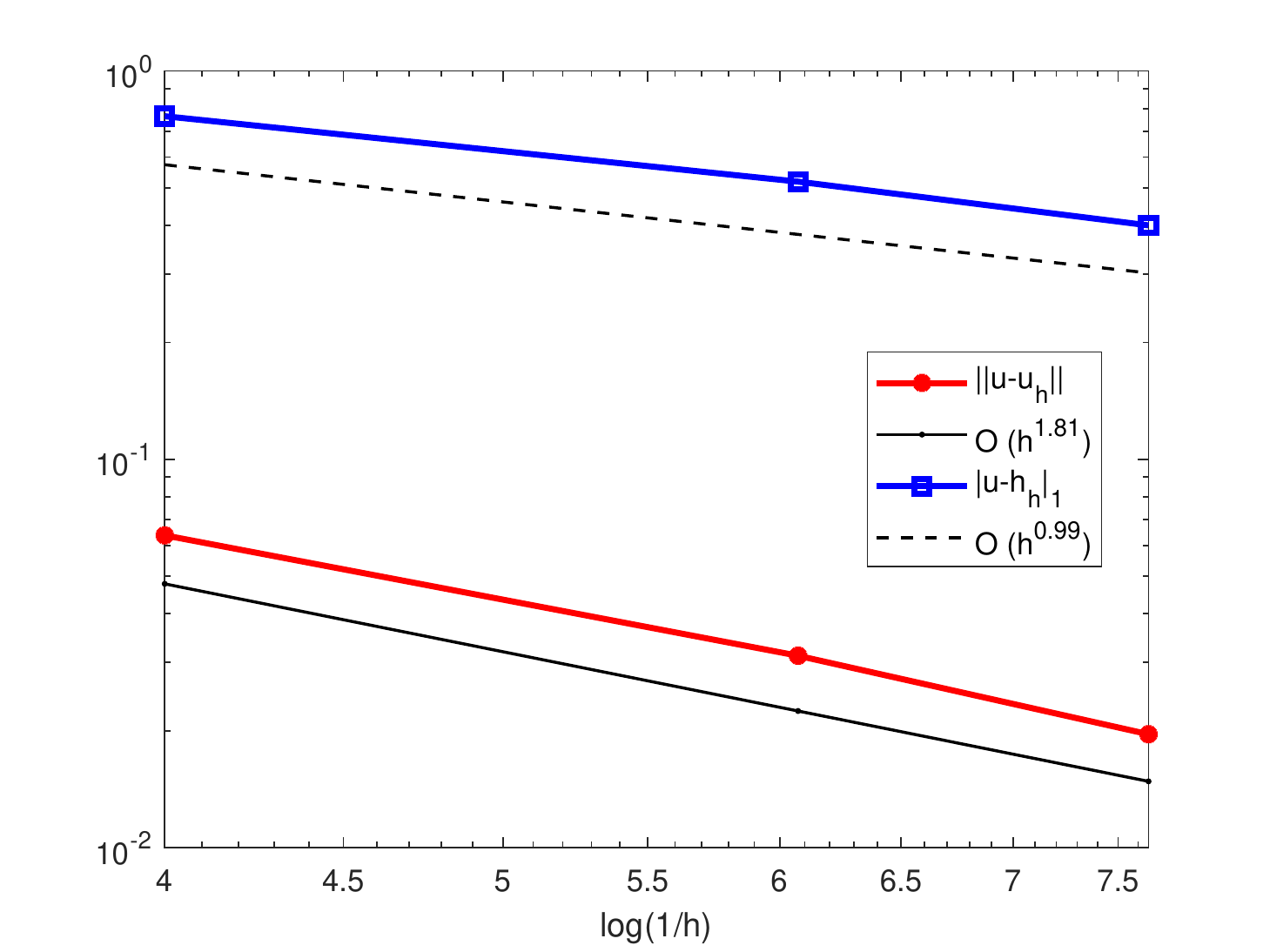}}\\
  \caption{Convergence rates in $L^2$ and $H^1$ norms for Example \ref{ex2:Poisson3D}}\label{fig:PoissonVEM3RatEx2}
\end{figure}

\begin{table}[!htb]
  \centering
  \caption{Discrete $L^2$ and $H^1$ errors of Example \ref{ex2:Poisson3D} for polyhedral meshes}\label{tab:PoissonVEM3ErrEx2}
  \begin{tabular}{@{}cccc@{}}
  \toprule
   $N$            & $h$       &  ErrL2  & ErrH1 \\
  \midrule
   170  & 2.500e-01  &   6.32804e-02   &  7.63490e-01 \\
   504  &  1.647e-01 &   3.09424e-02  &  5.18111e-01   \\
  1024  &  1.307e-01 &   1.94670e-02  &   4.00026e-01  \\
  \bottomrule
\end{tabular}
\end{table}

\begin{figure}[!htb]
  \centering
  \includegraphics[scale=0.45]{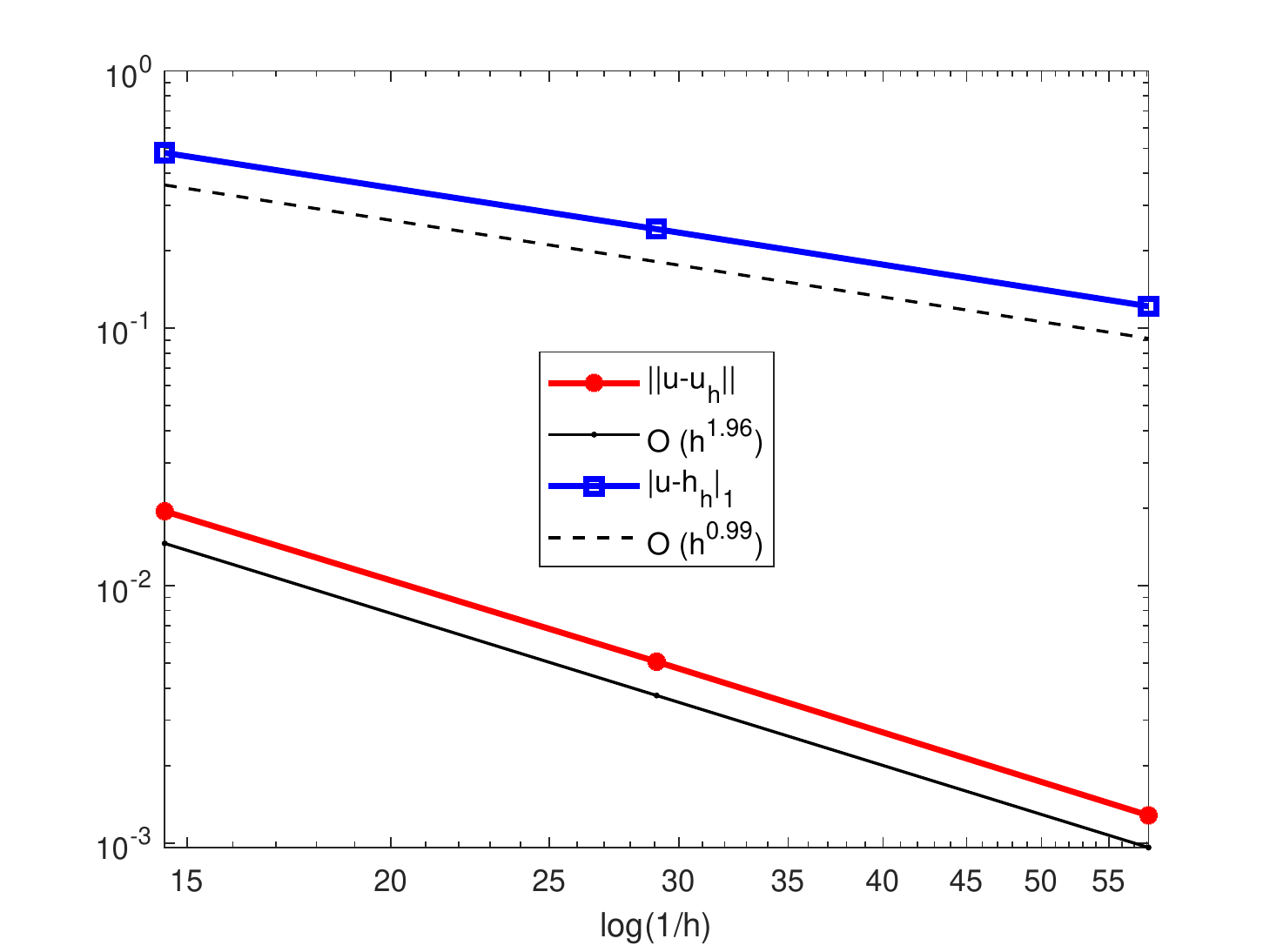}\\
  \caption{Convergence rates in $L^2$ and $H^1$ norms of Example \ref{ex2:Poisson3D} for the triangulation with smaller mesh sizes}\label{fig:PoissonVEM3RatEx2Refine}
\end{figure}

We still consider the problem with the Neumann boundary condition imposed on $x=0$ and repeat the numerical simulation in Example \ref{ex1:Poisson3D}. The uniform triangulation is generated by \mcode{cubemesh.m} and further refined by \mcode{uniformrefine3.m} in $i$FEM. By default, the initial mesh before refining has mesh size $h=1$.
 We display the rate of convergence and the discrete errors in Fig.~\ref{fig:PoissonVEM3RatEx2} and Tab.~\ref{tab:PoissonVEM3ErrEx2}, respectively. It is evident that the code gives satisfactory accuracy and the optimal rate of convergence is achieved for the $H^1$ norm. However, the order of the error in the discrete $L^2$ norm is not close to 2. In fact, this phenomenon is also observed for the classical linear finite element methods under the same conditions. The reason lies in the coarse mesh. To get the optimal convergence rate, one can run the code on a sequence of meshes with much smaller sizes. For instance, Fig.~\ref{fig:PoissonVEM3RatEx2Refine} shows the result for the triangulation with initial size 0.25, in which case the optimal rates of convergence are obtained for both norms. It should be pointed out that the linear virtual element method on a triangulation is exactly the standard finite element method since in this case the virtual element space happens to be the set of polynomials of degree $k\le 1$.
 Compared to the vectorized code in $i$FEM for the finite element methods, the current implementation is less efficient due to the large \mcode{for} loop.
For example, three uniform refinements of the above initial triangulation will yield a mesh of 196608 triangular elements, and hence 196608 loops over the elements.

\section{Concluding remarks} \label{sec:conclude}

In this paper, a MATLAB software package for the virtual element method was presented for various typical problems. The usage
of the library, named mVEM, was demonstrated through several examples.
Possible extensions of this library that are of interest include time-dependent problems, adaptive mixed VEMs, three-dimensional linear elasticity, polyhedral mesh generator, and nonlinear solid mechanics. Various applications such as the Cahn-Hilliard problem, Stokes-Darcy problem and Navier-Stokes are also very appealing.

For three-dimensional problems, the current code can be further vectorized to achieve comparable performance in MATLAB with respect to compiled languages, where the most time-consuming part lies in the evaluation of the large number of face integrals and element integrals, for example, the elementwise computation of the $L^2$ projection matrices for the reaction-diffusion problems. To spare the computational cost, one may divide the polyhedral element as a union of some tetrahedrons and compute the integrals on those elements with the same number of tetrahedrons. Such a procedure can be vectorized in MATLAB with an additional effort of the code design. In the current version, we still do the elementwise loop as in \cite{Sutton-2017} to make the code more transparent.
mVEM is free and open source software.



\end{document}